\documentclass{amsart}
\usepackage{amsmath, amsfonts, amssymb, amsthm, eufrak}
\usepackage[all]{xy}

\let\a\alpha
\let\b\beta

\let\d\delta

\let\k\kappa
\let\l\lambda
\let\m\mu

\let\p\pi
\let\r\rho
\let\s\sigma
\let\t\tau
\let\f\varphi

\let\L\Lambda
\let\La\Lambda

\let\Om\Omega

\def\lra{\longrightarrow}
\def\ra{\rightarrow}
\def\lthra{\twoheadrightarrow}
\def\A{\mathcal A}
\def\B{\mathcal B}
\def\C{\mathcal C}
\def\E{\mathcal E}
\def\F{\mathcal F}
\def\G{\mathcal G}
\def\I{\mathcal I}
\def\J{\mathcal J}
\def\K{\mathcal K}
\def\M{\mathcal M}
\def\P{\mathbb{P}}

\def\T{\mathcal T}

\def\O{\mathcal O}
\def\U{\mathcal U}

\def\Ext{\mathcal Ext}
\def\Hom{\mathcal Hom}

\def\q{\mathfrak{q}}
\def\ag{\mathfrak{a}} 
\def\mg{\mathfrak{m}} 

\def\tensor{\otimes}
\def\isom{\simeq}

\def\bdm{\begin{displaymath}}
\def\edm{\end{displaymath}}

\def\ba{\begin{array}}
\def\ea{\end{array}}

\newcommand{\tilda}{\widetilde}

\begin{document}

\title{on two notions of semistability}
\author{mario maican}
\maketitle

\tableofcontents


\section{Introduction}

The notion of a (Gieseker) semistable sheaf is well-established in
the literature and allows one to construct moduli spaces
of sheaves with fixed Hilbert polynomial on a projective
variety. The construction, carried out in \cite{simpson},
relies on the existence
theorems from Geometric Invariant Theory, more precisely,
it is shown that the moduli space occurs as the quotient
of a certain set of semistable points of a quotient scheme
modulo a reductive algebraic group.

To get a semistable quotient from a semistable sheaf
$\F$ we need to express $\F$ as a quotient
$m \O(-d) \lra \F \lra 0$ with large $m$ and $d$.
In general this procedure is quite abstract and of
little use for the purposes of describing concretely
the geometry of the moduli space.

Another approach for studying moduli spaces uses
monads. Let M$_{\P^2}(r,c_1,c_2)$ be the moduli space
of semistable (in the sense of Mumford-Takemoto) torsion-free
sheaves on $\P^2$ of rank $r$ and Chern classes $c_1$,
$c_2$. Assume that there exist locally free sheaves
$\E_1$, $\E_2$, $\E_3$ on $\P^2$ such that each $\F$
giving a point in M$_{\P^2}(r,c_1,c_2)$ is the cohomology
of a monad
\bdm
0 \lra \E_1 \lra \E_2 \lra \E_3 \lra 0.
\edm
The space $W$ of monads is acted upon in an obvious manner
by the algebraic group
$G=$ Aut$(\E_1)\times$ Aut$(\E_2)\times$ Aut$(\E_3)$.
Two fundamental questions now arise: firstly, is there a
semistability notion for $W$ such that a monad is semistable
precisely if its cohomology is semistable and secondly,
is M$_{\P^2}(r,c_1,c_2)$ a good quotient of the
set $W^{ss}$ of semistable monads modulo $G$ ?

The description of M$_{\P^2}(2,c_1,c_2)$ as a good quotient
was done in \cite{barth} for $c_1$ even and in \cite{hulek}
for $c_1$ odd. In \cite{chang} it was shown that a generic
stable bundle on $\P^3$ of rank 2, Chern classes
$c_1=0$, $c_2=4$ and $\a$-invariant 1 is the cohomology
of a self-dual monad. In \cite{drezet-1987} Dr\'ezet
described as quotients those M$_{\P^2}(r,c_1,c_2)$
for which $\Delta = \d$. He takes $\E_3=0$ and $\E_1$,
$\E_2$ direct sums of certain exceptional bundles.
In all these instances the group $G$ was reductive.
Quotients by nonreductive $G$ were considered by Dr\'ezet
in \cite{drezet-1991} where he studies M$_{\P^2}(r,c_1,c_2)$
of ``faible hauteur''. Again $\E_3=0$, so he is able to
express each semistable bundle as the cokernel of a
semistable morphism.

A notion of semistability for complexes of morphisms of sheaves
modulo nonreductive groups
was proposed by Dr\'ezet and Trautmann in
\cite{drezet-1991}, \cite{drezet-1998} and
\cite{drezet-trautmann}. Let us briefly explain the case of
morphisms of sheaves. In this paper we will not
need the notion of a semistable complex of length
3 or more.
Dr\'ezet and Trautmann consider sheaves $\E_1$
and $\E_2$ on $\P^n$ which are direct sums of simple sheaves,
e.g. direct sums of line bundles.
Thus Aut$(\E_1)\times$ Aut$(\E_2)$ is nonreductive
if $\E_1$ or $\E_2$ has more than one kind of simple sheaf
in its decomposition. This group acts on the vector space
$W=$ Hom$(\E_1,\E_2)$ and the set of semistable points
$W^{ss}$ is defined by means of polarizations which will
be not detailed here. We refer to section 3 for the precise
definition. In \cite{drezet-trautmann} as well as in
\cite{drezet-2000} it was shown that this notion of
semistability quite often leads to a theory similar to
the Geometric Invariant Theory.

Recently, in \cite{fr-trautmann}, Freiermuth and Trautmann
studied the moduli space of semistable (in the sense of
Gieseker) sheaves $\F$ on $\P^3$ with Euler characteristic 1
and with support curves of multiplicity 3.
They show that each $\F$ has a resolution
\bdm
0 \lra 2\O(-3) \stackrel{\psi}{\lra} \O(-1) \oplus 3\O(-2)
\stackrel{\f}{\lra} \O(-1) \oplus \O \lra \F \lra 0
\edm
with $\f$ semistable in the sense of \cite{drezet-trautmann}.
Moreover, the moduli space is a geometric quotient of the
parameter space of $(\psi,\f)$ modulo the action of the group
of automorphisms. \\

In this paper we are interested in semistable sheaves
on $\P^2$ with linear Hilbert polynomial.
Let M$_{\P^2}(r,\chi)$ denote the moduli space of
such sheaves $\F$ with fixed multiplicity $r$ and Euler
characteristic $\chi$. Motivated by \cite{fr-trautmann}
we will seek to express $\F$ as a cokernel
\bdm
\E_1 \stackrel{\f}{\lra} \E_2 \lra \F \lra 0
\edm
with $\E_1$ and $\E_1$ direct sums of line bundles
and $\f$ semistable in the sense of Dr\'ezet and Trautmann.
We carry this out in sections 4, 5, 6 for sheaves
satisfying certain cohomological conditions.
The picture we provide is far from complete
because we do not have a full list of resolutions
for all $\F$ giving a point in M$_{\P^2}(r,\chi)$
even in the case $r=4$ (the cases $r=1,2$ are trivial
while the case $r=3$ is completely understood).

Our cohomological conditions define locally closed
subvarieties in M$_{\P^2}(r,\chi)$ and in section 7
we address the question whether these subvarieties
are good or geometric quotients of the sets of
semistable morphisms $\f$ modulo the canonical
action of the group of automorphisms.
We find that when $r$, $\chi$ are mutually prime,
in other words when M$_{\P^2}(r,\chi)$
is a fine moduli space, we always have geometric
quotients. If the moduli space is not fine the problem
is more complicated and we can answer it only in some
cases.

In section 8 we compute the codimensions of all locally
closed subsets of M$_{\P^2}(r,\chi)$ under investigation.

In section 9 we prove a general duality result.
The dual of a sheaf $\F$ giving a point in
M$_{\P^2}(r,\chi)$ is $\F^D=\Ext^1(\F,\Om^2)(1)$.
Applying the map $\F \lra \F^D$ to a locally closed
subset $X$ in M$_{\P^2}(r,\chi)$ we get a locally
closed subset in M$_{\P^2}(r,r-\chi)$ denoted $X^D$.
At (9.5) we show that under certain conditions
$X$ and $X^D$ are isomorphic. In particular,
this is true for all sets $X$ under investigation
in this paper.
Our theorem is inspired from the result
present in \cite{fr-diplom}, that M$_{\P^2}(r,\chi)$ and
M$_{\P^2}(r,r-\chi)$ are birational if gcd$(r,\chi)=1$.
We show that this is also true for the following
choices of $(r,\chi)$: (6,4), (8,6), (9,6). \\

We summarize our results in the table from below.
The first column contains the cohomological conditions
defining a locally closed subset
$X \subset$ M$_{\P^2}(r,\chi)$. The second column contains
the codimension of $X$. When we write ``0'' we mean an
open dense subset. Each sheaf $\F$ giving a point in
$X$ has resolution of the kind featured in the row
below the semistability conditions.
We have even more detailed information about these
resolutions: the morphisms $\f$ having $\F$ as
cokernel form a subset $W_o$ inside the set
$W^{ss}(G,\L)$ of morphisms which are semistable
with respect to a polarization $\L$ and to the canonical
action of the group $G$ of automorphisms.
We refer to section 3 for the terminology.
The third column of our table contains the information
about $\L$ and the forth column tells us whether
$X$ is a quotient of $W_o$ by $G$. When we write
``good'' it is self-understood that the quotient is not
geometric. We wrote ``unknown'' whenever we could
not prove that a quotient exists.
The subset $W_o \subset W^{ss}(G,\L)$ is given by
the following conditions: for all the blocks different
than the last block in the table we require that
$\f$ be injective and that its scalar entries
(regarding it as a matrix) are zero.
For the last block we refer to (6.10) and (6.11).

\noindent \\ \\
{\small
\begin{tabular}{|l|c|c|l|}
\hline \hline
\multicolumn{4}{|l|}{ \begin{tabular}{l}
M$_{\P^2}(n+1,n) \qquad n \ge 1$
\end{tabular}
} \\
\hline
\begin{tabular}{l}
$h^0(\F(-1))=0$
\end{tabular}
& 0 &  $0 < \l_1 < \frac{1}{n}$  & 
\begin{tabular}{l}
geometric
\end{tabular} \\
\hline
\multicolumn{4}{|c|}{ $ 0 \lra \O(-2) \oplus (n-1)\O(-1) \lra n\O
\lra 0$ } \\
\hline \hline
\multicolumn{4}{|l|}{M$_{\P^2}(n+2,n) \qquad n=3,4,5,6$} \\
\hline
\begin{tabular}{l}
$h^0(\F(-1))=0$ \\
$h^1(\F)=0$ \\
$h^1(\F \tensor \Om^1(1))=0$
\end{tabular}
& 0 & $\frac{1}{2n} < \l_1 < \frac{1}{n}$ &
\begin{tabular}{l}
good for $n=4, 6$ \\
geometric for $n=3,5$
\end{tabular} \\
\hline
\multicolumn{4}{|c|}{$0 \lra 2\O(-2) \oplus (n-2)\O(-1)
\lra n\O \lra \F \lra 0$} \\
\hline
\begin{tabular}{l}
$h^0(\F(-1))=0$ \\
$h^1(\F)=0$ \\
$h^1(\F \tensor \Om^1(1))=1$
\end{tabular}
& $n-1$ &
\begin{tabular}{l}
$(\l_1,\m_1)$ in the interior of \\
the triangle with vertices \\
$(0,0),\ \left( \frac{1}{n+1}, \frac{1}{n+1} \right)$, \\
$\left( \frac{1}{n^2-n+2},\frac{2}{n^2-n+2} \right)$
\end{tabular} &
\begin{tabular}{l}
geometric for $n=3,5$ \\
unknown for $n=4,6$
\end{tabular} \\
\hline
\multicolumn{4}{|c|}{$0 \lra 2\O(-2) \oplus (n-1)\O(-1)
\lra \O(-1) \oplus n\O \lra \F \lra 0$} \\
\hline \hline
\multicolumn{4}{|l|}{M$_{\P^2}(4,2)$} \\
\hline
\begin{tabular}{l}
$h^0(\F(-1))=0$ \\
$h^1(\F)=0$ \\
$h^1(\F \tensor \Om^1(1))=0$
\end{tabular}
& 0 & $\l_1 = \frac{1}{2}$ &
\begin{tabular}{l}
good
\end{tabular} \\
\hline
\multicolumn{4}{|c|}{$0 \lra 2\O(-2) \lra 2\O \lra \F
\lra 0$} \\
\hline
\begin{tabular}{l}
$h^0(\F(-1))=0$ \\
$h^1(\F)=0$ \\
$h^1(\F \tensor \Om^1(1))=1$
\end{tabular}
& 1 &
\begin{tabular}{l}
$(\l_1,\m_1)$ in the interior of \\
the quadrilater with vertices \\
$(0,0),\ \left( \frac{1}{3},\frac{1}{3} \right),\
\left( \frac{1}{2}, 1 \right),\ (0,1)$
\end{tabular} &
\begin{tabular}{l}
unknown
\end{tabular} \\
\hline
\multicolumn{4}{|c|}{$0 \lra 2\O(-2) \oplus \O(-1) \lra
\O(-1) \oplus 2\O \lra \F \lra 0$} \\
\hline 
\end{tabular}

\noindent
\begin{tabular}{|l|c|c|l|}
\hline
\multicolumn{4}{|l|}{M$_{\P^2}(n+3,n) \qquad n=4,5,6$} \\
\hline
\begin{tabular}{l}
$h^0(\F(-1))=0$ \\
$h^1(\F)=0$ \\
$h^1(\F \tensor \Om^1(1))=0$
\end{tabular}
& 0 & $\frac{2}{3n} < \l_1 < \frac{1}{n}$ &
\begin{tabular}{l}
geometric for $n=4,5$ \\
unknown for $n=6$
\end{tabular} \\
\hline
\multicolumn{4}{|c|}{$0 \lra 3\O(-2) \oplus (n-3) \O(-1) \lra
n\O \lra \F \lra 0$} \\
\hline
\begin{tabular}{l}
$h^0(\F(-1))=0$ \\
$h^1(\F)=0$ \\
$h^1(\F \tensor \Om^1(1))=1$
\end{tabular}
& $n-2$ & see (4.7) &
\begin{tabular}{l}
geometric for $n=4,5$ \\
unknown for $n=6$
\end{tabular} \\
\hline
\multicolumn{4}{|c|}{$0 \lra 3\O(-2) \oplus (n-2)\O(-1) \lra
\O(-1) \oplus n\O \lra \F \lra 0$} \\
\hline \hline
\multicolumn{4}{|l|}{M$_{\P^2}(7,4)$} \\
\hline
\begin{tabular}{l}
$h^0(\F(-1))=0$ \\
$h^1(\F)=0$ \\
$h^1(\F \tensor \Om^1(1))=2$
\end{tabular}
& 6 &
\begin{tabular}{l}
$(\l_1,\m_1)$ in the interior of \\
the quadrilater with vertices \\
$(0,0),\ \left( \frac{1}{3},\frac{1}{2} \right),\
\left( \frac{17}{24}, 1 \right),\ (1,1)$
\end{tabular} &
\begin{tabular}{l}
geometric
\end{tabular} \\
\hline
\multicolumn{4}{|c|}{$0 \lra 3\O(-2) \oplus 3\O(-1) \lra
2\O(-1) \oplus 4\O \lra \F \lra 0$} \\
\hline \hline
\multicolumn{4}{|l|}{M$_{\P^2}(6,3)$} \\
\hline
\begin{tabular}{l}
$h^0(\F(-1))=0$ \\
$h^1(\F)=0$ \\
$h^1(\F \tensor \Om^1(1))=0$
\end{tabular}
& 0 & $\l_1 = \frac{1}{3}$ &
\begin{tabular}{l}
good
\end{tabular} \\
\hline
\multicolumn{4}{|c|}{$0 \lra 3 \O(-2) \lra 3\O \lra \F \lra 0$} \\
\hline
\begin{tabular}{l}
$h^0(\F(-1))=0$ \\
$h^1(\F)=0$ \\
$h^1(\F \tensor \Om^1(1))=1$
\end{tabular}
& 1 &
\begin{tabular}{l}
$(\l_1,\m_1)$ in the interior of \\
the segment with endpoints \\
$\left( \frac{1}{4}, \frac{1}{4} \right),\
\left( \frac{1}{5}, \frac{2}{5} \right)$
\end{tabular} &
\begin{tabular}{l}
unknown
\end{tabular} \\
\hline
\multicolumn{4}{|c|}{$0 \lra 3\O(-2) \oplus \O(-1) \lra
\O(-1) \oplus 3\O \lra \F \lra 0$} \\
\hline
\begin{tabular}{l}
$h^0(\F(-1))=0$ \\
$h^1(\F)=0$ \\
$h^1(\F \tensor \Om^1(1))=2$
\end{tabular}
& 4 &
\begin{tabular}{l}
$(\l_1,\m_1)$ in the interior of \\
the triangle with vertices \\
$(0,0),\ \left( \frac{1}{5}, \frac{1}{5} \right),\
\left( \frac{1}{3}, \frac{1}{2} \right)$
\end{tabular} & 
\begin{tabular}{l}
unknown
\end{tabular} \\
\hline
\multicolumn{4}{|c|}{$0 \lra 3\O(-2) \oplus 2\O(-1) \lra
2 \O(-1) \oplus 3\O \lra \F \lra 0$} \\
\hline
\begin{tabular}{l}
$h^0(\F(-1))=0$ \\
$h^1(\F)=1$
\end{tabular}
& 4 & $ 0 < \l_1 < \frac{1}{4}$ & 
\begin{tabular}{l}
geometric
\end{tabular} \\
\hline
\multicolumn{4}{|c|}{$0 \lra \O(-3) \oplus 3\O(-1) \lra 4\O
\lra \F \lra 0$} \\
\hline
\begin{tabular}{l}
$h^0(\F(-1))=1$ \\
$h^1(\F)=0$
\end{tabular}
& 4 & $0 < \m_2 < \frac{1}{4}$ &
\begin{tabular}{l}
geometric
\end{tabular} \\
\hline
\multicolumn{4}{|c|}{$0 \lra 4\O(-2) \lra 3\O(-1) \oplus \O(1)
\lra \F \lra 0$} \\
\hline \hline
\multicolumn{4}{|l|}{M$_{\P^2}(4,1)$} \\
\hline
\begin{tabular}{l}
$h^0(\F(-1))=0$ \\
$h^1(\F)=1$
\end{tabular}
& 2 & $0 < \l_1 < \frac{1}{2}$ &
\begin{tabular}{l}
geometric
\end{tabular} \\
\hline
\multicolumn{4}{|c|}{$0 \lra \O(-3) \oplus \O(-1) \lra
2\O \lra \F \lra 0$} \\
\hline \hline
\multicolumn{4}{|l|}{M$_{\P^2}(5,2)$} \\
\hline
\begin{tabular}{l}
$h^0(\F(-1))=0$ \\
$h^1(\F)=1$
\end{tabular}
& 3 & $0 < \l_1 < \frac{1}{3}$ &
\begin{tabular}{l}
geometric
\end{tabular} \\
\hline
\multicolumn{4}{|c|}{$0 \lra \O(-3) \oplus 2\O(-1) \lra
3\O \lra \F \lra 0$} \\
\hline \hline
\multicolumn{4}{|l|}{M$_{\P^2}(n,3) \qquad 4 \le n \le 15$} \\
\hline
\begin{tabular}{l}
$h^0(\F(-1))=1$ \\
$h^1(\F)=0$
\end{tabular}
& $n-2$ &
\begin{tabular}{l}
$(\l_1,\m_1)$ in the interior of \\
the triangle with vertices \\
$(0,0),\ \left( \frac{1}{n}, \frac{1}{n} \right),\
\left( \frac{1}{n-2}, \frac{1}{n-3} \right)$
\end{tabular} &
\begin{tabular}{l}
geometric for $n \neq 3\k$ \\
unknown for $n= 3\k$
\end{tabular} \\
\hline
\multicolumn{4}{|c|}{$0 \lra \O(-2) \lra (n-2)\O(-2) \oplus
3\O(-1) \lra (n-3)\O(-1) \oplus 3\O \lra \F \lra 0$} \\
\hline \hline
\end{tabular}
}

\noindent \\ 
We should mention that, by virtue of our duality
results (9.5) and (9.7), for each block in the table
there is a ``dual block'' obtained by replacing
M$_{\P^2}(r,\chi)$ with M$_{\P^2}(r,r-\chi)$,
$\F$ with $\F^D$ and $\f$ with $\Hom (\f,\Om^2)(1)$.
We did not feel the need to include another ``dual''
table, instead we only spell out at (9.8) the cases
of open dense subsets.

As a general remark, the kind of arguments from this
paper become very hard to carry out in the
case of large multiplicity. This is so because, when
the multiplicity becomes large, other than
semistability conditions on $\f$ enter into play.
Thus, for large multiplicity, Dr\'ezet and Trautmann's
notion of semistability is no longer satisfactory.

\noindent \\
\emph{Acknowledgements:} The author wishes to thank
J.-M. Dr\'ezet for many useful comments. 
The referee pointed out several improvements, including
a simplification of the proof of (4.3), for which the
author is grateful.


\section{Semistable Sheaves and Their Moduli}

From now on $k$ will be an algebraically closed field of
characteristic zero. All schemes over $k$ will be assumed to
be algebraic, meaning that they can be covered with finitely
many spectra of finitely generated $k$-algebras.
A separated algebraic scheme will also be called
an \emph{algebraic variety}. A \emph{variety} will be the
maximal spectrum of a reduced algebraic variety.
Our main reference for this section is \cite{hl}.

Let $X$ be a smooth projective variety of dimension $n$
with ample line bundle $\O_X(1)$. For a coherent sheaf
$\F$ on $X$ we denote by $\chi(\F)$ its Euler characteristic
given by
\begin{align*}
\chi(\F)=\sum_{i\ge 0} (-1)^{i} \text{dim}_k \, H^i (X,\F).
\end{align*}
The Euler characteristic of the twisted sheaf $\F(m)=\F \tensor
\O_X(m)$ is a polynomial expression in $m$. Thus, we can define
the \emph{Hilbert polynomial} $P_{\F}(m)$ of $\F$ by the formula
\begin{align*}
P_{\F}(m)=\chi(\F(m)).
\end{align*}
It is known that the degree of $P_{\F}(m)$ equals the dimension
of the topological support supp$(\F)$ of $\F$. We write
\begin{align*}
P_{\F}(m)= \sum_{i=0}^d \a_i (\F) \frac{m^i}{i!}.
\end{align*}
The coefficients $\a_i (\F)$ are integers, see \cite{hl}.
The dominant coefficient $\a_d (\F)$ is called the \emph{multiplicity}
of $\F$ and is positive because, by the Theorem B of Serre,
for $m$ large enough we have $P_{\F}(m)=$dim$_k H^0(X,\F(m))>0$.
It is known that $\a_d (\F)$ equals the degree of the scheme
Supp$(\F)$ which has supp$(\F)$ as underlying topological
space and $\O_X/{\mathcal Ann}(\F)$ as structure sheaf.
We define the \emph{reduced Hilbert polynomial}
of $\F$
\begin{align*}
p_{\F}=\frac{P_{\F}}{\a_d (\F)}.
\end{align*}

\noindent \\
{\bf (2.1) Definition:} Let $\F$ be a coherent sheaf on $X$.
Assume that Supp$(\F)$ is pure dimensional of dimension $d$. 
We say that $\F$ is \emph{(semi)stable} if the following two
conditions are satisfied:
\begin{enumerate}
\item[(i)] $\F$ does not have nonzero subsheaves $\F'$ with
Supp$(\F')$ having dimension smaller than $d$;
\item[(ii)] for any proper subsheaf $\F'\subset \F$ we have
\begin{align*}
p_{\F'} (\le) < p_{\F}
\end{align*}
meaning that for $m$ sufficiently large the following inequality
holds:
\begin{align*}
p_{\F'}(m) (\le) < p_{\F}(m).
\end{align*}
\end{enumerate}

\noindent \\
{\bf (2.2) Remark:} We will be interested in semistable sheaves on
$\P^2$ with linear Hilbert polynomial $P_{\F}(m)=rm+\chi$.
Such sheaves are supported on projective
curves $C$ and the conditions
from the above definition take
the form:
\begin{enumerate}
\item[(i)] $\F$ does not have zero dimensional torsion;
\item[(ii)] for any proper subsheaf $\F'\subset \F$ we have
\begin{align*}
\frac{\a_0 (\F')}{\a_1 (\F')} (\le) < \frac{\a_0 (\F)}{\a_1 (\F)}.
\end{align*}
\end{enumerate}
We point out that $\F$ is a torsion $\O_{\P^2}$-module
because at every point $x$ there is
a nonzero germ of $\O_x$ vanishing on the support of $\F$,
hence annulating $\F_x$. The zero-dimensional torsion
of a sheaf is its largest subsheaf supported on finitely
many points.

The positive integer $r$ is the so-called multiplicity of $\F$
while $\chi$ is its Euler characteristic.
The restriction of $\F$ to a generic line in $\P^2$ is a
sheaf of length $r$ supported at finitely many points;
$r$ is also equal to the degree of $C$.
Here are more facts about such sheaves (compare
\cite{fr-diplom}, thm. 3.1):

\noindent \\
{\bf (2.3) Proposition:} \emph{ Let $\F$ be a semistable sheaf
on $\P^n$ with Hilbert polynomial $P_{\F}(m)=r m + \chi,\
0\le \chi < r$. Let $C$ be its support. Then:
\begin{enumerate}
\item[(i)] $\F$ is Cohen-Macaulay;
\item[(ii)] $\F$ is locally free on the smooth part of $C$;
\item[(iii)] $C$ has no zero dimensional components and no embedded
points;
\item[(iv)] if gcd$(r,\chi)=1$ then $\F$ is stable;
\item[(v)] if $h^0(\F(-1))=0$ then $h^1(\F(i))=0$ for
$i\ge r-\chi -1$.
\end{enumerate} }

\noindent \\
As a generic plane curve is smooth, we see that a generic
$\F$ from (2.3) is a line bundle supported on a smooth curve
of degree $r$. Its degree can be computed with the Riemann-Roch
formula: deg$(\F)=g(C)-1+\chi=\frac{r(r-3)}{2}+\chi$.
Line bundles supported on smooth curves are clearly stable
because their quotient sheaves are supported on finitely many
points, hence their proper subsheaves have the same
multiplicity but strictly smaller Euler characteristic.
Other, well-known,  examples of stable sheaves with
one-dimensional support are the structure sheaves $\O_C$, where
$C$ is any curve in $\P^2$ given as the zero-set of a polynomial
of degree $r$. We can see this using, for instance, lemma
(6.8): any indeal sheaf $\I \subset \O_C$ has Hilbert
polynomial
\bdm
P_{\O_C}(t)-P_{\O_{C'}}(t)-a=rt-\frac{r(r-3)}{2}-
r't+\frac{r'(r'-3)}{2}-a,
\edm
where $a \ge 0$ and $r'<r$ are integers. Thus
\bdm
\frac{\a_0(\I)}{\a_1(\I)}=\frac{-r-r'+3}{2}-\frac{a}{r-r'}
< \frac{-r+3}{2} = \frac{\a_0 (\O_C)}{\a_1 (\O_C)}.
\edm

\noindent \\
{\bf (2.4) Definition:} Let $\F$ be a semistable sheaf on $X$.
A \emph{Jordan-H\"older filtration} of $\F$ is a filtration by
subsheaves
\begin{align*}
0=\F_0 \subsetneqq \F_1 \subsetneqq \ldots \subsetneqq \F_r=\F
\end{align*}
such that all quotients $\F_i/\F_{i-1}$ are stable with reduced
Hilbert polynomial $p_{\F}$.
Two semistable sheaves $\F$ and $\G$ on $X$ are said to be
\emph{stable equivalent} if they posess Jordan-H\"older 
filtrations with isomorphic quotients. By this we mean the following:
there is a bijection between the set of quotients of the
filtration of $\F$ and the set of quotients of the filtration
of $\G$ such that the quotients corresponding via this bijection
are isomorphic.\\

Note that, in the case of stable sheaves, stable equivalence means
isomorphism. Any semistable sheaf $\F$ has at least one
Jordan-H\"older filtration. $\F$ may have more than one
filtration, however, the direct
sum $\oplus \F_i/\F_{i-1}$ does not depend on the filtration.

The moduli space of semistable sheaves on $X$ parametrizes stable
equivalence classes with fixed Hilbert polynomial.
It is a coarse moduli space for a certain moduli problem
which can be defined by means of the following
functor: Let us fix a numerical polynomial $P(m)$, i.e. a polynomial
with rational coefficients which takes integer values on the
integers. For any scheme $S$ over $k$ we define $\M_X (P)(S)$
as the set of equivalence classes $[\F]$ of $S$-flat coherent
sheaves $\F$ on $S\times X$ whose restriction $\F_s$ to any fiber
$\p^{-1}(s)$ is a semistable sheaf with Hilbert polynomial $P$.
Here $\p :S\times X\lra S$ is the projection onto the first
factor. Two sheaves $\F$ and $\G$ on $S\times X$ are said to be
equivalent
if there is a line bundle ${\mathcal L}$ on $S$ such that
$\F$ is isomorphic to $\G \tensor \p^* {\mathcal L}$.
Given a morphism $f:T\lra S$ of schemes over $k$ we have a map
\begin{align*}
\M_X (P)(f): \M_X (P)(S)\lra \M_X (P)(T)
\end{align*}
given by the pull-back
\begin{align*}
\M_X (P)(f)([\F])=[(f\times 1)^*\F].
\end{align*}
We have thus defined a countervariant functor $\M_X (P)$
from the category of schemes over $k$ to the category of
sets.

\noindent \\
{\bf (2.5) Definition:} A scheme M over $k$ is called
a \emph{coarse moduli space} of semistable sheaves on $X$
with Hilbert polynomial $P$ if there is a natural transformation
of functors
\begin{align*}
\M_X (P)(\_ )\stackrel{\t}{\lra} \text{Mor}(\_ ,\text{M})
\end{align*}
satisfying the following properties:
\begin{enumerate}
\item[(i)] the map
\begin{align*}
\t(\text{Spec}(k)):\M_X (P)(\text{Spec}(k))\lra
\text{Mor}(\text{Spec}(k),\text{M})
\end{align*}
is a bijection. In other words the set of closed points of
M is in a one-to-one correspondence with the set of
stable equivalence classes of semistable sheaves on $X$ with
Hilbert polynomial $P$;
\item[(ii)] given a scheme N and a natural transformation
\begin{align*}
\M_X (P)(\_ ) \stackrel{\t'}{\lra} \text{Mor}(\_ ,\text{N})
\end{align*}
there is a unique morphism $f:$ M$\lra $N such that
$f\circ \t (S)=\t' (S)$ for all $S$.
\end{enumerate}

\noindent \\
{\bf (2.6) Theorem:} \emph{Let $X$ be a smooth projective
variety with ample line bundle $\O_X (1)$. Let $P$ be a
numerical polynomial. Then:
\begin{enumerate}
\item[(i)] there exists a coarse moduli space \emph{M}$_X (P)$
of semistable sheaves on $X$ with Hilbert polynomial $P$;
\item[(ii)] \emph{M}$_X (P)$ is a projective scheme;
\item[(iii)] there is an open subscheme \emph{M}$^s_X(P)$ of
\emph{M}$_X (P)$ whose closed points parametrize the isomorphism
classes of stable sheaves with Hilbert polynomial $P$.\\
\end{enumerate}}

The theorem in its full generality was proven in \cite{simpson}.
Let us recall the way M$_X (P)$ is constructed:
For a suitably large integer $m$ let $V$ be a vector space
of dimension $P(m)$. We consider the quotient scheme
\begin{displaymath}
Q= \text{Quot} (X, V\tensor \O_X(-m), P)
\end{displaymath}
of coherent sheaves $\F$ on $X$ with Hilbert polynomial
$P$ which occur as quotients $V \tensor \O_X (-m) \lthra
\F$. The reductive group SL$(V)$ acts on $Q$ by its
action on the first component of $V \tensor \O_X (-m)$.
Inside $Q$ there is the open and SL$(V)$-invariant
subset $R$ of semistable points. 
Semistability here is meant in the sense of
Geometric Invariant Theory, which also guarantees the existence
of a categorical quotient (see (7.1)) of $R$ by SL$(V)$.
This quotient is the moduli space M$_X(P)$.

We now turn to the question under which circumstances the
moduli space of stable sheaves M$^s_X (P)$ is fine.
Fine moduli spaces represent certain functors, so we
define the countervariant functor
$\M^s_X (P)$ from the category of schemes over $k$ to the
category of sets in the same way as $\M_X (P)$ was defined,
with the difference that we now require each
restriction $\F_s$ to be stable.

\noindent \\
{\bf (2.7) Definition:} We say that M$^s_X (P)$ is a
\emph{fine moduli space} of stable sheaves on $X$ with Hilbert
polynomial $P$ if the natural transformation
\begin{align*}
\M^s_X (P)(\_ ) \stackrel{\t^s}{\lra}
\text{Mor}(\_ ,\text{M}^s_X (P))
\end{align*}
induced by $\t$ is an isomorphism of functors.
If this is true, let $\U$ be the sheaf on M$_X^s \times X$
whose class $[\U]\in \M_X^s (P)($M$^s_X(P))$ corresponds under
$\t$ to the identity map of M$_X^s (P)$. We say that
$\U$ is a \emph{universal family} on M$_X^s (P)$.

\noindent \\
{\bf (2.8) Remark:} The inverse of $\t^s (S)$ for a scheme $S$
is given by $\t^s (S)^{-1}(f)=[f^*\U]$. In fact, M$^s_X(P)$
is a fine moduli space if and only if there exists a coherent
sheaf $\U$ on M$^s_X(P)\times X$ such that:
\begin{enumerate}
\item[(i)] $\U$ is flat over M$_X^s(P)$;
\item[(ii)] for any point $[ \F ] \in $ M$^s_X (P)$ the restriction
of $\U$ to the fiber $[\F]\times X$ is isomorphic to $\F$;
\item[(iii)] $\U$ has the following universality property:
for any scheme $S$ over $k$ and any $S$-flat coherent family $\F$
of semistable sheaves on $S\times X$ with Hilbert polynomial
$P$, there exists a unique morphism
$f:S\lra $ M$^s_X(P)$ such that $\F \simeq f^* \U \tensor \p^*
{\mathcal L}$ for some line bundle ${\mathcal L}$ on $S$.
\end{enumerate}

\noindent \\
{\bf (2.9) Theorem:} \emph{Consider the numerical polynomial
\begin{align*}
P(m)=\sum_{i=0}^d \a_i {m+i-1 \choose i}.
\end{align*}
Assume that gcd$(\a_0,\ldots,\a_d)=1$. Then \emph{M}$^s_{\P^n}(P)$
is a fine moduli space for any $n\ge d$.}

\noindent \\
We refer to \cite{hl} for the proof of this theorem.
In this paper we will focus on moduli spaces of sheaves
on $\P^2$ with linear Hilbert polynomial $P_{\F}(m)=rm+\chi$.
To bring us closer to
the notations from \cite{lepotier}, where such moduli spaces
were systematically studied, we also write
M$_{\P^n}(r,\chi)$ instead of M$_{\P^n}(rm+\chi)$,
respectively M$^s_{\P^n}(r,\chi)$ instead of
M$^s_{\P^n}(rm+\chi)$. Combining (2.9) with (2.3)(iv)
we obtain:

\noindent \\
{\bf (2.10) Proposition:} \emph{Assume that gcd$(r,\chi)=1$.
Then \emph{M}$_{\P^n}(r,\chi)=$\emph{M}$^s_{\P^n}(r,\chi)$
is a fine moduli space.}

\noindent \\
From theorem 3.19(2) in \cite{lepotier} we learn that
M$^s_{\P^2}(r,\chi)$ is not a fine moduli space in the case
when $r$ and $\chi$ are not mutually prime. Let us quote the
precise statement (which is stronger):

\noindent \\
{\bf (2.11) Theorem:} \emph{If $r$ and $\chi$ are not mutually
prime, then for any open subset $U \subset $\emph{M}$_{\P^2}(r,\chi)$
there is no universal sheaf on $U \times \P^2$.} \\

As the spaces M$_{\P^2}(r,\chi)$ and M$_{\P^2}(r,r+\chi)$
are isomorphic, we will assume henceforth that $0 < \chi \le r$.
Theorem 3.1 and proposition 2.3 from \cite{lepotier}
yield the following:

\noindent \\
{\bf (2.12) Theorem:} \emph{For any integers $r \ge 1$ and
$\chi$ the moduli space \emph{M}$_{\P^2}(r,\chi)$ is irreducible,
of dimension $r^2+1$ and smooth on the open dense set
represented by stable sheaves.} \\

We finish this section with an easy observation about subsets
of moduli spaces:

\noindent \\
{\bf (2.13) Remark:} Let $\E$ be a locally free sheaf on $X$.
For any integers $i,\, j \ge 0$
the subset of M$_X(P)$ of stable equivalence classes
of sheaves $\F$ with $h^i (X,\F \tensor \E) \ge j$ is a closed
algebraic subset.

\noindent \\
\emph{Proof:} Using the notations preceeding (2.7),
we consider the universal
quotient $V \tensor \O_X (-m) \lthra \tilda{\F}$ on
$Q \times X$. The sheaf $\tilda{\F} \boxtimes \E$ is flat over
$Q$ so, according to the semicontinuity theorem, the
set $Y$ of equivalence classes of quotients
$V \tensor \O_X (-m) \lthra \F$ with $h^i (X,\F \tensor \E) \ge j$
is a closed algebraic subset in $Q$. Notice that $Y$ is
SL$(V)$-invariant so, by virtue of (7.2)(iv),
its image under the good quotient map
$R \lra $ M$_X(P)$ is closed. This image is precisely the subset
from the remark.

If M$_X (P)$ is a fine moduli space, then the remark
follows directly from the semicontinuity theorem applied to
$\U \boxtimes \E$.


\section{Semistable Morphisms of Sheaves}

Given coherent sheaves $\E$ and $\F$
on $\P^n$ the affine space $W =$Hom$(\E,\F)$ is acted upon by
the algebraic group $G=$Aut$(\E)\times$Aut$(\F)/k^*$.
Here $k^*$ is embedded as the group of homotheties
$\{ (t\cdot 1_{\E},\, t\cdot 1_{\F}),\ t\in k^*\}$. The action
is given by $(g,h).w=h\circ w \circ g^{-1}$.
If $G$ is reductive then Geometric Invariant Theory distinguishes
an open subset $W^{ss}\subset W$ of so-called semistable morphisms
and constructs a categorical quotient $W^{ss}//G$.
Our difficulty is that in general $G$ is not reductive.
A notion of semistability in the context of nonreductive
groups has been studied in \cite{drezet-trautmann} and
its usefullness has been made clear in the work \cite{drezet-2000}
of Dr\'ezet.

This section introduces Dr\'ezet and Trautmann's notion of
semistability and is mainly a reproduction of notations from
\cite{drezet-trautmann}. Let us fix sheaves
\begin{displaymath}
\E= \oplus_{1\le i\le r} M_i \tensor \E_i,\qquad
\F= \oplus_{1\le l\le s} N_l \tensor \F_l
\end{displaymath}
where $M_i,\ N_l$ are vector spaces over $k$ of dimensions
$m_i,\, n_l$ while $\E_i,\ \F_l$ are simple sheaves on $\P^n$,
meaning that their only endomorphisms are homotheties.
For our purposes $\E_i$ and $\F_l$ will be line bundles.
We assume that Hom$(\E_i,\E_j)=0$
when $i>j$ and Hom$(\F_l,\F_m)=0$ when $l>m$.
We denote
\begin{eqnarray*}
H_{li} & = & \text{Hom}(\E_i,\F_l), \\
A_{ji} & = & \text{Hom}(\E_i,\E_j), \\
B_{ml} & = & \text{Hom}(\F_l,\F_m).
\end{eqnarray*}
The group $G$ consists of pairs of matrices $(g,h)$,
\begin{align*}
g=\left[
\begin{array}{ccccc}
g_1     &  0    &  \cdots  &  \cdots  & 0 \\
u_{21}  &  g_2  &  \cdots  &  \cdots  & 0 \\
\vdots  &       &  \ddots  &          & \vdots \\
\vdots  &       &          & \ddots   & \vdots \\
u_{r1}  &       &  \cdots  &  u_{r,r-1} & g_r
\end{array}
\right]
,\qquad
h=\left[
\begin{array}{ccccc}
h_1     &  0    &  \cdots  &  \cdots  & 0 \\
v_{21}  &  h_2  &  \cdots  &  \cdots  & 0 \\
\vdots  &       &  \ddots  &          & \vdots \\
\vdots  &       &          & \ddots   & \vdots \\
v_{s1}  &       &  \cdots  &  v_{s,s-1} & h_s
\end{array}
\right]
\end{align*}
with
\begin{align*}
g_i\in \text{GL}(M_i),\ \ h_l \in \text{GL}(N_l),\ \
u_{ji}\in \text{Hom}(M_i, M_j \tensor A_{ji}), \ \
v_{ml}\in \text{Hom}(N_l, N_m \tensor B_{ml}).
\end{align*}
The conditions $u_{ij}=0$ and $v_{ml}=0$ define a
reductive subgroup $G_{red}$ inside $G$.

For fixed positive integers $\l_i,\, \m_l$
we consider the character $\chi$ of $G$ given by
\begin{align*}
\chi (g,h)=\prod_{1\le i\le r} \text{det}(g_i)^{-\l_i} \cdot
\prod_{1\le l \le s} \text{det}(h_l)^{\m_l}.
\end{align*}
Since $\chi$ must be trivial on the subgroup of homotheties
$k^*$, we impose the relation
\begin{align*}
\sum_{1\le i \le r} m_i \l_i = \sum_{1\le l \le s} n_l \m_l
\end{align*}
and we denote by $d$ this sum.
We will call a \emph{polarization} the tuple
\begin{align*}
\L =(\l_1,\ldots,\l_r, \m_1,\ldots,\m_s).
\end{align*}

\noindent \\
{\bf (3.1) Definition:} Let $\L$ be a fixed polarization.
A point $w\in W$ is called:

\begin{enumerate}
\item[(i)] \emph{semistable} with respect to $G_{red}$ and $\L$
if there are $n\ge 1$ and a polynomial $f\in k[W]$ satisfying
$f(g.x)= \chi^n (g) f(x)$ for all $g\in G_{red},\ x\in W$,
such that $f(w)\neq 0$;
\item[(ii)] \emph{stable} with respect to $G_{red}$ and $\L$
if Stab$_{G_{red}}(w)$ is zero dimensional and there is $f$
as above but with the additional property that the action
of $G_{red}$ on $W_f=\{ x\in W,\ f(x)\neq 0 \}$ is closed;
\item[(iii)] \emph{properly semistable} if it is semistable
but not stable. \\
\end{enumerate}

The question is now how to define semistability with
respect to $G$. The key is the following observation from
Geometric Invariant Theory: let $T$ be a maximal torus
inside a reductive algebraic group which acts on a projective
variety. Then a point on the variety is semistable if and
only if all points in its orbit are semistable with respect
to $T$. In our context the subgroup of diagonal matrices
is a maximal torus inside both $G_{red}$ and $G$.
This justifies the following:

\noindent \\
{\bf (3.2) Definition:} A point $w\in W$ is called \emph{(semi)stable}
with respect to $G$ and $\L$ if $g.w$ is (semi)stable with respect to
$G_{red}$ and $\L$ for all $g\in G$. We denote by $W^{ss}(G,\L),\
W^s (G,\L)$ the corresponding sets. \\

To describe the sets of semistable points in concrete situations
we will use a very special case of King's criterion of semistability
as formulated in \cite{drezet-trautmann}. Let us write
\begin{align*}
\E = \oplus_j \, \E_j',
\qquad
\F = \oplus_m \, \F_m'
\end{align*}
where $\E_j',\, \F_m'$ are line bundles. We represent $w$ by a
matrix $(w_{mj})$ with $w_{mj} \in \ $Hom$(\E_j',\F_m')$.
We put
\begin{align*}
\l_j'= \l_i \quad \text{if} \quad \E_j'\isom \E_i,
\qquad \qquad
\m_m'= \m_l \quad \text{if} \quad \F_m' \isom \F_l.
\end{align*}

\noindent \\
{\bf (3.3) Proposition:} \emph{A morphism $w \in W$ is (semi)stable
with respect to $G$ and $\L$
if and only if for all $g\in G$ and for any zero submatrix
$((g.w)_{m,j})_{m\in M, j\in J}$ we have}
\begin{align*}
\sum_{m\in M} \m_m' (\le) < \sum_{j \notin J} \l_j'.
\end{align*}

\noindent \\
For convenience we replace each $\l_i$ with $\l_i/d$ and each
$\m_l$ with $\m_l/d$. Thus our polarization $\L$ will be a tuple
of rational numbers satisfying
\begin{align*}
\tag{3.4}
\sum_{i=1}^r m_i \l_i = 1 = \sum_{l=1}^s n_l \m_l.
\end{align*}
The set of polarizations can be realized as an open subset
of the Euclidean space of dimension $r+s-2$.


\section{Sheaves $\F$ with $h^0(\F(-1))=0$ and $h^1(\F)=0$}

The main technical toll that we will use in this paper
is the Beilinson complex. Given a coherent sheaf $\F$ on
$\P^2$ there is a sequence of sheaves
\begin{align*}
0\lra \C^{-2} \lra \C^{-1} \lra \C^0 \lra \C^1 \lra \C^2 \lra 0
\end{align*}
which is exact except in the middle where the cohomology is
$\F$. The sheaves $\C^i$ are given by
\begin{eqnarray*}
\C^{-2} & = & H^0(\F\tensor \Om^2 (2))\tensor \O(-2), \\
\C^{-1} & = & H^0(\F\tensor \Om^1(1))\tensor \O(-1)
                \oplus H^1(\F\tensor \Om^2(2))\tensor \O(-2), \\
\C^0    & = & H^0(\F)\tensor \O \oplus H^1(\F\tensor \Om^1(1))
                  \tensor \O(-1) \oplus H^2(\F\tensor \Om^2(2))
                     \tensor \O(-2), \\
\C^1    & = & H^1(\F)\tensor \O \oplus H^2(\F\tensor \Om^1(1))
                \tensor \O(-1), \\
\C^2    & = & H^2 (\F) \tensor \O.
\end{eqnarray*}
The sheaves $\F$ we are interested in are supported on curves, so
\begin{align*}
H^2(\F)=0,\quad H^2 (\F\tensor \Om^1(1))=0,\quad H^2 (\F\tensor
\Om^2(2))=0.
\end{align*}
Also, on $\P^2$ we have $\Om^2(2)=\O(-1)$. The Beilinson sequence
that we will use takes the form
\begin{align*}
\tag{4.1}
0\lra \C^{-2} \lra \C^{-1} \lra \C^0 \lra \C^1 \lra 0
\end{align*}
where
\begin{eqnarray*}
\C^{-2} & = & H^0(\F(-1))\tensor \O(-2), \\
\C^{-1} & = & H^0(\F\tensor \Om^1(1))\tensor \O(-1)
                \oplus H^1(\F(-1))\tensor \O(-2), \\
\C^0    & = & H^0(\F)\tensor \O \oplus H^1(\F\tensor \Om^1(1))
                  \tensor \O(-1), \\
\C^1    & = & H^1(\F)\tensor \O.
\end{eqnarray*}
The morphisms 
\begin{eqnarray*}
H^0(\F(-1))\tensor \O(-2) & \lra & 
H^1(\F(-1))\tensor \O(-2), \\
H^0(\F\tensor \Om^1(1)) \tensor \O(-1) & \lra &
H^1(\F\tensor \Om^1(1)) \tensor \O(-1), \\
H^0(\F) \tensor \O & \lra &
H^1(\F) \tensor \O
\end{eqnarray*}
from above are all zero. Indeed, each of these morphisms
can be represented by a matrix with scalar entries.
Performing Gaussian elimination on these matrices
we arrive at a complex like (4.1) in which the
cohomology vector spaces get replaced by subspaces.
Using standard methods in cohomology theory we can show
that the dimension of each of these subspaces equals
the dimension of the corresponding cohomology space.
In other words no Gaussian elimination was performed;
the matrices were zero to begin with.

Apart from the question of the semistability of the morphism $\f$,
the resolutions of generic sheaves giving points in
M$_{\P^2}(r,\chi)$ apper first in \cite{fr-diplom}.
For the sake of completeness we have
included the really simple arguments here without quoting every
time the above work.

\noindent \\
{\bf (4.2) Claim:} \emph{Let $\F$ be a sheaf on $\P^2$ with
$h^0(\F(-1))=0$. Let $n\ge 2$ be an integer and assume that $\F$ has
Hilbert polynomial $P_{\F}(t)=(n+1)t+n$. Then $\F$ is semistable
if and only if it has a resolution}
\begin{align*}
0\lra \O(-2)\oplus (n-1)\O(-1)\stackrel{\f}{\lra} n\O \lra
\F \lra 0
\end{align*}
\emph{with $\f$ not equivalent to a matrix of the form}
\begin{align*}
\left[
\begin{array}{cc}
\star & \psi \\
\star & 0
\end{array}
\right]
\end{align*}
\emph{where $\psi :m\O(-1) \lra m\O,\ 1 \le m \le n-1$.}

\emph{Equivalently, $\F$ is semistable if and only if it
has a resolution as above with $\f$ semistable with respect
to $\La$. Here $\La=(\l_1,\l_2,\m_1)$ is any polarization
satisfying $0< \l_1 < \frac{1}{n}$.}

\noindent \\
\emph{Proof:} We have $h^0 (\F(-1))=0,\ h^1(\F(-1))=1,\
h^0(\F)=n,\ h^1(\F)=0$ because of (2.3).
Thus (4.1) gives the following resolution with $\f_{12}=0$:
\begin{align*}
0 \lra \O(-2)\oplus (m+n-1)\O(-1) \stackrel{\f}{\lra}
m\O(-1) \oplus n\O \lra \F \lra 0.
\end{align*}
Here $m$ is an integer and, since $\f$ is injective, we
can only have $m=0$ or $m=1$.

Assume that $\F$ is semistable. Then $m \neq 1$, otherwise
$\F$ would have a subsheaf $\F'$ with resolution
\begin{align*}
0 \lra n\O(-1) \lra n\O \lra \F' \lra 0.
\end{align*}
We have $P_{\F'}(t)=nt+n$, hence such a subsheaf would
destabilize $\F$. Thus far we have obtained a resolution
\begin{align*}
0 \lra \O(-2) \oplus (n-1)\O(-1) \stackrel{\f}{\lra}
n\O \lra \F \lra 0.
\end{align*}
Let us mention that this resolution was first obtained
in \cite{maican} and it is also present in \cite{fr-diplom}.
The matrix $\f$ cannot be equivalent to a matrix of the form
\begin{align*}
\left[
\begin{array}{cc}
\star & \psi \\
\f_{21} & 0
\end{array}
\right]
\end{align*}
otherwise we would get an exact commutative diagram
\begin{displaymath}
\xymatrix
{
 & 0 \ar[d] & 0 \ar[d] \\
0 \ar[r] & m \O(-1) \ar[r]^{\psi} \ar[d]_{\left[ \begin{array}{cc}
   \small{0} \\ I_m  \end{array} \right]} & m\O \ar[r] \ar[d]^{\left[
        \begin{array}{c} I_m \\ 0 \end{array} \right]} & \F'
           \ar[r] & 0 \\
0 \ar[r] & \O(-2) \oplus (n-1)\O(-1) \ar[r]^{\quad \qquad \f} \ar[d]_{\left[
   \begin{array}{cc} I & 0 \end{array} \right] } & n\O \ar[r]
      \ar[d]^{\left[ \begin{array}{cc} 0 & I \end{array} \right]}
            & \F \ar[r] & 0 \\
0 \ar[r] & \O(-2) \oplus (n-m-1)\O(-1) \ar[r]^{\qquad \quad \f_{21}} \ar[d] 
    & (n-m)\O \ar[d] \\
 & 0 & 0
}
\end{displaymath}
in which $\F'$ is a destabilizing subsheaf of $\F$.

Conversely, we assume that $\F$ has a resolution as in the
statement of the claim, and we try to show that the conditions
from (2.2) are satisfied. At every point $x$ in the support of
$\F$ we have
\begin{align*}
\text{depth}_x \, \F = 2 - pd_x \, \F \ge 1,
\end{align*}
showing that $\F$ does not have zero-dimensional torsion.
Assume now that $\F$ has a subsheaf $\F'$ which contradicts
(2.2)(ii), in other words which satisfies
\begin{align*}
\frac{\a_0 (\F')}{\a_1 (\F')} > \frac{n}{n+1}.
\end{align*}
The multiplicity $m=\a_1 (\F')$ cannot exceed the multiplicity
of $\F$. Thus $\a_0 (\F') \ge m$. Since $h^0(\F'(-1)) \le
h^0(\F(-1)) = 0$ we have $P_{\F'}(-1)=-h^1(\F'(-1)) \le 0$
forcing $\a_0 (\F') \le m$. So far we have obtained
$P_{\F'}(t)=mt+m$ for some integer $1 \le m \le n$.
Now let us notice that $\F$ is generated by global sections,
so we must have $h^0 (\F')\le n-1$, forcing $m \le n-1$.
We have
\begin{align*}
h^0 (\F'(-2))=0, \quad h^1 (\F'(-2))=m, \quad
h^0 (\F'(-1))=0, \quad h^1 (\F'(-1))=0.
\end{align*}
The Beilinson sequence of $\F'(-1)$ gives a resolution
\begin{align*}
0 \lra m\O(-2) \lra m \O(-1) \lra \F'(-1) \lra 0.
\end{align*}
This yields a commutative diagram
\begin{displaymath}
\xymatrix
{
0 \ar[r] & m \O(-1) \ar[r]^{\psi} \ar[d]_{\b} & m\O \ar[r]
      \ar[d]_{\a} & \F' \ar[r] \ar[d] & 0 \\
0 \ar[r] & \O (-2) \oplus (n-1)\O(-1) \ar[r]^{\qquad \quad \f}
    & n\O \ar[r] & \F \ar[r] & 0
}.
\end{displaymath}
The map $\a$ is injective because it is injective on the level
of global sections. Hence also $\b$ is injective.
It is clear now that $\f$ is equivalent to a matrix of the form
\begin{align*}
\left[
\begin{array}{cc}
\star & \psi \\
\star & 0 
\end{array}
\right] .
\end{align*}
This contradicts the hypothesis and finishes the proof of
the first part of the claim.

The second part of the claim follows from (3.3). Namely,
King's criterion says that $\f$ is semistable with respect
to $\La$ if and only if whenever
\begin{align*}
\f \sim 
\left[
\begin{array}{cc}
\star & \psi \\
\star & 0 
\end{array}
\right]
\qquad \text{with} \quad
\psi : p\O(-2)\oplus q\O(-1) \lra m\O
\end{align*}
we have $m\m_1 \ge p\l_1 + q \l_2$. Thus, we need to find
$\La$ satisfying the conditions
\begin{align*}
m\m_1 < p\l_1 + q\l_2 \quad \text{if and only if}
\quad q \ge m.
\end{align*}
Here $0\le m \le n,\ 0 \le p \le 1,\ 0 \le q \le n-1$.
These conditions are the same as
\begin{eqnarray*}
m\m_1 & < & m \l_2 \quad \text{for} \quad 1 \le m \le n-1,\\
m\m_1 & \ge & \l_1 + (m-1)\l_2 \quad \text{for} \quad
                                        1 \le m \le n.
\end{eqnarray*}
Using relations (3.4)
\begin{align*}
\m_1 = \frac{1}{n}, \qquad \l_2 = \frac{1-\l_1}{n-1},
\end{align*}
we arrive at the conditions
\begin{align*}
\l_1 < \frac{1}{n},
\qquad
\frac{n-m}{n(n-1)} \ge \l_1 \frac{n-m}{n-1}.
\end{align*}
The conditions on $\l_1$ are precisely those of the claim.
Q.e.d.

\noindent \\
{\bf (4.3) Claim:} \emph{Let $\F$ be a sheaf on $\P^2 =\P(V)$
with $h^0(\F(-1))=h^1(\F)=0$. Assume that $\F$ has Hilbert
polynomial $P(t)=(n+2)t+n$ where $n\in \{3,\, 4,\, 5,\, 6\}$.
Then $\F$ is semistable if and only if it has a resolution}
\begin{align*}
\tag{i}
0 \lra 2\O(-2) \oplus (n-2)\O(-1) \stackrel{\f}{\lra}
n\O \lra \F \lra 0
\end{align*}
\emph{with $\f$ not equivalent to a matrix of the form}
\begin{align*}
\left[
\begin{array}{cc}
\star & \psi \\
\star & 0
\end{array}
\right]
\end{align*}
\emph{where $\psi : m\O(-1)\lra m\O,\ 1 \le m \le n-2$
or $\psi : \O(-2)\oplus (m-1)\O(-1)\lra m\O,\ n/2
< m \le n-1$, or it has a resolution}
\begin{align*}
\tag{ii}
0 \lra 2\O(-2) \oplus (n-1)\O(-1) \stackrel{\f}{\lra}
\O(-1) \oplus n\O \lra \F \lra 0
\end{align*}
\emph{whith $\f_{11}$ having linearly independent entries,
$\f_{12}=0$, $\f_{21}\neq 0$ and $\f_{22}$ not equivalent to a matrix
of the form}
\begin{align*}
\left[
\begin{array}{cc}
\star & \psi \\
\star & 0
\end{array}
\right]
\end{align*}
\emph{where $\psi :m\O(-1) \lra m\O,\ 1 \le m \le n-1$.}\\

\emph{The maps $\f$ occuring in (i) are precisely those maps
$\f \in W^{ss}(G,\L)$ with nonzero determinant.
Here $\L=(\l_1,\l_2,\m_1)$ is any polarization satisfying
$\frac{1}{2n} \le \l_1 < \frac{1}{n}$.}

\emph{The maps $\f$ occuring in (ii) are precisely those maps
$\f \in W^{ss}(G,\L)$ with det$(\f)\neq 0$ and $\f_{12}=0$.
Here $\L=(\l_1,\l_2,\m_1,\m_2)$ is any polarization
for which the pair $(\l_1,\m_1)$ is in the interior of
the triangle with vertices}
\begin{align*}
(0,0), \quad \left( \frac{1}{n+1}\, ,\, \frac{1}{n+1} \right),
\quad \left( \frac{1}{n^2-n+2}\, , \, \frac{2}{n^2-n+2} \right).
\end{align*}

\emph{When $n \ge 7$ solely the ``only if'' part of the above
statement is true. Thus, all we can say in the case $n \ge 7$,
is that each semistable sheaf $\F$ occurs as the cokernel
of a semistable $\f$, but there are semistable
morphisms $\f$ whose cokernel is not a semistable sheaf.}

\noindent \\
\emph{Proof:} One direction is clear, cf. the proof of (4.2).
Conversely, suppose that $\F'\subset \F$ is a destabilizing
subsheaf. Arguing as in (4.2) we see that the Hilbert polynomial
of $\F'$ is either $mt+m$ with $1 \le m \le n$ or $(m+1)t+m$
with $\frac{n}{2} < m \le n$. In the case $P_{\F'}(t)=mt+m$
we deduce that $\f$ is equivalent to a matrix of the form
\begin{align*}
\left[
\begin{array}{cc}
\star & \psi \\
\star & 0
\end{array}
\right]
\end{align*}
where $\psi :m\O(-1) \lra m\O$. Assume now that
$P_{\F'}(t)=(m+1)t+m$ with $\frac{n}{2} < m \le n$. We have
\begin{align*}
h^0(\F'(-1))=0, \quad h^1(\F'(-1))=1, \quad
h^0(\F'(-2))=0, \quad h^1(\F'(-2))=m+2.
\end{align*}
From $\Om^1 \subset 3 \O(-1)$ and $h^0(\F'(-1))=0$ we
get $h^0(\F' \tensor \Om^1)=0$. The Beilinson
sequence of $\F'(-1)$, which has $\F'(-1)$ as middle
cohomology, takes the form
\begin{align*}
0 \lra (m+2)\O(-2) \lra (m+3)\O(-1) \stackrel{\eta}{\lra}
\O \lra 0.
\end{align*}
Since $\eta$ is surjective it must be equivalent to a matrix
of the form
\begin{align*}
(X, Y, Z, 0, \ldots, 0).
\end{align*}
So far we have arrived at the following resolution of $\F'$:
\begin{align*}
0 \lra (m+2)\O(-1) \lra \Om^1(1) \oplus m\O \lra \F' \lra 0.
\end{align*}
From this we get $h^1(\F')=0$, $h^0(\F')=m$.
Writing $p=h^1(\F' \tensor \Om^1(1))$, the sequence (4.1)
gives the resolution
\begin{align*}
0 \lra \O(-2) \oplus (m+p-1) \O(-1) \stackrel{\psi}{\lra}
p\O(-1) \oplus m\O \lra \F' \lra 0,
\end{align*}
with $\psi_{12}=0$. From the injectivity of $\psi$ we see
that we can only have $p=0$ or $p=1$. In the latter case
$\F'$ has a subsheaf $\F''$ with resolution
\bdm
0 \lra m\O(-1) \lra m\O \lra \F'' \lra 0.
\edm
This situation has been examined before.
Thus we arrive at the resolution
\begin{align*}
0 \lra \O(-2) \oplus (m-1)\O(-1) \stackrel{\psi}{\lra}
m\O \lra \F' \lra 0.
\end{align*}
We get the following exact commutative diagrams in case (i)
\begin{displaymath}
\xymatrix
{
0 \ar[r] & \O(-2) \oplus (m-1)\O(-1) \ar[r]^{\ \ \ \ \ \ \ \ \ \ \ \psi}
\ar[d]^{\b} & m\O \ar[r] \ar[d]^{\a} & \F' \ar[r] \ar[d] & 0 \\
0 \ar[r] & 2\O(-2) \oplus (n-2)\O(-1) \ar[r]^{\ \ \ \ \ \ \ \ \ \ \ \f}
& n\O \ar[r] & \F \ar[r] & 0
} ,
\end{displaymath}
and in case (ii)
\begin{displaymath}
\xymatrix
{
0 \ar[r] & \O(-2) \oplus (m-1)\O(-1) \ar[r]^{\ \ \ \ \ \ \ \ \ \ \psi}
\ar[d]^{\b} & m\O \ar[r] \ar[d]^{\a} & \F' \ar[r] \ar[d] & 0 \\
0 \ar[r] & 2\O(-2) \oplus (n-1)\O(-1) \ar[r]^{\ \ \ \ \ \ \f} & \O(-1) \oplus 
n\O \ar[r] & \F \ar[r] & 0
} .
\end{displaymath}
The map $\a$ is injective because it is injective on global
sections. It follows that $\b$ is also injective.
If $\b_{11} \neq 0$, which can happen only in case (i),
we get the contradiction
\begin{align*}
\f \sim \left[
\begin{array}{cc}
\star & \psi \\
\star & 0
\end{array}
\right] .
\end{align*}
In case (ii) we have $\f_{11}\b_{11} = \a_{11}\psi_{11} =0$
forcing $\b_{11}=0$ because, by hypothesis, the entries of
$\f_{11}$ are linearly independent.
Assume from now on that $\b_{11}=0$.
The case
\begin{align*}
\left[
\begin{array}{cc}
\b_{21} & \b_{22}
\end{array}
\right] \sim \left[
\begin{array}{cccc}
X & 0 & \cdots & 0 \\
Y & 0 & \cdots & 0 \\
Z & 0 & \cdots & 0 \\
0 & 1 & \cdots & 0 \\
\vdots & \vdots & \ddots & \vdots \\
0 & 0 & \cdots & 1 \\
0 & 0 & \cdots & 0 \\
\vdots & \vdots & & \vdots \\
0 & 0 & \cdots & 0
\end{array}
\right]
\end{align*}
leads to $n-1 \ge 3+m-1$, so $n \ge 3+m > 3 + \frac{n}{2}$
which contradicts the hypothesis $n \le 6$.
The case
\begin{align*}
\left[
\begin{array}{cc}
\b_{21} & \b_{22}
\end{array}
\right] \sim \left[
\begin{array}{cccc}
X & 0 & \cdots & 0 \\
0 & 1 & \cdots & 0 \\
\vdots & \vdots & \ddots & \vdots \\
0 & 0 & \cdots & 1 \\
0 & 0 & \cdots & 0 \\
\vdots & \vdots & & \vdots \\
0 & 0 & \cdots & 0
\end{array}
\right]
\end{align*}
is excluded from the fact that ${\mathcal Coker}(\b)$,
as a subsheaf of the torsion-free sheaf ${\mathcal Coker}(\a)$,
must be torsion-free. We are left with the case
\begin{align*}
\left[
\begin{array}{cc}
\b_{21} & \b_{22}
\end{array}
\right] \sim \left[
\begin{array}{cccc}
X & 0 & \cdots & 0 \\
Y & 0 & \cdots & 0 \\
0 & 1 & \cdots & 0 \\
\vdots & \vdots & \ddots & \vdots \\
0 & 0 & \cdots & 1
\end{array}
\right] . \quad \text{We get} \qquad
\f \sim \left[
\begin{array}{cc}
\star & \psi' \\
\star & 0
\end{array}
\right]
\end{align*}
with $\psi'$ an $(m+1)\times (m+1)$-matrix with entries in $V^*$.
This again contradicts the hypothesis and shows that $\F$ is
semistable. \\

The part of the claim concerning the semistability of $\f$
follows from (3.3). Namely, in case (i), we are looking for
$\L$ satisfying
\begin{eqnarray*}
m\m_1 & < & m\l_2 \quad \text{for} \quad 1 \le m \le n-2, \\
m\m_1 & < & \l_1 + (m-1)\l_2 \quad \text{for} \quad \frac{n}{2}
< m \le n-1, \\
m\m_1 & \ge & (m-1)\l_2 \quad \text{for} \quad \frac{n}{2}
< m \le n-1, \\
m\m_1 & \ge & \l_1 + (m-1)\l_2 \quad \text{for} \quad
1 \le m \le \frac{n}{2}, \\
m\m_1 & \ge & 2\l_1 +(m-2)\l_2 \quad \text{for} \quad
2 \le m \le n.
\end{eqnarray*}
Using relations (3.4) the above conditions become
$\frac{1}{2n} \le \l_1 < \frac{1}{n}$. Similarly, in case (ii),
we need to find $\L$ satisfying the conditions
\begin{eqnarray*}
\m_1 & < & 2\l_1, \\
m\m_2 & < & m\l_2 \quad \text{for} \quad 1 \le m \le n-1, \\
\m_1 & > & \l_1, \\
m\m_2 & \ge & (m-1)\l_2 \quad \text{for} \quad 1 \le m \le n, \\
\m_1 + m\m_2 & < & \l_1 + m\l_2 \quad \text{for} \quad
1 \le m \le n-1, \\
\m_1 + m\m_2 & \ge & 2 \l_1 + (m-1)\l_2 \quad \text{for}
\quad 1 \le m \le n.
\end{eqnarray*}
Using (3.4) the above conditions can be translated into the
conditition that $(\l_1,\m_1)$ is in the interior of a triangle
as in the statement of the claim. Q.e.d.

\noindent \\
{\bf (4.4) Observation:} If $n \ge 7$ then the situation
\begin{align*}
\left[
\begin{array}{cc}
\b_{21} & \b_{22}
\end{array}
\right] \sim \b_0 = \left[
\begin{array}{cccc}
X & 0 & \cdots & 0 \\
Y & 0 & \cdots & 0 \\
Z & 0 & \cdots & 0 \\
0 & 1 & \cdots & 0 \\
\vdots & \vdots & \ddots & \vdots \\
0 & 0 & \cdots & 1 \\
0 & 0 & \cdots & 0 \\
\vdots & \vdots & & \vdots \\
0 & 0 & \cdots & 0
\end{array}
\right]
\end{align*}
is feasible. Thus, to ensure the semistability of $\F$,
we would have to exclude, say in case (i), matrices of the form
\begin{align*}
\left[
\begin{array}{ccc}
\star & 0 & \psi'' \\
\star & \psi' & 0 \\
\star & 0 & 0
\end{array}
\right]
\qquad \text{with} \quad
\psi' \sim \left[
\begin{array}{ccc}
Y & X & 0 \\
Z & 0 & X \\
0 & -Z & Y
\end{array}
\right]
\end{align*}
and $\psi''$ an $m \times (m+1)$-matrix with entries in $V^*$.
From (3.3) we see that such conditions cannot be formulated
in terms of semistability, so they are beyond the interest of this paper.
Indeed, according to (3.3), semistability
conditions on a matrix specify that the matrix should not, up to
equivalence, have certain zero submatrices.

\noindent \\
{\bf (4.5) Claim:} \emph{Let $\F$ be a sheaf on $\P^2 =\P(V)$
with $h^0(\F(-1))=h^1(\F)=0$. Assume that $\F$ has Hilbert
polynomial $P(t)=4t+2$.
Then $\F$ is semistable if and only if it has a resolution}
\begin{align*}
\tag{i}
0 \lra 2\O(-2) \lra 2\O \lra \F \lra 0
\end{align*}
\emph{or it has a resolution}
\begin{align*}
\tag{ii}
0 \lra 2\O(-2) \oplus \O(-1) \stackrel{\f}{\lra}
\O(-1) \oplus 2\O \lra \F \lra 0
\end{align*}
\emph{with}
\begin{align*}
\f = \left[
\begin{array}{ccc}
X_1 & X_2 & 0 \\
\star & \star & Y_1 \\
\star & \star & Y_2
\end{array}
\right]
\end{align*}
\emph{where $X_1,\, X_2 \in V^*$ are linearly independent
one-forms and, likewise, $Y_1,\, Y_2 \in V^*$ are linearly
independent. These morphisms are precisely the morphisms
$\f \in W^{ss}(G,\L)$ with nonzero determinant and $\f_{12}=0$.
Here $\L = (\l_1,\l_2,\m_1,\m_2)$ is any polarization for which
the pair $(\l_1,\m_1)$ belongs to the interior of the quadrilater
with vertices}
\begin{align*}
(0\, , \, 0), \quad \left( \frac{1}{3} \, , \, \frac{1}{3} \right),
\quad \left( \frac{1}{2} \, , \, 1 \right), \quad (0 \, , \, 1).
\end{align*}

\noindent \\
\emph{Proof:} For the first part of the claim the proof is
the same as the proof of (4.3) so we omit it.
For the part of the claim concerning the semistability of
$\f$ we arrive at the following conditions on $\L$:
\begin{align*}
\m_1 > \l_1 \ \ \ \text{and} \ \ \ \m_2 < \l_2 \ \ \
\text{which is the same as} \ \ \ \m_1 > 4\l_1 -1 .
\end{align*}
They describe the quadrilater from the claim. Q.e.d. \\

In the remaining part of this section we will be concerned
with sheaves $\F$ on $\P^2$
having Hilbert polynomial $P(t)= (n+3)t+n,\ n\ge 3$,
and satisfying $h^0(\F(-1))=0$, $h^1(\F)=0$.
Such a sheaf $\F$ has one of the following resolutions:
\begin{align*}
\tag{i}
0 \lra 3\O(-2) \oplus (n-3)\O(-1) \stackrel{\f}{\lra}
n\O \lra \F \lra 0,
\end{align*}
\begin{align*}
\tag{ii}
0 \lra 3\O(-2) \oplus (n-2)\O(-1) \stackrel{\f}{\lra}
\O(-1) \oplus n\O \lra \F \lra 0,
\end{align*}
\begin{align*}
\tag{iii}
0 \lra 3\O(-2) \oplus (n-1)\O(-1) \stackrel{\f}{\lra}
2\O(-1) \oplus n\O \lra \F \lra 0,
\end{align*}
with $\f_{12}=0$ in cases (ii) and (iii).

\noindent \\
{\bf (4.6) Claim:} \emph{Let $\F$ be a sheaf with resolution (i)
and $3 \le n \le 6$. Then $\F$ is semistable if and only if
$\f$ is not equivalent to a matrix of the form}
\begin{align*}
\left[
\begin{array}{cc}
\star & \psi \\
\star & 0
\end{array}
\right]
\end{align*}
\emph{where}
\begin{eqnarray*}
& & \psi : m\O(-1) \lra m\O, \qquad 1 \le m \le n-3,
\quad \text{\emph{or}} \\
& & \psi : \O(-2) \oplus (m-1)\O(-1) \lra m\O, \qquad \frac{n}{3} < m
\le n-2, \quad \text{\emph{or}} \\
& & \psi : 2\O(-2) \oplus (m-2)\O(-1) \lra m\O, \qquad \frac{2n}{3}
< m \le n-1. \\
\end{eqnarray*}
\emph{Thus, any sheaf $\F$ with Hilbert polynomial $6t+3$
and resolution}
\begin{align*}
0 \lra 3\O(-2) \stackrel{\f}{\lra} 3\O \lra \F \lra 0
\end{align*}
\emph{is semistable. The morphisms $\f$ occuring form the open
subset of $W^{ss}(G,\L)$ given by the condition det$(\f)\neq 0$.
Here $\L$ is the only admissible polarization, namely
$\L=(\frac{1}{3},\frac{1}{3})$.}

\emph{A sheaf $\F$ with Hilbert polynomial $(n+3)t+n,\quad
n\in \{ 4,5,6 \}$, and resolution}
\begin{align*}
0 \lra 3\O(-2) \oplus (n-3)\O(-1) \stackrel{\f}{\lra}
n\O \lra \F \lra 0
\end{align*}
\emph{is semistable if and only if $\f$ is semistable with
respect to any polarization $\L=(\l_1,\l_2,\m_1)$ satisfying
$\frac{2}{3n} \le \l_1 < \frac{1}{n}$.}

\noindent \\
\emph{Proof:} One direction is clear. For the other direction
suppose that $\F' \subset \F$ is a destabilizing subsheaf.
The Hilbert polynomial of $\F'$ must be one of the following:
$mt+m$, with $1\le m \le n$, $(m+1)t+m$ with $\frac{n}{3}
< m \le n$, $ (m+2)t+m$ with $\frac{2n}{3} < m \le n$.

In the case $P_{\F'}(t)=mt+m$ we deduce, as in the proof
of (4.2), that
\begin{align*}
\f \sim \left[
\begin{array}{cc}
\star & \psi \\
\star & 0 
\end{array}
\right] \qquad \text{with} \quad \psi :m\O(-1)\lra m\O.
\end{align*}
In the case $P_{\F'}(t)=(m+1)t+m$ we arrive, as in the proof
of (4.3), at
\begin{align*}
\f \sim \left[
\begin{array}{cc}
\star & \psi \\
\star & 0 
\end{array}
\right]
\end{align*}
with $\psi : \O(-2) \oplus (m-1)\O(-1) \lra m\O$ or
$\psi :(m+1)\O(-1) \lra (m+1)\O$.

Assume now that $P_{\F}(t)=(m+2)t+m$ with $\frac{2n}{3}
< m \le n$. Since $\F$ is generated by global sections we must
have $h^0(\F') \le n-1$. Thus
\begin{align*}
n-1 \ge m+h^1(\F') > \frac{2n}{3}+h^1(\F'), \quad
\frac{n}{3}-1 > h^1(\F') \quad \text{forcing} \quad
h^1(\F')=0.
\end{align*}
The Beilinson sequence (4.1) of $\F'$ gives one of the following
resolutions:
\begin{align*}
0 \lra 2\O(-2) \oplus (m-2)\O(-1) \lra m\O \lra \F' \lra 0,
\end{align*}
\begin{align*}
0 \lra 2\O(-2) \oplus (m-1)\O(-1) \stackrel{\psi}{\lra}
\O(-1) \oplus m\O \lra \F' \lra 0,
\end{align*}
\begin{align*}
0 \lra 2\O(-2) \oplus m\O(-1) \stackrel{\psi}{\lra}
2\O(-1) \oplus m\O \lra \F' \lra 0,
\end{align*}
with $\psi_{12}=0$.
In the third case $\F'$ has a subsheaf $\F''$ with resolution
\begin{align*}
0 \lra m\O(-1) \lra m\O \lra \F'' \lra 0.
\end{align*}
This situation has been examined before. In the first case
we get an exact commutative diagram
\begin{displaymath}
\xymatrix
{
0 \ar[r] & 2\O(-2) \oplus (m-2)\O(-1)
\ar[r]^{\ \ \ \ \ \ \ \ \ \ \ \ \ \psi}
\ar[d]^{\b} & m\O \ar[r] \ar[d]^{\a} & \F' \ar[r] \ar[d] & 0 \\
0 \ar[r] & 3\O(-2) \oplus (n-3)\O(-1)
\ar[r]^{\ \ \ \ \ \ \ \ \ \ \ \ \ \ \ \f}
& n\O \ar[r] & \F \ar[r] & 0
}
\end{displaymath}
with $\a, \ \b$ injective. We have $m-2 \le n-3$ because $\b_{22}$
is injective. On the other hand, $m\ge n-1$ by hypothesis.
Thus $\b_{22}$ is an isomorphism forcing rank$(\b_{11})=2$.
In consequence
\begin{align*}
\f \sim \left[
\begin{array}{cc}
\star & \psi \\
\star & 0
\end{array}
\right] .
\end{align*}
Finally, assume that $\F'$ has the second resolution.
We get the commutative diagram
\begin{displaymath}
\xymatrix
{
0 \ar[r] & 2\O(-2) \oplus (m-1)\O(-1) \ar[r]^{\ \ \ \ \ \ \ \psi}
\ar[d]^{\b} & \O(-1) \oplus m\O \ar[r] \ar[d]^{\a} & \F' \ar[r] \ar[d] & 0 \\
0 \ar[r] & 3\O(-2) \oplus (n-3)\O(-1) \ar[r]^{\ \ \ \ \ \ \ \ \ \ \ \f}
& n\O \ar[r] & \F \ar[r] & 0
} .
\end{displaymath}
The map $\a_{12}$ is injective because $\a$ is injective on
global sections. We have $\f_{12}\b_{22}= \a_{12}\psi_{22}$.
The latter map is injective, hence $\b_{22}$ is injective, too.
Thus $n-3 \ge m-1$, $n \ge m+2 > \frac{2n}{3} +2$, $n>6$,
contradiction.

The rest of the proof is as in (4.3), so it will be omitted.

\noindent \\
{\bf (4.7) Claim:} \emph{Let $\F$ be a sheaf on $\P^2=\P(V)$ with
resolution (ii) and $3\le n \le 6$. Then $\F$ is semistable if and
only if the entries of $\f_{11}$ span a subspace of $V^*$ of
dimension at least two and $\f$ is not equivalent to a matrix
of the form}
\begin{align*}
\left[
\begin{array}{cc}
\star & \psi \\
\star & 0
\end{array}
\right]
\end{align*}
\emph{where}
\begin{eqnarray*}
& & \psi : m\O(-1) \lra m\O, \qquad 1 \le m \le n-2, \quad
                                \text{\emph{or}} \\
& & \psi : \O(-2)\oplus (m-1)\O(-1) \lra m\O, \qquad
\frac{n}{3} < m \le n-1, \quad \text{\emph{or}} \\
& & \psi : 2\O(-2) \oplus (m-1)\O(-1) \lra \O(-1) \oplus m\O,
\qquad \frac{2n}{3} < m \le n-1.
\end{eqnarray*}
\emph{Equivalently, $\F$ is semistable if and only if $\f$ is
semistable with respect to $\L$. Here $\L=(\l_1,\l_2,\m_1,\m_2)$
is a polarization satisfying the property that the pair $(\l_1,\m_1)$}\\
(i) \emph{is in the interior of the segment with endpoints
$(\frac{1}{4},\frac{1}{4})$ and $(\frac{1}{5},\frac{2}{5})$,
in the case $n=3$;} \\
(ii) \emph{is in the interior of the triangle bounded by the lines
$\m_1 = \l_1,\ \m_1 = 1-4\l_1,\ \m_1 = 4\l_1 - \frac{1}{3}$,
in the case $n=4$;} \\
(iii) \emph{is in the interior of the triangle bounded by the lines
$\m_1 = \l_1,\ \m_1 = \frac{1}{6},\ \m_1 = \frac{15}{4}\l_1 -
\frac{1}{4}$, in the case $n=5$;} \\
(iv) \emph{is in the interior of the segment with endpoints
$(\frac{1}{7},\frac{1}{7})$ and $(\frac{3}{29},\frac{5}{29})$,
in the case $n=6$.}

\noindent \\
\emph{Proof:} The proof is the same as the proof of (4.6) with the
only difference that when $n\neq 3$ it is possible to have
$h^1(\F')=1$ and $P_{\F'}(t)=(m+2)t+m$. For such a sheaf the
Beilinson sequence takes the form
\begin{align*}
0 \lra 2\O(-2) \oplus (p+m-2)\O(-1) \stackrel{\rho}{\lra}
p\O(-1) \oplus (m+1)\O \stackrel{\eta}{\lra} \O \lra 0.
\end{align*}
As $\r_{12}=0$ we may assume that
$\eta = (X,Y,Z,0,\ldots,0)$.
Thus $p \ge 3$ and $\F'$ has resolution 
\begin{align*}
0 \ra 2\O(-2) \oplus (p+m-2)\O(-1) \lra \Om^1 \oplus (p-3)\O(-1)
\oplus (m+1)\O \ra \F' \ra 0
\end{align*}
from which we get the resolution
\begin{align*}
0 \ra \O(-3) \oplus 2\O(-2) \oplus (p+m-2)\O(-1)
\stackrel{\f'}{\lra} 3\O(-2) \oplus (p-3)\O(-1) \oplus (m+1)\O
\end{align*}
\begin{align*}
\lra \F' \lra 0
\end{align*}
with $\f_{13}'=0$ and $\f_{23}'=0$.
Since $\f'$ is injective we must have $p=3$.
But then $\F'$ has a subsheaf $\F''$ with resolution
\begin{align*}
0 \lra (m+1)\O(-1) \lra (m+1)\O \lra \F'' \lra 0.
\end{align*}
This situation has already been examined. Q.e.d.

\noindent \\
{\bf (4.8) Claim:} \emph{Let $\F$ be a sheaf on $\P^2$ with
Hilbert polynomial $P_{\F}(t)=6t+3$ and resolution}
\begin{align*}
0 \lra 3\O(-2) \oplus 2 \O(-1) \stackrel{\f}{\lra}
2\O(-1) \oplus 3\O \lra \F \lra 0,
\end{align*}
\emph{with $\f_{12}=0$. Then $\F$ is semistable if and only if}
\begin{align*}
\f_{11} \nsim \left[
\begin{array}{ccc}
\star & \star & 0 \\
\star & \star & 0
\end{array}
\right], \quad \f_{11} \nsim \left[
\begin{array}{ccc}
\star & 0 & 0 \\
\star & \star & \star
\end{array}
\right], \quad \f_{22} \nsim \left[
\begin{array}{cc}
\star & 0 \\
\star & 0 \\
\star & \star
\end{array}
\right], \quad \f_{22} \nsim \left[
\begin{array}{cc}
0 & 0 \\
\star & \star \\
\star & \star
\end{array}
\right] .
\end{align*}
\emph{Equivalently, $\F$ is semistable if and only if $\f$ is
semistable with respect to $\L$, where $\L =(\l_1,\l_2,\m_1,\m_2)$
is any polarization for which $(\l_1,\m_1)$ is in the interior
of the triangle with vertices $(0,0),\, (\frac{1}{5},\frac{1}{5}),\,
(\frac{1}{3},\frac{1}{2})$.}

\noindent \\
{\bf (4.9) Claim:} \emph{Let $\F$ be a sheaf on $\P^2$ with
Hilbert polynomial $P_{\F}(t)=7t+4$ and resolution}
\begin{align*}
0 \lra 3\O(-2) \oplus 3\O(-1) \lra 2\O(-1) \oplus 4\O \lra \F
\lra 0,
\end{align*}
\emph{with $\f_{12}=0$. Then $\F$ is semistable if and only if
$\f_{11}$ satisfies the same conditions as in (4.8) and, in
addition, $\f_{22}$ is not equivalent to a matrix of the form}
\begin{align*}
\left[
\begin{array}{cc}
\star & \psi \\
\star & 0
\end{array}
\right] \qquad \text{\emph{with}} \quad
\psi : m\O(-1) \lra m\O, \quad m=1,2,3.
\end{align*}
\emph{Equivalently, $\F$ is semistable if and only if $\f$ is
semistable with respect to $\L$, where $\L =(\l_1,\l_2,\m_1,\m_2)$
is any polarization for which $(\l_1,\m_1)$ is in the interior
of the quadrilater with vertices}
\begin{align*}
(0,\, 0), \quad \left( \frac{1}{3},\, \frac{1}{2} \right), \quad
\left( \frac{17}{24},\, 1 \right), \quad (1,\, 1).
\end{align*}


\section{Sheaves $\F$ with $h^0(\F(-1))=0$ and $h^1(\F)=1$}

In this section $\F$ will be a sheaf on $\P^2$ with $h^0(\F(-1))=0$,
$h^1(\F)=1$ and Hilbert polynomial $P_{\F}(t)=at+b$, $0 \le b < a$.
From the Beilinson complex (4.1) we deduce that $\F$ must have one
of the following resolutions:\\
(i) when $a\le 2b$
\begin{align*}
0 \ra (a-b)\O(-2)\oplus (m+2b-a)\O(-1) \stackrel{\r}{\lra}
\Om^1 \oplus (m-3)\O(-1) \oplus (b+1)\O \ra \F \ra 0,
\end{align*}
(ii) when $a>2b$
\begin{align*}
0 \ra (a-b)\O(-2)\oplus m\O(-1) \stackrel{\r}{\lra}
\Om^1 \oplus (m+a-2b-3)\O(-1) \oplus (b+1)\O \ra \F \ra 0,
\end{align*}
where $m$ is an integer and $\r_{12}=0,\ \r_{22}=0$. Combining these
with the exact sequence
\begin{align*}
0 \lra \O(-3) \lra 3\O(-2) \lra \Om^1 \lra 0
\end{align*}
we get one of the following resolutions for $\F$:
\begin{align*}
\tag{i}
0 \lra \O(-3)\oplus (a-b)\O(-2)\oplus (m+2b-a)\O(-1)
\qquad \qquad \qquad \qquad
\end{align*}
\begin{align*}
\qquad \qquad \qquad \qquad
\stackrel{\psi}{\lra}
3\O(-2) \oplus (m-3)\O(-1) \oplus (b+1)\O \lra \F \lra 0,
\end{align*}
\begin{align*}
\tag{ii}
0 \lra \O(-3)\oplus (a-b)\O(-2)\oplus m\O(-1) \qquad \qquad \qquad \qquad
\end{align*}
\begin{align*}
\qquad \qquad \qquad \qquad 
\stackrel{\psi}{\lra}
3\O(-2) \oplus (m+a-2b-3)\O(-1) \oplus (b+1)\O \lra \F \lra 0,
\end{align*}
where
\begin{align*}
\psi = \left[
\begin{array}{ccc}
X \\
Y & \psi_{12} & 0 \\
Z \\
0 & \r_{21} & 0 \\
0 & \r_{31} & \r_{32}
\end{array}
\right].
\end{align*}

\noindent \\
{\bf (5.1) Claim:} \emph{There are no semistable sheaves $\F$
on $\P^2$ with $h^0(\F(-1))=0$, $h^1(\F)=1$ and Hilbert polynomial}
\begin{align*}
P_{\F}(t)=(n+1)t+n \quad \text{or} \quad
P_{\F}(t)=(n+2)t+n, \quad n \ge 0.
\end{align*}

\noindent \\
\emph{Proof:} If $a-b=1$ then $\psi$ cannot be injective.
If $a-b=2$, then we must have $m=3$ in case (i), respectively
$m=2b+3-a$ in case (ii). It follows that $\F$ has a subsheaf
$\F'$ with resolution
\begin{align*}
0 \lra (n+1)\O(-1) \stackrel{\r_{32}}{\lra} (n+1)\O \lra \F'
\lra 0.
\end{align*}
This subsheaf destabilizes $\F$.

If $P_{\F}(t)=(n+1)t+n$, the claim already follows from
the fact that $h^0(\F(-1))=0$ implies $h^1(\F)=0$, cf.
(2.3)(v). Above we have an alternate argument.
\footnote{The referee pointed out that for any sheaf on
$\P^2$, semistable or not, with Hilbert polynomial
$P_{\F}(t)=(n+1)t+n$ and $h^0(\F(-1))=0$, we must have
$h^1(\F)=0$. Indeed, as in the proof of (4.3) with $\F$
instead of $\F'$, we have a monad
\bdm
0 \lra (n+2)\O(-1) \lra (n+3)\O \stackrel{\eta}{\lra} \O(1)
\lra 0
\edm
with cohomology $\F$. As $(n+3)\O$ is generated by global
sections, it follows that $\O(1)$ is generated by their
images under $\eta$. But $\O(1)$ cannot be generated by
fewer than 3 linearly independent sections. Thus
$\eta$ is surjective on the level of global sections.
We get $h^1({\mathcal Ker}(\eta))=0$ and, a fortiori,
$h^1 (\F)=0$.}

In the remaining part of this section we will assume that $\F$
is a semistable sheaf with Hilbert polynomial
$P_{\F}(t)=(n+3)t+n,\ n \ge 0$.
We have the resolution
\begin{align*}
0 \ra \O(-3) \oplus 3\O(-2) \oplus (m+n)\O(-1)
\stackrel{\psi}{\lra}
3\O(-2) \oplus m\O(-1) \oplus (n+1)\O \ra \F \ra 0.
\end{align*}
We must have $m\le 1$ to ensure that $\psi$ is injective.
But if $m=1$ then $\F$ has a destabilizing sheaf $\F'$ as above.
Thus $m=0$ and we arrive at the resolution
\begin{align*}
0 \lra \O(-3) \oplus 3\O(-2) \oplus n\O(-1)
\stackrel{\psi}{\lra}
3\O(-2) \oplus (n+1)\O \lra \F \lra 0.
\end{align*}
Since $\psi$ is injective we must have rank$(\psi_{12})\ge 2$.
If rank$(\psi_{12})=2$ then $\F$ has a destabilizing subsheaf
$\F'$ with resolution
\begin{align*}
0 \lra \O(-2) \oplus n\O(-1) \lra (n+1)\O \lra \F' \lra 0.
\end{align*}
In conclusion rank$(\psi_{12})=3$ and $\F$ has the resolution
\begin{align*}
0 \lra \O(-3) \oplus n\O(-1) \stackrel{\f}{\lra} (n+1)\O
\lra \F \lra 0.
\end{align*}

If $n>3$ some of the semistability conditions on $\F$ cannot be translated
into semistability conditions on $\f$ because one of the conditions
on $\f$ would have to be that there is no commutative exact
diagram
\begin{displaymath}
(*) \quad \quad
\xymatrix
{
0 \ar[r] & \O(-2) \oplus (m-1)\O(-1) \ar[r]^{\ \ \ \ \ \ \ \ \ \ \ \f'}
\ar[d]^{\b} & m\O \ar[r] \ar[d]^{\a} & \F' \ar[r] \ar[d] & 0 \\
0 \ar[r] & \O(-3) \oplus n\O(-1) \ar[r]^{\ \ \ \ \ \ \ \ \ \ \ \f}
& (n+1)\O \ar[r] & \F \ar[r] & 0
}
\end{displaymath}
with $\a,\ \b$ injective and $3m>n$.
If $n>3$ then $\b$ may have the form
\begin{align*}
\left[
\begin{array}{c}
0 \\
\b_0
\end{array}
\right]
\end{align*}
with $\b_0$ as in (4.4). In this case the condition $\f \b =\a \f'$
cannot be translated in terms of semistability of $\f$.

\noindent \\
{\bf (5.2) Claim:} \emph{Let $\F$ be a sheaf on $\P^2$ with
$h^0(\F(-1))=0,\ h^1(\F)=1$ and Hilbert polynomial}
\begin{align*}
P_{\F}(t)=(n+3)t+n, \qquad n=1,2,3.
\end{align*}
\emph{Then $\F$ is semistable if and only if it has a resolution}
\begin{align*}
0 \lra \O(-3)\oplus n\O(-1) \stackrel{\f}{\lra}
(n+1)\O \lra \F \lra 0
\end{align*}
\emph{with $\f$ not equivalent to a matrix of the form}
\begin{align*}
\left[
\begin{array}{cc}
\star & \psi \\
\star & 0
\end{array}
\right] \qquad \text{\emph{where}} \quad
\psi : m\O(-1) \lra m\O, \quad 1\le m \le n.
\end{align*}
\emph{The morphisms $\f$ occuring above are precisely those
morphisms semistable with respect to $\L$ with nonzero
determinant. Here $\L=(\l_1,\l_2,\m_1)$ is any polarization
satisfying $0 < \l_1 < \frac{1}{n+1}$. If $n>3$ solely the
``only if'' part of the above statement remains true.}

\noindent \\
\emph{Proof:} One direction follows from the discussion before
the claim. Conversely, we assume that $\F$ has a resolution as above
and we try to prove that $\F$ is semistable. As in the proof of
(4.6), a destabilizing subsheaf $\F'$ of $\F$ must have one
of the following Hilbert polynomials: $mt+m$ with $1 \le m \le n$,
$(m+1)t+m$ with $\frac{n}{3} < m \le n$, $(m+2)t+m$ with
$\frac{2n}{3} < m \le n$. In the first case we get the contradiction
\begin{align*}
\f \sim \left[
\begin{array}{cc}
\star & \psi \\
\star & 0
\end{array}
\right] \qquad \text{with} \quad
\psi : m\O(-1) \lra m\O, \quad 1\le m \le n.
\end{align*}
In the second case we have the exact commutative diagram (*)
from above with injective $\a$ and $\b$. Since ${\mathcal Coker}
(\b)$ is torsion-free as a subsheaf of the torsion-free
sheaf ${\mathcal Coker}(\a)$, and since
\begin{align*}
\b \nsim \left[
\begin{array}{c}
0 \\
\b_0
\end{array}
\right], \qquad \text{we must have} \quad
\b = \left[
\begin{array}{c}
0 \\
X \\
Y
\end{array}
\right] \quad \text{or} \quad \b = \left[
\begin{array}{cc}
0 & 0 \\
X & 0 \\
Y & 0 \\
0 & 1
\end{array}
\right] .
\end{align*}
We get
\begin{align*}
\f \sim \left[
\begin{array}{cc}
\star & \psi \\
\star & 0
\end{array}
\right] \qquad \text{with} \quad
\psi : (m+1)\O(-1) \lra (m+1)\O.
\end{align*}
Finally, in the third case, we must have $m=n$ and, as in the
proof of (4.6), an exact commutative diagram
\begin{displaymath}
\xymatrix
{
0 \ar[r] & 2\O(-2) \oplus (n-2)\O(-1) \ar[r]^{\ \ \ \ \ \ \ \ \ \ \ \f'}
\ar[d]^{\b} & n\O \ar[r] \ar[d]^{\a} & \F' \ar[r] \ar[d] & 0 \\
0 \ar[r] & \O(-3) \oplus n\O(-1) \ar[r]^{\ \ \ \ \ \ \ \ \f}
& (n+1)\O \ar[r] & \F \ar[r] & 0
}
\end{displaymath}
or a diagram
\begin{displaymath}
\xymatrix
{
0 \ar[r] & 2\O(-2) \oplus (n-1)\O(-1) \ar[r]^{\ \ \ \ \ \ \ \ \f'}
\ar[d]^{\b} & \O(-1) \oplus n\O \ar[r]
\ar[d]^{\a} & \F' \ar[r] \ar[d] & 0 \\
0 \ar[r] & \O(-3) \oplus n\O(-1) \ar[r]^{\ \ \ \ \ \ \ \ \f}
& (n+1)\O \ar[r] & \F \ar[r] & 0
} .
\end{displaymath}
In the first case $\a$ is injective because it is injective on
global sections. Thus $\b$ is injective. But then ${\mathcal Coker}(\b)$
has a direct summand supported on a conic. This contradicts
the fact that ${\mathcal Coker}(\b)$ is a subsheaf of
${\mathcal Coker}(\a)\isom \O$.

In the second diagram we have $\f_{12}\b_{22}=\a_{12}\f_{22}'$.
But $\a_{12}$ is injective because $\a$ is injective on global
sections. Also, $\f_{22}'$ is injective because $\f_{12}'=0$
and $\f'$ is injective. Thus $\f_{12}\b_{22}$ is injective,
forcing $\b_{22}$ to be injective. It follows that ${\mathcal Ker}(\b)
\subset 2\O(-2)$. If $\a$ is not injective then we get the
contradiction
\begin{align*}
\O(-1) \isom {\mathcal Ker}(\a) \isom {\mathcal Ker}(\b)
\subset 2\O(-2).
\end{align*}
Thus $\a$ is injective forcing ${\mathcal Coker}(\a)$ to be
supported on a line. But this is impossible, because 
$\O(-3) \subset {\mathcal Coker}(\b) \subset {\mathcal Coker}(\a)$.
Q.e.d.

\noindent \\
{\bf (5.3) Remark:} The sheaves from (5.2) with Hilbert
polynomial $6t+3$ are stable. Indeed, assume that $\F$
has a subsheaf $\F'$ with Hilbert polynomial $2t+1$.
It must be stable, hence the structure sheaf of a conic.
We arrive at a commutative diagram
\bdm
\xymatrix
{
0 \ar[r] & \O(-2) \ar[r] \ar[d]^{\b} & \O \ar[r] \ar[d]^{\a}
& \F' \ar[r] \ar[d] & 0 \\
0 \ar[r] & \O(-3) \oplus 3\O(-1) \ar[r]^{\quad \quad \quad \f}
& 4\O \ar[r] & \F \ar[r] & 0
}
\edm
with injective $\a$ and $\b$. After performing column
operations on $\f$ we may assume that three among the
rows of $\f_{12}\b_{21}$ are zero. But, according to
(5.2), $\f_{12}$ is semistable with respect to the only
admissible polarization on the vector space of morphisms
$3\O(-1) \lra 4\O$. From remark (5.4) we get $\b_{21}=0$,
so $\b=0$, contradiction.

Assume now that $\F$ has a quotient sheaf $\F''=\F/\F'$ with
Hilbert polynomial $2t+1$. $\F''$ must be stable, hence
it is the structure sheaf of a conic, hence it is generated
by one global section. Thus the map $\F \lra \F''$
is surjective on global sections, forcing $h^0(\F')=2$.
Thus $h^1(\F')=0$ which, together with $h^1(\F'')=0$
implies that $h^1(\F)=0$. Contradiction.

\noindent \\
{\bf (5.4) Remark:} Let $\f$ be a $4\times 3$-matrix with
entries in $V^*$ which is semistable:
Modulo operations on rows and columns, $\f$ is not equivalent
to a matrix having a zero row, a zero $2\times 2$-submatrix,
or a zero $3\times 1$-submatrix. Then one of the maximal
minors of $\f$ is not zero.

\noindent \\
\emph{Proof:} Assume that all maximal minors of $\f$ are
zero. Each $3\times 3$-submatrix $\psi$ of $\f$ satisfies the
hypotheses of remark (5.5), hence it is equivalent to
$\psi_1$ or $\psi_2$. We can choose $\psi$ to have a zero entry,
thus ruling out $\psi_2$. From the assumption that all
minors of $\f$ are zero it is easy to deduce that the
row of $\f$ which is not part of $\psi$ is a linear
combination of the rows of $\psi$, cf. the proof of (6.7).
This contradicts the semistability of $\f$.

\noindent \\
{\bf (5.5) Remark:} Let $\psi$ be a $3\times 3$-matrix with
entries in $V^*$ and zero determinant. Assume that $\psi$ is
equivalent to neither of the following matrices:
\begin{align*}
\left[
\begin{array}{ccc}
0 & \star & \star \\
0 & \star & \star \\
0 & \star & \star
\end{array}
\right] \quad \text{or} \quad \left[
\begin{array}{ccc}
0 & 0 & \star \\
0 & 0 & \star \\
\star & \star & \star
\end{array}
\right] \quad \text{or} \quad \left[
\begin{array}{ccc}
0 & 0 & 0 \\
\star & \star & \star \\
\star & \star & \star
\end{array}
\right] .
\end{align*}
Then $\psi$ is equivalent to one of the following matrices:
\begin{align*}
\psi_1 = \left[
\begin{array}{ccc}
X & Y & 0 \\
Z & 0 & Y \\
0 & -Z & X
\end{array}
\right] \quad \text{or} \ \ \psi_2 = \left[
\begin{array}{lll}
X & Y & Z \\
Y & a_1 X + a_2 Y & a_3 X + a_4 Y + a_5 Z \\
Z & a_6 X + a_7 Y + a_8 Z & a_9 X + a_{10}Z
\end{array}
\right]
\end{align*}
with $a_1,\ldots,a_{10} \in k^*$. Modulo operations on rows and
columns $\psi_2$ is not equivalent to a matrix having a zero
entry.

\noindent \\
\emph{Proof:} We distinguish two cases:
either $\psi$ has one zero entry
or, modulo equivalence,
all entries of $\psi$ are nonzero. In the second case
$\psi \sim \psi_2$. In the first case we may assume that
$\psi_{11}=X,\ \psi_{12}=Y,\ \psi_{13}=0$. We now consider two
subcases: span$\{ \psi_{23},\psi_{33} \}$ is equal to or is
different from span$\{ X,Y \}$. In the second subcase we may write
\begin{align*}
\psi = \left[
\begin{array}{ccc}
X & Y & 0 \\
\psi_{21} & \psi_{22} & Y \\
\psi_{31} & \psi_{32} & Z
\end{array}
\right] .
\end{align*}
We have
\begin{align*}
\text{det}(\psi)= XZ\psi_{22}+ Y^2 \psi_{31} -XY\psi_{32}
      -YZ \psi_{21} = 0
\end{align*}
forcing $\psi_{22}=aY$. Performing operations on rows we may
assume that $\psi_{22}=0$. Thus $Y\psi_{31}-X\psi_{32}
-Z\psi_{21}=0$. We get
\begin{align*}
\psi_{31}=cZ \quad \text{modulo}\ \ X, \qquad
\psi_{21}=cY \quad \text{modulo}\ \ X.
\end{align*}
But then
\begin{align*}
\psi \sim \left[
\begin{array}{ccc}
X & Y & 0 \\
0 & \star & Y \\
0 & \star & Z
\end{array}
\right] .
\end{align*}
From det$(\psi)=0$ we get
\begin{align*}
\psi \sim \left[
\begin{array}{ccc}
X & Y & 0 \\
0 & 0 & Y \\
0 & 0 & Z
\end{array}
\right] ,
\end{align*}
contradiction. This eliminates the second subcase.

Finally, we may assume that
\begin{align*}
\psi = \left[
\begin{array}{ccc}
X & Y & 0 \\
\psi_{21} & \psi_{22} & Y \\
\psi_{31} & \psi_{32} & X
\end{array}
\right] .
\end{align*}
We have
\begin{align*}
\text{det}(\psi)=X^2\psi_{22}+Y^2 \psi_{31}-XY\psi_{32}
                      -XY\psi_{21}=0
\end{align*}
hence
\begin{align*}
\psi_{22}=aY,\quad \psi_{31}=bX,\quad \psi_{21}+\psi_{32}
     =aX+bY.
\end{align*}
Performing operations on rows we may assume that $a=b=0$.
Denoting $Z=\psi_{21}$ we arrive at $\psi\sim \psi_1$. Q.e.d.


\section{Sheaves $\F$ with $h^0(\F(-1))\neq 0$ and $h^1(\F)=0$}

Let $\F$ be a sheaf on $\P^2$ with $h^0(\F(-1))=p\neq 0$,
$h^1(\F)=0$ and Hilbert polynomial $P_{\F}(t)=at+b,\ 0 \le b < a$.
From the Beilinson complex we deduce that $\F$ has to have one
of the following resolutions: \\
(i) when $a<2b$
\begin{align*}
0 \ra p\O(-2) \stackrel{\psi}{\lra} (p+a-b)\O(-2)\oplus (m+2b-a)\O(-1)
\stackrel{\f}{\lra} m\O(-1) \oplus b\O \ra \F \ra 0,
\end{align*}
(ii) when $a \ge 2b$
\begin{align*}
0 \ra p\O(-2) \stackrel{\psi}{\lra} (p+a-b)\O(-2)\oplus m \O(-1)
\stackrel{\f}{\lra} (m+a-2b)\O(-1) \oplus b\O \ra \F \ra 0,
\end{align*}
where $m$ is an integer, $\psi_{11}=0$ and $\f_{12}=0$. In case (ii) we
must have $m\ge 2$ because, if $m=1$, then we get the contradiction
$\psi = 0$. We obtain the following exact commutative diagram, say
in case (ii):
\begin{displaymath}
\xymatrix
{
& & 0 \ar[d] & 0 \ar[d] \\
0 \ar[r] & p\O(-2) \ar[d]^{=} \ar[r]^{\psi_{21}} & m\O(-1)
\ar[d]^{\left[ \begin{array}{c} 0 \\ I \end{array} \right]}
\ar[r]^{\f_{22}} & b\O \ar[r] \ar[d]^{\left[ \begin{array}{c} 0 \\ I
\end{array} \right]} & \C \ar[r] \ar[d] & 0 \\
0 \ar[r] & p\O(-2) \ar[d] \ar[r]^{\psi \qquad \qquad }
& (p+a-b)\O(-2) \oplus
m\O(-1) \ar[r]^{\f} \ar[d]^{[I,\, 0]} & (m+a-2b)\O(-1)\oplus b\O
\ar[r]^{\qquad \qquad \eta} \ar[d]^{[I,\, 0]} & \F \ar[d] \ar[r] & 0 \\
0 \ar[r] & \K \ar[r] & (p+a-b)\O(-2) \ar[d] \ar[r]^{\f_{11}} &
(m+a-2b)\O(-1) \ar[d] \ar[r] & \G \ar[d] \ar[r] & 0 \\
& & 0 & 0 & 0
} .
\end{displaymath}
The above induces the exact sequence
\begin{align*}
\tag{6.1}
0 \lra \K \lra \C \lra \F \lra \G \lra 0.
\end{align*}
Note that $\K$ is torsion-free or zero,
as a subsheaf of the torsion-free sheaf $(p+a-b)\O(-2)$.

\noindent \\
{\bf (6.2) Remark:} Assume that $\F$ is semistable. Then $\C$
does not have zero-dimensional torsion and is not supported on
a curve.

\noindent \\
\emph{Proof:} Let $\T$ be the zero-dimensional torsion of $\C$.
As $\F$ has no zero-dimensional torsion it follows that the
induced map $\T \lra \F$ is zero. Thus $\T$ is a subheaf
of $\K$. The latter is torsion-free, so $\T=0$.

Assume now that $\C$ is supported on a curve.
Then $m=p+b$ and $\K=0$. Thus $\C$ is a subsheaf of $\F$.
We have
\begin{align*}
P_{\C} (t)=b\, {t+2 \choose 2} -(p+b){t+1 \choose 2}+
p {t \choose 2} = (b-p)t+b.
\end{align*}
But $\frac{b}{b-p} > \frac{b}{a}$ which shows that $\C$ violates
the semistability of $\F$.

\noindent \\
{\bf (6.3) Remark:} Assume that $\F$ is semistable. Then, in
case (i), we either have $m+2b-a < b$ or $m+2b-a \ge b$ and
all maximal minors of $\f_{22}$ are zero. Similarly, in case (ii),
either $m<b$ or $m\ge b$ and all maximal minors of $\f_{22}$ are zero.
This follows from (6.2).

\noindent \\
{\bf (6.4) Remark:} $\eta$ is an isomorphism on global sections.
As a consequence, if $\F$ is semistable, then $\f_{22}$ cannot
have the form
\begin{align*}
\left[
\begin{array}{ccccc}
0 & 0 & \star & \cdots & \star \\
\vdots & \vdots & \vdots & & \vdots \\
0 & 0 & \star & \cdots & \star \\
X & Y & \star & \cdots & \star
\end{array}
\right] .
\end{align*}
Indeed, if $\f_{22}$ had the above form, then we would
get the commutative diagram
\begin{displaymath}
\xymatrix
{
2\O(-1) \ar[r]^{\ \ [X,\, Y]}
\ar[d]^{\left[ \begin{array}{c} I \\ 0 \end{array} \right]}
& \O \ar[r]
\ar[d]^{\left[ \begin{array}{c} 0 \\ I \end{array} \right]}
& {\mathbb C}_{x} \ar[d] \ar[r] & 0 \\
m\O(-1) \ar[r]^{\f_{22}} & b\O \ar[r]^{\eta_{12}} & \F
} .
\end{displaymath}
Here ${\mathbb C}_{x}$ is the structure sheaf of the point
$x=(0:0:1)$. But the map ${\mathbb C}_{x} \lra \F$ is zero
because $\F$ does not have zero-dimensional torsion.
This shows that $\eta_{12}$ has nontrivial kernel.
This contradicts the fact that $\eta$ is an isomorphism on
global sections.

\noindent \\
{\bf (6.5) Claim:} \emph{There are no semistable sheaves
$\F$ on $\P^2$ with $h^0(\F(-1))\neq 0$, $ h^1(\F)=0$ and
Hilbert polynomial}
\begin{align*}
P_{\F}(t)=nt+1, \quad n \ge 2 \quad \text{\emph{or}} \quad
P_{\F}(t)=nt+2, \quad n \ge 4.
\end{align*}

\noindent \\
\emph{Proof:} The case $P_{\F}(t)=nt+1$ follows directly from
(6.4) because $\f_{22}$ must have the form
\begin{align*}
[ X,\, Y,\, \star, \cdots, \star].
\end{align*}
In the case $P_{\F}(t)= nt+2$ all $2\times 2$-minors of $\f_{22}$
are zero, cf. (6.3). It follows that $\f_{22}$ has the form
\begin{align*}
\left[
\begin{array}{ccccc}
0 & 0 & \star & \cdots & \star \\
X & Y & \star & \cdots & \star
\end{array}
\right] .
\end{align*}
The claim follows from (6.4).

\noindent \\
{\bf (6.6) Remark:} Let $\a=(\a_{ij})$ be a morphism of sheaves on 
$\P^n = \P (V)$:
\begin{align*}
\a : (m+1)\O \lra m\O(l).
\end{align*}
Assume that at least one of the maximal minors of $\a$ is a
nonzero polynomial. Then ${\mathcal Ker}(\a)\isom \O(-d)$
where $d$ is an integer satisfying $0 \le d \le ml$.
More precisely, let $\a_i$, $1\le i \le m+1$, denote the
minor obtained from $\a$ by erasing the $i^{\text{th}}$ column.
Let
\begin{align*}
\b = (\b_1,\ldots , \b_{m+1}), \qquad \text{where} \quad
\b_i = \frac{\a_i}{\text{g.c.d.}(\a_1,\ldots,\a_{m+1})},
\quad 1 \le i \le m+1.
\end{align*}
Let $d$ be the degree of the entries of $\b$. Then we have
the exact sequence
\begin{align*}
0 \lra \O(-d) \stackrel{\b}{\lra} (m+1)\O \stackrel{\a}{\lra} m\O(l).
\end{align*}

\noindent \\
{\bf (6.7) Claim:} \emph{Let $\F$ be a semistable sheaf on
$\P^2=\P(V)$ with $h^0(\F(-1))\neq 0$, $h^1(\F)=0$ and Hilbert
polynomial $P_{\F}(t)= nt+3$, $n \ge 4$. Then $h^0(\F(-1))=1$
and $\F$ has a resolution}
\begin{align*}
0 \lra \O(-2) \stackrel{\psi}{\lra} (n-2)\O(-2)\oplus 3\O(-1)
\stackrel{\f}{\lra} (n-3)\O(-1) \oplus 3\O \lra \F \lra 0
\end{align*}
\emph{with $\f_{12}=0,\ \psi_{11}=0,\ \f_{21}\neq 0$,}
\begin{align*}
\psi_{21} \sim \left[
\begin{array}{c}
X \\
Y \\
Z
\end{array}
\right], \qquad \f_{22} \sim \left[
\begin{array}{rrr}
-Y & X & 0 \\
-Z & 0 & X \\
0 & -Z & Y
\end{array}
\right] , \qquad \f_{11} \nsim \left[
\begin{array}{cc}
\f' & 0 \\
\star & \star
\end{array}
\right]
\end{align*}
\emph{where $\f'$ is an $m\times m$-matrix with entries in $V^*$,
$1 \le m \le n-3$. Moreover, $\F$ is an extension of the form}
\begin{align*}
0 \lra \O_{C}(1) \lra \F \lra \G \lra 0
\end{align*}
\emph{where $C$ is a curve of degree $d$, $4\le d \le n$, and the
map $\F \lra \G$ is zero on global sections. If $n\ge 7$ then
$d \ge 5$.}

\noindent \\
\emph{Proof:} Assume $n \ge 6$ so that we are in case (ii).
If $m \ge 5$ then $\f_{22}$ has the form
\begin{align*}
\left[
\begin{array}{cccccc}
\star & \star & \star & 0 & \cdots & 0 \\
\star & \star & \star & \star & \cdots & \star \\
\star & \star & \star & \star & \cdots & \star
\end{array}
\right] = \left[
\begin{array}{cccccc}
\star & \star & \star      \\
\star & \star & \star & \f'\\
\star & \star & \star
\end{array}
\right] .
\end{align*}
From (6.3) we see that all $2\times 2$-minors of $\f'$
are zero. Since $\f_{22}$ cannot have a zero column it
follows that
\begin{align*}
\f_{22} \sim \left[
\begin{array}{cccccc}
\star & \star & \star & 0 & \cdots & 0 \\
\star & \star & \star & 0 & \cdots & 0 \\
\star & \star & \star & \star & \cdots & \star
\end{array}
\right] .
\end{align*}
By virtue of (6.4) this is impossible.

Assume now that $m=4$. Firstly, we notice that $\f_{22}$
cannot have a zero row because, if
\begin{align*}
\f_{22} = \left[
\begin{array}{ccc}
0 & \cdots & 0 \\
& \f' 
\end{array}
\right] ,
\end{align*}
then, arguing as at (6.2), we get that $\f'$ has all maximal
minors equal to zero hence $\f'$ has one row identically zero.
This, again, contradicts (6.4). Secondly, using the same
kind of arguments, we notice that $\f_{22}$ cannot have the
form
\begin{align*}
\left[
\begin{array}{cccc}
X & 0 & 0 & 0 \\
\star & \star & \star & \star \\
\star & \star & \star & \star
\end{array}
\right] .
\end{align*}
Now $\f_{22}$ has nontrivial kernel in $\oplus_4 V^*$
by hypothesis. No element in the kernel can have the form
\begin{align*}
\left[
\begin{array}{c}
X \\ Y \\ 0 \\ 0
\end{array}
\right]
\end{align*}
otherwise we would arrive at a matrix excluded by (6.4):
\begin{align*}
\f_{22} \sim \left[
\begin{array}{rrrr}
0 & 0 & \star & \star \\
0 & 0 & \star & \star \\
-Y & X & \star & \star
\end{array}
\right] .
\end{align*}
Performing operations on the columns of $\f_{22}$
we may assume that
\begin{align*}
\left[
\begin{array}{c}
X \\ Y \\ Z \\ 0
\end{array}
\right]
\end{align*}
is in the kernel of $\f_{22}$. Performing operations
on the rows of $\f_{22}$ we may assume that
\begin{align*}
\f_{22} = \left[
\begin{array}{rrrr}
-Y & X & 0 & u \\
-Z & 0 & X & v \\
0 & -Z & Y & w
\end{array}
\right] .
\end{align*}
From
\begin{align*}
0 = \left|
\begin{array}{rrr}
-Y & X & u \\
-Z & 0 & v \\
0 & -Z & w
\end{array}
\right| = Z^2 u -YZv + XZw
\end{align*}
we get $Zu-Yv+Xw=0$ which shows that the third column
of $\f_{22}$ is a linear combination of the first two
columns. Thus $\f$ is equivalent to a matrix having
a zero column, contradiction.

The case $m=2$ is excluded by using (6.4). We conclude that
$m=3$ and, from what was said above, that we have
\begin{align*}
\psi_{21} \sim \left[
\begin{array}{c}
X \\ Y \\ Z
\end{array}
\right] , \qquad \qquad
\f_{22} \sim \left[
\begin{array}{rrr}
-Y & X & 0 \\
-Z & 0 & X \\
0 & -Z & Y
\end{array}
\right] .
\end{align*}

Thus far we have obtained the desired resolution of $\F$
in the cases $n\ge 6$. The cases $n=4$ and $n=5$ are completely
analogous. From our concrete description of $\f_{22}$ we
see that $\C \isom \O(1)$. Since $\F$ surjects onto $\G$,
the latter has support of dimension zero or one. Thus, at least
one of the maximal minors of $\f_{11}$ must be a nonzero
polynomial. We can apply (6.6) to conclude that
${\mathcal Ker}(\f_{11})$$\isom \O(-d+1)$ for some integer
$d \ge 3$. We have
\begin{align*}
P_{\G}(t)= (n-3){t+1 \choose 2} -(n-2){t \choose 2}
+ {t+3-d \choose 2} = (n-d)t+ \frac{(d-2)(d-3)}{2}.
\end{align*}
The sheaf $\G$ violates the semistability of $\F$
precisely when
\begin{align*}
\frac{(d-2)(d-3)}{2(n-d)} < \frac{3}{n}, \qquad
\text{i.e.} \quad
n(d-5)< -6.
\end{align*}
Thus, we cannot have $d=3$ and, if $d=4$, then $n \le 6$.
We conclude that $\F$ is an extension
\begin{align*}
0 \lra \O_{C}(1) \lra \F \lra \G \lra 0
\end{align*}
with deg$(C)=d \ge 4$, respectively deg$(C) \ge 5$ in the case
$n \ge 7$. Finally, we cannot have
\begin{align*}
\f_{11} \sim \left[
\begin{array}{cc}
\f' & 0 \\
\star & \star
\end{array}
\right]
\qquad \text{with} \quad
\f' : m\O(-2) \lra m\O(-1), \quad 1 \le m \le n-3.
\end{align*}
Indeed, if this were the case, then, since $\f_{11}$
has at least one nonzero maximal minor, we would get
det$(\f')\neq 0$ and a surjection $\F \lra {\mathcal Coker}(\f')$
onto a sheaf with Hilbert polynomial $P(t)=mt$.
This would contradict the semistability of $\F$. Q.e.d.

\noindent \\
{\bf (6.8) Lemma:} \emph{Let $C \subset \P^2$ be a curve given
by the equation $f=0$, where $f(X,Y,Z)$ is a homogeneos polynomial.
Let $\I \subset \O_C$ be a sheaf of ideals. Then there is
a homogeneous polynomial $g(X,Y,Z)$ dividing $f$ such that
the sheaf of ideals $\J \subset \O_C$ generated by $g$ satisfies:
$\I \subset \J$ and $\J / \I$ is supported on finitely many
points.}

\noindent \\
\emph{Proof:} Dehomogenizing in a suitable open affine subset
we reduce the problem to the following: let $f(X,Y)$
be a polynomial in $k[X,Y]$. Let $I \subset k[X,Y]$
be an ideal containing $f$. Then there is a polynomial
$g(X,Y)$ dividing $f$ such that $I \subset <g>$ and
$<g>/I$ is supported on finitely many points.

Let $f=f_1^{n_1} \cdot \ldots \cdot f_{\k}^{n_{\k}}$ be the
decomposition of $f$ into irreducible factors.
Let
\begin{align*}
I = \q_1 \cap \ldots \cap \q_m \cap \ag_1 \cap \ldots \cap \ag_l
\end{align*}
be a primary decomposition of $I$. Here $m \le \k$,
$\q_i$ is a primary ideal associated to $<f_i>$ and
$\ag_1, \ldots, \ag_l$ are primary ideals associated
to maximal ideals $\mg_1, \ldots, \mg_l$. Let us put
\begin{align*}
\q = \q_1 \cap \ldots \cap \q_m.
\end{align*}
We notice that $\q / I$ is supported on $\mg_1, \ldots,
\mg_l$. For $1 \le i \le m$ let $r_i$ be the largest integer
such that $\q_i \subset <f_i^{r_i}>$. We claim that
$g= f_1^{r_1} \cdot \ldots \cdot f_m^{r_m}$ is the desired
polynomial. To prove this it is enough to show that
$<g>/\q$ is supported on finitely many points. Since
localization commutes with intersections it is enough to
show that each $<f_i^{r_i}>/\q_i$ is supported on finitely
many points.

So far we have reduced the problem to the following:
let $f\in k[X,Y]$ be an irreducible polynomial.
Let $\q \subset k[X, Y]$ be a primary ideal associated
to $<f>$. Let $r \ge 1$ be the largest integer such
that $\q \subset <f^r>$. Then $<f^r>/\q$ is supported on
finitely many points.

We may assume that $\q$ is not a power of $<f>$.
Let $s$ be the smallest integer such that $<f^s>\subset \q$.
We will prove the above statement by induction on $s$.
If $s=r+1$ then
\begin{align*}
<f^r>/\q \isom \frac{<f^r>/<f^{r+1}>}{\q / <f^{r+1}>}
\end{align*}
can be regarded as the structure sheaf of a proper subscheme
of the scheme $X\subset \P^2$ given by $\{ f=0 \}$.
This is so because
\begin{align*}
<f^r>/<f^{r+1}> \isom k[X,Y]/<f> \qquad
\text{as} \quad k[X,Y] \text{- modules.}
\end{align*}
But $X$ is an irreducible scheme of dimension one, hence
any proper subscheme has dimension zero.

Assume now that $s > r+1$ and the statement is true for any
ideal $\q'$ satisfying $\q' \subsetneqq <f^r>,\ \ 
\q' \nsubseteq  <f^{r+1}>,\ \ <f^{s-1}> \subsetneqq \q'$.
Such an ideal is $\q'= \q +<f^{s-1}>$. By the induction
hypothesis we know that $<f^r>/\q'$ is supported on finitely
many points. To finish the proof it is enough to show that
$\q'/\q$ is supported on finitely many points.
But
\begin{align*}
\q'/\q \isom <f^{s-1}>/\q \ \cap <f^{s-1}>.
\end{align*}
If $\q \ \cap <f^{s-1}> \neq <f^s>$ then the right-hand side
is supported on finitely many points by the first step in
the induction argument. Let us now choose $h \in $ $\q \ \setminus
<f^{r+1}>$. Then $f^{s-r-1}h \in \q \ \cap <f^{s-1}> \setminus
<f^s>$. This finishes the proof of the lemma. \\

In the remaining part of this section we will seek more precise
information about the morphisms occuring in (6.7). 
For a start, let us assume that $\F$ is an arbitrary sheaf
having a resolution as in (6.7), and let us determine which
subsheaves $\F' \subset \F$ are destabilizing.
Let $\G'$ be the image of $\F'$ in $\G$ and let $\I(1)$
be the preimage of $\F'$ in $\O_C(1)$.
Here $\I$ is the ideal sheaf of a subscheme of $C$.
By (6.8) we can find a curve $C'\subset C$ such that the
ideal sheaf $\J$ of $C'$ contains $\I$ and $P_{\I(1)}(t)
= P_{\J(1)}(t)-c$, where $c$ is a nonnegative integer.
From the exact sequence
\begin{align*}
0 \lra \I(1) \lra \F' \lra \G' \lra 0
\end{align*}
we get
\begin{align*}
P_{\F'}(t) = P_{\I(1)}(t) + P_{\G'}(t)
= P_{\J(1)}(t)+ P_{\G'}(t) -c.
\end{align*}
Let us put $\k =$ deg$(C')$. We allow $\k =0$ for the case
$\J =\O_C$. From the exact sequence
\begin{align*}
0 \lra \O (-d+1) \lra \O(-\k +1)\lra \J(1) \lra 0
\end{align*}
we see that $h^0(\J(1))=0$ if $\k \ge 2$. But then
$h^0 (\I(1))=0$, forcing the map $H^0 (\F')\lra H^0 (\G')$
to be injective. Since the map $\F \lra \G$ is zero on global
sections we see that $h^0 (\F')=0$. It follows that $\F'$
does not violate the semistability of $\F$.

In the case $\k =0$ we have $P_{\F/\F'}(t)= c+ P_{\G/\G'}(t)$,
hence $\F'$ violates the semistability of $\F$ if
and only if $\a_1 (\G/\G') >0$ and
\begin{align*}
\frac{\a_0 (\G/\G')+c}{\a_1 (\G/\G')} < \frac{3}{n}.
\end{align*}
Assume now that $\k = 1$. We have
\begin{align*}
P_{\F/\F'}(t)= P_{\O_C (1)/\I (1)}(t) + P_{\G / \G'}(t)
= t+2+c+ P_{\G / \G'},
\end{align*}
hence $\F'$ violates the semistability if and only if
\begin{align*}
\frac{2+c+\a_0 (\G/\G')}{1+ \a_1 (\G / \G')} < \frac{3}{n}.
\end{align*}
Now the exact sequence
\begin{align*}
0 \lra \O_C (1) \lra \F \lra \G \lra 0
\end{align*}
together with the hypothesis $h^1 (\F(i)) = 0$ for $i \ge 0$
give $h^1 (\G(i))=0$ for $i \ge 0$. This, together with the
exact sequence
\begin{align*}
0 \lra \G' \lra \G \lra \G/\G' \lra 0
\end{align*}
yield $h^1 (\G/\G'(i))=0$ for $i \ge 0$. In particular
$\a_0 (\G/\G') =$$ h^0 (\G/\G') \ge 0$.
This eliminates the case $\a_1 (\G / \G')=0$ from above.
Let us summarize our findings so far:

\noindent \\
{\bf (6.9) Remark:} $\F$ is semistable if and only if there
are no quotients sheaves $\E$ of $\G$ satisfying
\begin{align*}
h^1(\E)=h^1(\E(1))=0
\qquad \text{and} \qquad
0 \le \a_0 (\E) < \frac{3}{n} \, \a_1(\E) \neq 0.
\end{align*}
One direction was proved in the discussion above.
The other direction follows by taking $\k=0$ and $c=0$,
in other words taking $\F'$ to be the preimage of $\G'$,
where $\G'$ is the kernel of the surjection $\G \lra \E$.

\noindent \\
{\bf (6.10) Claim:} \emph{Let $\F$ be a sheaf on $\P^2$ with
resolution}
\begin{align*}
0 \lra \O(-2) \lra (n-2)\O(-2)\oplus 3\O(-1)
\lra (n-3)\O(-1) \oplus 3\O \lra \F \lra 0
\end{align*}
\emph{satisfying the properties from (6.7). Assume that
$n \in \{ 4,5,6,7 \}$. Then $\F$ is semistable.}

\noindent \\
\emph{Proof:} Assume that there is $\E$ as in (6.9).
We must have $\a_0 (\E)=0$, otherwise
\begin{align*}
1 < \frac{3}{n} \, \a_1 (\E) \le \frac{3}{n} \, \a_1 (\G)
= \frac{3}{n} (n-d) \quad \text{forcing} \quad
\frac{3d}{2} < n.
\end{align*}
This gives $n \ge 8$, contradicting our hypothesis.

The Beilinson sequence of $\E(1)$ leads to the following
resolution:
\begin{align*}
0 \lra m\O(-2) \lra m\O(-1) \lra \E \lra 0
\end{align*}
for some integer $m \le n-d$. In the case $n=4$ there
is no such $\E$. In the case $n=5$ we have $m=1$.
In the cases $n \in \{ 6,7 \}$ we have $m \in \{ 1,2 \}$.
We obtain a commutative exact diagram
\begin{displaymath}
\xymatrix
{
& (n-2)\O(-2) \ar[r]^{\f_{11}} \ar[d]^{\b} & (n-3)\O(-1)
\ar[r] \ar[d]^{\a} & \G \ar[r] \ar[d] & 0 \\
0 \ar[r] & m\O(-2) \ar[r]^{\f'} & m\O(-1) \ar[r] & \E \ar[r]
\ar[d] & 0 \\
& & & 0
}
\end{displaymath}
with $\a \neq 0$ because the following composition is surjective:
\begin{align*}
(n-3)\O(-1) \stackrel{\a}{\lra} m\O(-1) \lra \E.
\end{align*}
After performing operations on the rows and columns of
$\f_{11}$ and $\f'$ it is always possible to write
\begin{align*}
\a = \left[
\begin{array}{rr}
I_r & 0 \\
0 & 0 
\end{array}
\right], \qquad \b = \left[
\begin{array}{rr}
I_s & 0 \\
0 & 0 
\end{array}
\right] .
\end{align*}
We arrive at
\begin{align*}
\f_{11} \sim \left[
\begin{array}{cc}
\f'' & 0 \\
\star & \star 
\end{array}
\right]
\end{align*}
with $\f''$ an $r\times s$-matrix. But
\begin{align*}
r \ge \text{rank} (\a \f_{11}) = \text{rank} (\f'\b) =s.
\end{align*}
This contradicts the assumption on $\f_{11}$ and finishes the
proof of the claim.

\noindent \\
{\bf (6.11) Claim:} \emph{Let $\F$ be a sheaf on $\P^2$ with
$P_{\F}(t)=nt+3$, $h^1(\F)=0$, $h^0 (\F(-1))\neq 0$.
Assume that $8 \le n \le 15$. Then $\F$ is semistable if and only
if it has a resolution}
\begin{align*}
0 \lra \O(-2) \lra (n-2)\O(-2)\oplus 3\O(-1)
\stackrel{\f}{\lra} (n-3)\O(-1) \oplus 3\O \lra \F \lra 0
\end{align*}
\emph{satisfying the properties from (6.7) and, in addition,
the following property: $\f_{11}$ is not equivalent to a matrix
of the form}
\begin{align*}
\left[
\begin{array}{cc}
\f' & 0 \\
\star & \star
\end{array}
\right] \qquad \text{\emph{where}} \quad
\f' : (m+1)\O(-2) \lra m\O(-1)
\end{align*}
\emph{is a morphism having kernel $\O(-3)$ and $m$ is an integer
satisfying $m > \frac{n}{3}+1$.}

\noindent \\
\emph{Proof:} One direction is clear: if $\F$ is semistable then
it has a resolution as in (6.7). If $\f$ did not satisfy the
``additional property'' then $\F$ would surject onto a sheaf
$\E$ with resolution
\begin{align*}
0 \lra \O(-3) \lra (m+1)\O(-2) \lra m\O(-1) \lra \E \lra 0.
\end{align*}
But
\begin{align*}
P_{\E}(t)=m {t+1 \choose 2} - (m+1) {t \choose 2} +
{t-1 \choose 2}=(m-1)t+1,
\end{align*}
which shows that $\E$ violates the semistability of $\F$
precisely when $\frac{1}{m-1} < \frac{3}{n}$, that is
$\frac{n}{3}+1 < m$.

Conversely, we assume that $\F$ has the resolution from the claim
and let $\E$ be a sheaf as in (6.9). Our aim is to arrive at
a contradiction. Since $n \le 15$ we must have $\a_0 (\E)=0$
or $\a_0 (\E)=1$. In the first case the argument is the same
as at (6.10). 
Assume now that $\a_0(\E)=1$ and let us write $P_{\E}(t)=mt+1$.
We have
\begin{align*}
h^0 (\E)=1,\quad h^1(\E)=0,\quad h^0(\E(1))=m+1, \quad
h^1(\E(1))=0
\end{align*}
so the Beilinson sequence (4.1) of $\E(1)$ gives the resolution 
\begin{align*}
0 \lra \O(-2) \lra (p+m+2)\O(-1) \stackrel{\r}{\lra}
p\O(-1) \oplus (m+1)\O \lra \E(1) \lra 0.
\end{align*}
Here $p$ is some integer and from the fact that $\E$ is
supported on a curve we get rank$(\r_{11})=p$.
So far we obtain a resolution
\begin{align*}
0 \lra \O(-3) \lra (m+2) \O(-2) \stackrel{\f'}{\lra} (m+1)\O(-1)
\lra \E \lra 0
\end{align*}
which fits into an exact commutative diagram
\begin{displaymath}
\xymatrix
{
0 \ar[r] & \O(-d) \ar[r] & (n-2)\O(-2) \ar[r]^{\f_{11}}
\ar[d]^{\b} & (n-3) \O(-1) \ar[r] \ar[d]^{\a} & \G \ar[r] \ar[d]
& 0 \\
0 \ar[r] & \O(-3) \ar[r] & (m+2)\O(-2) \ar[r]^{\f'} \ar[d] &
(m+1)\O(-1) \ar[r] \ar[d] & \E \ar[r] \ar[d] & 0 \\
& & (m+2-s)\O(-2) \ar[r]^{\f''} \ar[d] & (m+1-r) \O(-1) \ar[r]
\ar[d] & 0 \\
& & 0 & 0
} .
\end{displaymath}
Here $r,\ s$ are the ranks of $\a,\ \b$. Since $\f''$ is
surjective we must have either $m+2-s > m+1-r > 0$, i.e.
$m+1 > r > s-1$, or $r=m+1$. In the first case
\begin{align*}
\f_{11} \sim \left[
\begin{array}{cc}
\psi & 0 \\
\star & \star 
\end{array}
\right]
\end{align*}
with $\psi$ an $r\times s$-matrix. Since at least one of the
maximal minors of $\f_{11}$ is nonzero we must have $r=s$.
But then our assumption on $\f_{11}$ is contradicted.

Assume now that $r=m+1$, i.e. that $\a$ is surjective.
If $\b$ is not surjective we get the same contradiction
as above. Finally, if $\b$ is surjective then
\begin{align*}
\f_{11} \sim \left[
\begin{array}{cc}
\f' & 0 \\
\star & \star 
\end{array}
\right] .
\end{align*}
Also, $1= \a_0 (\E) < \frac{3}{n} \a_1 (\E) = \frac{3m}{n}$
forces $\frac{n}{3}+1< m+1$, so our assumption on $\f_{11}$
is contradicted. Q.e.d.

\noindent \\
{\bf (6.12) Claim:} \emph{Let $W$ be the space of morphisms}
\begin{align*}
\f: (n-2)\O(-2) \oplus 3\O(-1) \lra (n-3)\O(-1) \oplus 3\O.
\end{align*}
\emph{Let $\L=(\l_1,\l_2,\m_1,\m_2)$ be a polarization satisfying}
\begin{align*}
\l_1 < \m_1 < \frac{n-2}{n-3} \, \l_1, \qquad
\frac{3(n-2)\l_1 -1}{2(n-3)} < \m_1.
\end{align*}
\emph{Equivalently, $\L$ is such that the pair $(\l_1,\m_1)$
is in the interior of the triangle with vertices
$(0,0),\ (\frac{1}{n}, \frac{1}{n}),\ (\frac{1}{n-2}, \frac{1}{n-3})$.}
\emph{Then $\f$ is semistable with respect to $\L$ if and only
if $\f$ is not equivalent to a matrix having one of the
following forms:}
\begin{align*}
\left[
\begin{array}{ccccc}
\star & \cdots & \star & 0 & 0 \\
\vdots & & \vdots & \vdots & \vdots \\
\star & \cdots & \star & 0 & 0 \\
\star & \cdots & \star & \star & \star
\end{array}
\right], \quad \left[
\begin{array}{cccccc}
\star & \cdots & \star & 0 & 0 & 0\\
\vdots & & \vdots & \vdots & \vdots & \vdots \\
\star & \cdots & \star & 0 & 0 & 0 \\
\star & \cdots & \star & \star & \star & \star \\
\star & \cdots & \star & \star & \star & \star
\end{array}
\right], \quad \left[
\begin{array}{cccccc}
\star & \cdots & \star & \star & \star & \star \\
\vdots & & \vdots & \vdots & \vdots & \vdots \\
\star & \cdots & \star & \star & \star & \star \\
0 & \cdots & 0 & \star & \star & \star \\
0 & \cdots & 0 & \star & \star & \star \\
0 & \cdots & 0 & \star & \star & \star
\end{array}
\right],
\end{align*}
\begin{align*}
\left[
\begin{array}{cccc}
\star & \cdots & \star & \star \\
\vdots & & \vdots & \vdots \\
\star & \cdots & \star & \star \\
0 & \cdots & 0 & \star
\end{array}
\right], \qquad \left[
\begin{array}{ccccc}
\star & \cdots & \star & \star & \star \\
\vdots & & \vdots & \vdots & \vdots \\
\star & \cdots & \star & \star & \star \\
0 & \cdots & 0 & \star & \star \\
0 & \cdots & 0 & \star & \star
\end{array}
\right], \qquad \left[
\begin{array}{cccccc}
\star & \cdots & \star & 0 & \cdots & 0 \\
\vdots & & \vdots & \vdots & & \vdots \\
\star & \cdots & \star & 0 & \cdots & 0 \\
\star & \cdots & \star & \star & \cdots & \star \\
\vdots & & \vdots & \vdots & & \vdots \\
\star & \cdots & \star & \star & \cdots & \star
\end{array}
\right] .
\end{align*}
\emph{Here the zero submatrix in the last matrix has $m$
rows and $n+1-m$ columns, $1\le m \le n-3$.}

\noindent \\
\emph{Proof:} Using (3.3) we translate the conditions that
$\f$ be not equivalent to the above matrices into conditions
on $\L$. Then we express those conditions on $\L$ which
allow $\f$ to be equivalent to any matrix having a zero
submatrix which does not appear in the statement of the claim.
We arrive at 18 inequalities which, after simplifications,
reduce to the inequalities from the claim.

\noindent \\
{\bf (6.13) Claim:} \emph{Let $\L$ be a polarization as at
(6.12). If $4 \le n \le 15$ then the morphisms $\f$
from (6.10) and (6.11) form a thin constructible
subset of $W^{ss}(G,\L)$.}

\noindent \\
\emph{Proof:} First we notice that the morphisms $\f$ from
(6.10) and (6.11) are in the closed subset of $W^{ss}(G,\L)$
given by the conditions $\f_{12}=0$ and det$(\f_{22})=0$.
The condition
\begin{align*}
\f_{22} \sim
\left[
\begin{array}{rrr}
-Y & X & 0 \\
-Z & 0 & X \\
0 & -Z & Y
\end{array}
\right]
\end{align*}
is a locally closed condition because any orbit with respect
to the action of an algebraic group is a locally closed set.
The condition ${\mathcal Ker}(\f)\isom \O(-2)$ gives a
constructible set as can be seen from the sequel.
The condition
\begin{align*}
\f_{11} \nsim \left[
\begin{array}{cc}
\f' & 0 \\
\star & \star 
\end{array}
\right] \qquad \text{with} \quad {\mathcal Ker}(\f')
\isom \O(-3)
\end{align*}
gives a constructible set. To see this we only need to prove
that the condition ${\mathcal Ker}(\f')\isom \O(-3)$
gives a constructible subset inside the set of
$m \times (m+1)$-matrices with entries linear forms.
This follows from the following observation: let $G$
be an algebraic group acting on a variety $X$. Let
$Y \subset X$ be a constructible subset. Then $G.Y$
is constructible, too. To see this apply Chevalley's theorem,
stating that the image of a constructible set under
an algebraic map is again constructible, to the multiplication
map $G \times X \lra X$.

To finish the argument we need to show that the
condition ${\mathcal Ker}(\f')\isom \O(-3)$ gives a
constructible set. We represent $\f'$ by a
$m \times (m+1)$-matrix $\a = (\a_{ij})$ with entries in $V^*$.
Using the notations from (6.6) we see that
${\mathcal Ker}(\f ')\isom \O(-3)$ if and only if
deg(g.c.d.$(\a_1,\ldots,\a_{m+1})$) $=m-1$.
This, furthermore, is equivalent to the following two
conditions:
\begin{enumerate}
\item[(i)] at least two among $\a_1, \ldots, \a_{m+1}$
are linearly independent;
\item[(ii)] the system $\a_i f_j = \a_j f_i$,
$ \ 1 \le i < j \le m+1$, has a nontrivial solution
$f=(f_1,\ldots, f_{m+1})$, $\ f_j \in V^*$.
\end{enumerate}

In view of (6.6) condition (i) is equivalent to saying that
${\mathcal Ker}(\a)$ is not isomorphic to $\O(-2)$.
This is equivalent to saying that
\begin{align*}
\a \nsim \left[
\begin{array}{cc}
& 0 \\
\a' & \vdots \\
& 0
\end{array}
\right] \qquad \text{with det}(\a') \neq 0.
\end{align*}
This condition gives a constructible set because the
matrices on the right-hand side form a locally closed
subset and the smallest invariant subset containing
a locally closed subset must be constructible, as observed
above.

Condition (ii) is a closed condition. Indeed, the above
system can be written as a linear system with unknowns
the coefficients of $f_j,\ 1 \le j \le m+1$, and coefficients
the coefficients of $\a_i,\ 1 \le i \le m+1$.
Such a linear system has a nontrivial solution if
and only if the associated matrix has vanishing maximal
minors. These minors are polynomials in the coefficients
of $\a_i,\ 1 \le i \le m+1$, so we get closed conditions
on $\a$.


\section{Applications to Moduli Spaces}

Thus far, for certain classes of semistable sheaves $\F$,
we have found presentations
\bdm
\E_1 \stackrel{\f}{\lra} \E_2 \lra \F \lra 0
\edm
with decomposable vector bundles $\E_1$ and $\E_2$.
In this section we will describe some locally closed subsets
inside the corresponding moduli spaces M$_{\P^2}(r,\chi)$,
defined by means of cohomological conditions as in remark
(2.13). The question we will try to answer is whether
such a subset is a good quotient of the set of morphisms
$\f$ modulo the action by conjugation of
Aut$(\E_1)\times$Aut$(\E_2)$. The difficulty here is that
Aut$(\E)$ is a nonreductive group if $\E$ has a direct
summand of the form $\O(a) \oplus \O(b)$ with $a \neq b$.

Whenever we are dealing with a fine moduli space
we can show the existence of quotients by using the universal
family to construct local sections, cf. the proof of (7.6).
If the moduli space is not fine we need to have a quotient
already constructed as, say, in the main theorem from
\cite{drezet-2000}. We apply this theorem at (7.12) to
describe open dense subsets of M$_{\P^2}(6,4)$
and M$_{\P^2}(8,6)$. At (7.13) we construct the
quotient ad hoc as a fiber bundle over a projective
variety.

Unfortunately, Dr\'ezet and Trautmann's
theory of quotients modulo nonreductive groups is still
incomplete. For instance, the main theorem
from \cite{drezet-2000} does not cover the quotients
from (4.7). Thus, we are not able to describe as a quotient
an open dense subset of M$_{\P^2}(9,6)$. Also, we do not
know if quotients exist for morphisms of type (2,2).
This accounts for the ``unknowns'' in the table from
the introduction.

We begin by recalling the notions of good and geometric
quotients. 
Let $G$ be a linear algebraic group acting on a
variety $X$. The action is algebraic, i.e. the map $G\times X
\lra X$ given by $(g,x) \lra g.x$ is a morphism of varieties.

\noindent \\
{\bf (7.1) Definition:} A \emph{categorical quotient} of $X$ by $G$
is a pair $(Y,\p)$ where $Y$ is a variety and $\p :X\lra Y$ is
a $G$-equivariant morphism satisfying the following universal
property: for any other $G$-equivariant morphism $\eta: X \lra Z$
there exists a unique morphism $\r :Y\lra Z$ making the diagram
commute:
\begin{displaymath}
\xymatrix
{
X \ar[d]_{\p} \ar[rd]^{\eta} \\
Y \ar[r]^{\r} & Z
}.
\end{displaymath}
We write $Y=X//G$. If, in addition, the fibers of $\p$ are orbits,
then $Y$ is called an \emph{orbit space} and is denoted $X/G$.
Note that a categorical quotient is unique up to isomorphism.

\noindent \\
{\bf (7.2) Definition:} A \emph{good quotient} of $X$ by $G$ is a
pair $(Y,\p)$ where $Y$ is a variety and $\p: X\lra Y$ is a morphism
satisfying:
\begin{enumerate}
\item[(i)] $\p$ is $G$-equivariant;
\item[(ii)] $\p$ is surjective;
\item[(iii)] for any open subset $U\subset Y$ the pull-back map
$\p^*$ gives an isomorphism of $\O_Y (U)$ onto the ring
of regular functions on $\p^{-1}(U)$ which are constant
on the $G$-orbits;
\item[(iv)] if $W\subset X$ is closed and $G$-invariant, then $\p (W)$
is closed;
\item[(v)] if $W_1, W_2 \subset X$ are closed, $G$-invariant
and disjoint, then $\p (W_1)$ and $\p (W_2)$ are also disjoint;
\item[(vi)] $\p$ is affine, i.e. it returns open affine sets to
affine sets.
\end{enumerate}
If, in addition, the fibers of $\p$ are orbits, then $(Y,\p)$ is called
a \emph{geometric quotient}.

\noindent \\
{\bf (7.3) Remark:} Let $G$ act on $X$ as above and assume
the existence of an affine surjective morphism $\p :X\lra Y$
whose fibers are orbits. Assume that $\p$ admits local sections, i.e.
for any $y\in Y$ there is an open neighbourhood $U$ of $y$ and a
morphism $\s: U\lra X$ satisfying $\p \circ \s = 1$.
Then $(Y,\p)$ is a geometric quotient. \\

Definition (7.2) is important because good quotients are
categorical quotients while geometric quotients are orbit spaces:

\noindent \\
{\bf (7.4) Proposition:} \emph{Let $(Y,\p)$ be a good
quotient of $X$ by $G$. Then:
\begin{enumerate}
\item[(i)] $(Y,\p)$ is a categorical quotient;
\item[(ii)] $\p (x_1)=\p (x_2)$ if and only if $\overline{G.x_1}$
intersects $\overline{G.x_2}$;
\item[(iii)] if the $G$-orbits in $X$ are closed, then $(Y,\p)$ is an
orbit space;
\item[(iv)] Let $X_o$ denote the subset of points $x\in X$ with
$G.x$ closed and of maximal dimension among the $G$-orbits.
Then there is an open subset $Y_o\subset Y$ such that $\p^{-1}(Y_o)
=X_o$ and $(Y_o,\p)$ is a geometric quotient of $X_o$ by $G$.
\end{enumerate}}

The main technical tool that we will use in this section is the
relative Beilinson complex. Given a variety $X$ and a coherent
sheaf $\F$ on $X \times \P^2$ there is a sequence
\bdm
0 \lra \C^{-2} \lra \C^{-1} \lra \C^0 \lra \C^1 \lra \C^2 \lra 0
\edm
of sheaves on $X \times \P^2$ which is exact, except in the
middle, where the cohomology is $\F$. On each fiber
$\{ x \} \times \P^2$ this sequence restricts to the Beilinson
complex of the restricted sheaf $\F_x$. Let $p:X\times \P^2
\lra X$ be the projection onto the first component.
The sheaves $\C^i$ are defined by means of the higher direct
images of $\F$:
\bdm
\C^i = \oplus_j R^j p_* (\F \tensor \Om^{j-i}_{X \times \P^2/X}
(j-i)) \boxtimes \O_{\P^2}(i-j).
\edm
In our applications $\F$ will be flat over $X$ and its restrictions
$\F_x$ onto the fibers $\{ x \} \times \P^2$ will have one-dimensional
supports. Thus $H^2(\F_x)=0$ for all $x \in X$. From the Base
Change Theorem on p. 11 in \cite{okonek} we get $R^2 p_* (\F)=0$.
Analogously, the other second direct images occuring above are zero.
The relative Beilinson complex now takes the form
\begin{align*}
\tag{7.5}
0 \lra \C^{-2} \lra \C^{-1} \lra \C^0 \lra \C^1 \lra 0
\end{align*}
with
\begin{eqnarray*}
\C^{-2} & = & p_* (\F(-1)) \boxtimes \O(-2), \\
\C^{-1} & = & p_* (\F \tensor \Om^1_{X \times \P^2/X}(1))
\boxtimes \O(-1) \oplus R^1 p_* (\F(-1)) \boxtimes \O(-2), \\
\C^0 & = & p_* (\F) \boxtimes \O \oplus
R^1 p_* (\F \tensor \Om^1_{X \times \P^2/X}(1)) \boxtimes \O(-1),\\
\C^1 & = & R^1 p_* (\F) \boxtimes \O.
\end{eqnarray*}

\noindent \\
{\bf (7.6) Proposition:} \emph{Let $n \ge 1$ be an integer and
let $W$ be the vector space of morphisms $\f$ of sheaves on
$\P^2$ of the form}
\bdm
\O(-2) \oplus (n-1)\O(-1) \stackrel{\f}{\lra} n\O.
\edm
\emph{With the notations from section 3 assume that the
polarization $\L = (\l_1, \l_2, \m_1)$ satisfies $0 < \l_1 <
\frac{1}{n}$. Let $W_o$ be the open subset of $W^{ss}(G,\L)$
given by the condition det$(\f) \neq 0$. Then $W_o$ admits
a geometric quotient modulo $G$ which is isomorphic to the
open dense subset of \emph{M}$_{\P^2}(n+1,n)$ given by the condition
$h^0 (\F(-1))\neq 0$.}

\noindent \\
\emph{Proof:} Let us consider the coherent sheaf $\tilda{\F}$
on $W_o \times \P^2$ given by the exact sequence
\bdm
\O_{W_o} \boxtimes \O_{\P^2}(-2) \oplus \O_{W_o} \boxtimes
(n-1) \O_{\P^2}(-1) \stackrel{\Phi}{\lra} n \O_{W_o \times \P^2}
\lra \tilda{\F} \lra 0.
\edm
On each fiber $\{ \f \} \times \P^2$ this sequence restricts to
\bdm
0 \lra \O(-2) \oplus (n-1) \O(-1) \stackrel{\f}{\lra} n\O
\lra \tilda{\F}_{\f} \lra 0.
\edm
According to (4.2) each restriction $\tilda{\F}_{\f}$ is
semistable with Hilbert polynomial $P(t)=(n+1)t+n$.
As the Hilbert polynomial is independent of $\f$, and as the
base $W_o$ is reduced, the sheaf $\tilda{\F}$ is flat over
$W_o$. By the definition (2.5) of a coarse moduli space,
$\tilda{\F}$ gives rise to a morphism
\bdm
\eta : W_o \lra \text{M}_{\P^2}(n+1,n)
\edm
which sends $\f$ to the stable equivalence class of
$\tilda{\F}_{\f}$.

By (4.2) the image of $\eta$ is the subset M$_o$ of
M$_{\P^2}(n+1,n)$ given by the condition $h^0(\F(-1))=0$.
By (2.13) this subset is open and, as the moduli space
is irreducible, it must be dense.

The fibers of $\eta$ are $G$-orbits. Indeed, an isomorphism
$f$ between two cokernels $\F_1$ and $\F_2$ of $\f_1$ and $\f_2$
from $W_o$ must fit into a commutative diagram
\bdm
\xymatrix
{
0 \ar[r] & \O(-2) \oplus (n-1)\O(-1)
\ar[r]^{\quad \quad \quad \f_1} \ar[d]^g
& n\O \ar[r] \ar[d]^h & \F_1 \ar[r] \ar[d]^f & 0 \\
0 \ar[r] & \O(-2) \oplus (n-1) \O(-1)
\ar[r]^{\quad \quad \quad \f_2}
& n\O \ar[r] & \F_2 \ar[r] & 0
}
\edm
in which $g$ and $h$ are isomorphisms. Here $h$ is defined
in such a way as to coincide with $f$ on the level of global
sections, while $g$ is the induced map on the kernels.

To prove that $\eta$ is a geometric quotient map it is enough
to construct local sections as in (7.3). For this we will use
the fact that M$_{\P^2}(n+1,n)$ is a fine moduli space,
cf. (2.10), so it has a universal family. Let $\U$ denote the
restriction of the universal family to M$_o \times \P^2$.
Let $p :$ M$_o \times \P^2 \lra $ M$_o$ be the projection
onto the first component. $\U$ is flat over M$_o$ and all its
restrictions to the fibers of $p$ have Beilinson resolution
(4.2). In view of the Base Change Theorem the higher direct
images
\bdm
p_* (\U(-1)), \quad \quad
R^1 p_* (\U(-1)), \quad \quad
p_* (\U), \quad \quad
R^1 p_* (\U),
\edm
\bdm
p_* (\U \tensor \Om^1_{\text{M}_o \times \P^2/\text{M}_o}(1)),
\quad
R^1 p_* (\U \tensor \Om^1_{\text{M}_o \times \P^2/\text{M}_o}(1))
\edm
are locally free of ranks 0, 1, $n$, 0, $n-1$, 0. Let us cover
M$_o$ with open subsets $S$ on which the above sheaves are free.
On $S \times \P^2$, and relative to fixed trivializations
of the higher direct images, the Beilinson complex (7.5)
gives the resolution
\bdm
0 \lra \O_{S} \boxtimes \O_{\P^2}(-2) \oplus \O_{S} \boxtimes
(n-1) \O_{\P^2}(-1) \stackrel{\f}{\lra} n \O_{S \times \P^2}
\lra \U \lra 0.
\edm
We put $\s (x)=\f_x$ for $x \in S$ and notice that $\s : S \lra
W_o$ is a local section of $\eta$. Q.e.d. \\

The sets $W_o$ are nonempty for all $n$. Indeed, it is easy to
construct an $n \times (n-1)$-matrix $\psi$ with entries in
$V^*$ whose maximal minors are linearly independent, and
which has the form
\bdm
\left[
\ba{cccc}
\star & \star & \cdots & \star \\
\star & \star & \cdots & \star \\
0 & \star & \cdots & \star \\
\vdots & \vdots & \ddots & \vdots \\
0 & 0 & \cdots & \star
\ea
\right] . \qquad \text{For example} \quad \left[
\ba{cccc}
Y & Z & Y & Z \\
X & 0 & 0 & 0 \\
0 & Y & 0 & 0 \\
0 & 0 & Z & 0 \\
0 & 0 & 0 & X 
\ea
\right]
\edm
is such a matrix for $n=5$. But it now becomes clear that
the following matrix is semistable and has nonzero determinant:
\bdm
\left[
\ba{cc}
X^2 & 0 \\
0 & \psi
\ea
\right].
\edm

The existence of the geometric quotient $W_o/G$ can be put into
a broader context if we realize that $W^{ss}(G,\L)$ itself has
a geometric quotient, as J.-M. Dr\'ezet pointed out to the author:

\noindent \\
{\bf (7.7) Proposition:} \emph{Let $W^{ss}(G,\L)$ be the set
of morphisms of sheaves on $\P^2$,}
\bdm
\O(-2) \oplus (n-1) \O(-1) \stackrel{\f}{\lra} n\O,
\edm
\emph{which are semistable with respect to a polarization $\L$
satisfying $0 < \l_1 < \frac{1}{n}$. Then there exists a
geometric quotient $W^{ss}(G,\L)/G$ which is a fiber bundle
with fiber $\P^{3n+2}$ and base a projective variety of
dimension $n^2-n$.}

\noindent \\
\emph{Proof:} Let us represent $\f$ as a pair $(\f_1,\f_2)$
where $\f_1$ is an $n \times 1$-matrix with entries in $S^2 V^*$,
while $\f_2$ is an $n \times (n-1)$-matrix with entries in $V^*$.
Let $W_i$ denote the vector space of matrices $\f_i$.

The reductive group $G_2=$GL$(n-1)\times$GL$(n)/k^*$ acts
on $W_2$ by conjugation. Here $k^*$ is embedded as the subgroup
of homotheties. The only admissible polarization on $W_2$
is $\left( \frac{1}{n-1},\frac{1}{n} \right)$ and, as $n-1$
and $n$ are mutually prime, equality cannot be achieved in
(3.3). This shows that the set of semistable points
$W_2^{ss}$ for the action of $G_2$ coincides with the set of
stable points. By the classical Geometric Invariant Theory
there is a geometric quotient $W_2^{ss}/G_2$ which is
a projective variety of dimension $n^2-n$. Let us denote it
by $N$.

We view $W$ as a trivial bundle with fiber $W_1$ and base $W_2$.
Let $U$ be the trivial bundle on $W_2$ with fiber the space
of $(n-1) \times 1$-matrices with entries in $V^*$.
We consider the morphism of bundles $f:U \lra W$ given at
every point $\f_2$ by left-multiplication with $\f_2$.
It is easy to see that $f$ is injective at every semistable point
$\f_2$, hence the restriction of ${\mathcal Coker}(f)$ to
$W_2^{ss}$ is a vector bundle of rank $3n+3$, denoted by $E$.
$\P(E)$ carries a $G_2$ action which is compatible with the
action on $W_2^{ss}$.
At (8.1) below we will prove that
for any $\f_2 \in W^{ss}_2$ the isotropy
group Stab$_{G_2}(\f_2)$ is trivial, so it acts trivially
on $\P(E_{\f_2})$.
It follows that $\P(E)$ descends to a fiber bundle $F$
on $N$, see 4.2.15 in \cite{hl}.

We notice that the semistability conditions on $\f$ read
as follows: $\f_2$ is in $W_2^{ss}$ and $\f_1'\neq 0$
for all $\f'$ in the same $G$-orbit as $\f$.
In other words, $W^{ss}(G,\L)$ can be identified with
the complement of ${\mathcal Im}(f)$ inside
${W_1}_{|W_2^{ss}}$. The map
\bdm
{W_1}_{|W_2^{ss}} \setminus {\mathcal Im}(f) \lra \P(E)
\edm
admits local sections because $E$ is a bundle. In view
of (7.3) this map is a geometric quotient modulo the
action of the subgroup of $G$ given by the conditions
$h_1=1$, $g_2=1$, see the notations preceeding (3.1).
Combining this with the fact that the map
$\P(E) \lra F$ is a geometric quotient modulo $G_2$,
we easily deduce that the map $W^{ss}(G,\L) \lra F$
is a geometric quotient modulo $G$. Q.e.d. \\

Our construction is similar
to, though much less elegant than, the construction
from 10.2 in \cite{drezet-trautmann} which addresses
morphisms on $\P^n$ of the form
\bdm
m_1 \O(-2) \oplus m_2 \O(-1) \stackrel{\f}{\lra} n_1 \O.
\edm
The polarization satisfies $0 < \l_1 < \l_{\text{min}}$,
where $\l_{\text{min}}$ is the smallest positive number
such that for $\l_1$ varying in the interval
$(0, \l_{\text{min}})$ the set of semistable points
remains unchanged. In the context of the above
proposition $\l_{\text{min}}=\frac{1}{n}$.
They show that if certain conditions on the integers
$m_1, m_2, n_1$ are satisfied, then there exists
a geometric quotient which is a Grassmann bundle 
Grass$(m_1, p_* \E (2))$ with base $N=W_2^{ss}/G_2$.
Here $p : N \times \P^2 \lra N$ is the projection onto the
first component and 
$\E$ is the universal sheaf on $N \times \P^2$
which restricts to ${\mathcal Coker}(\f_2)$
on each fiber $\{ [\f_2 ] \} \times \P^2$.
One of their conditions, having to do
with the injectivity of $\f_2$ regarded as map from
$m_2 \O(-1)$ to $n_1 \O$, is $n_1 \ge n m_2$.
Thus Dr\'ezet and Trautmann's construction addresses only
the case $n=2$ of the above proposition. \\

There is yet another, more direct way of constructing
the quotient in the case $n=2$. The semistability
conditions on a morphism $\f : \O(-2) \oplus \O(-1)
\lra 2\O$ read: det$(\f)\neq 0$ and $\f_{12}, \f_{22}$
are linearly independent in $V^*$. The map
\bdm
W^{ss}(G,\L) \lra \text{Grass}(2,V^*) \times \P(S^3 V^*)
\isom \P^2 \times \P(S^3 V^*)
\edm
given by
\bdm
\f \lra (\text{span}(\f_{12},\f_{22}), <\text{det}(f)>)
\edm
has fibers $G$-orbits and has image the universal cubic
\bdm
C = \{ (x,<f>) \in \P^2 \times \P(S^3 V^*), \quad f(x)=0 \}.
\edm
It was first noticed in \cite{maican} that the map
$W^{ss}(G,\L) \lra C$ has local sections:
Choose a point $(x,<f>)$ in $C$, say $x=(0:0:1)$.
As $f$ does not contain the monomial $Z^3$,
there are unique quadratic polynomials $q_1 (X,Y,Z)$
and $q_2 (X,Z)$ such that $f=q_1 Y-q_2 X$.
We put
\bdm
\s (x,<f>)= \left[
\begin{array}{cc}
q_1 & X \\
q_2 & Y
\end{array}
\right].
\edm
Since all processes involved in defining $\s$ are
algebraic, we see that $\s$ extends to a section
of the map $W^{ss}(G,\L)\lra C$ defined on a neighbourhood
of $(x,<f>)$. Thus $W^{ss}(G,\L)/G \isom C$.
A more sophisticated proof of this isomorphism can be found in
\cite{fr-diplom}.\\

In the simplest case $n=1$, $W^{ss}(G,\L)$ is just the set
of nonzero morphisms $\O(-2) \lra \O$ and
$W^{ss}(G,\L)/G$ is $\P(S^2 V^*)$. \\

As noticed, in the cases $n=1,2$ we have $W_o = W^{ss}(G,\L)$,
hence $W_o/G$ is complete, hence the set M$_o$ from (7.6)
is complete, hence M$_o$ is the entire moduli space.
We have obtained the well-known fact that every semistable
sheaf on $\P^2$ with Hilbert polynomial $P(t)=2t+1$ is
the structure sheaf of a conic; in other words
M$_{\P^2}(2,1) \isom \P(S^2 V^*)$.
In the case $n=3$ we have rediscovered one of Le Potier's
result from \cite{lepotier} to the effect that
M$_{\P^2}(3,2)$ is isomorphic to the universal cubic.

If $n \ge 3$ $W_o$ is a proper open subset of the set of
semistable points, hence $W_o/G$ is not complete, hence
M$_o$ is a proper open subset of M$_{\P^2}(n+1,n)$.
Indeed, it is easy to construct semistable morphisms
with zero determinant; for example, in the case $n=3$,
\bdm
\left[
\ba{ccc}
0 & X & Y \\
XY & Z & 0 \\
-X^2 & 0 & Z
\ea
\right].
\edm
Thus, at most we can say at this time is the following:

\noindent \\
{\bf (7.8) Corollary:} \emph{For $n \ge 3$ the projective
varieties $W^{ss}(G,\L)/G$ and \emph{M}$_{\P^2}(n+1,n)$ are
birational.}

\noindent \\
\emph{Proof:} From (7.2)(iv) we see that the image of
$W_o$ under the quotient map
\bdm
W^{ss}(G,\L) \lra W^{ss}(G,\L)/G
\edm
is an open set $U$. In fact, $W_o$ is the
preimage of $U$. Clearly, the properties from (7.2) are
satisfied for the map $W_o \lra U$. This proves that
$W_o/G \isom U$ and so we have isomorphic open dense
subsets of $W^{ss}(G,\L)/G$ and of M$_{\P^2}(n+1,n)$. \\

The same proof as at (7.6) can be used to show that
for all fine moduli spaces M$_{\P^2}(r,\chi)$ occuring
in sections 4, 5, 6 the locally closed subsets described
by cohomological conditions are geometric quotients
$W_o/G$ of the corresponding sets $W_o \subset W^{ss}(G,\L)$.
We have summarized the results in the table
from the introduction. For the quotients in section 5
we should mention that $\f_x$ depends in an algebraic
manner on the maps from the Beilinson complex of $\U_x$,
hence it depends in an algebraic manner on $x$;
see (7.14) for the details.

The assumption that M$_{\P^2}(r,\chi)$ be fine,
i.e. the assumption that a universal family exists,
is needed for the construction
of the local sections of $\eta$. The proof of (7.6) does
not apply if the moduli space is not fine because,
according to (2.11), there is no universal family on
any open subset of such a moduli space.

Two of the quotients from sections 5 and 6 have very
concrete descriptions. First we consider the case
$n=1$ from (5.2). The set of morphisms
\bdm
\f :\O(-3) \oplus \O(-1) \lra 2\O
\edm
semistable with respect to a polarization satisfying
$0 < \l_1 < \frac{1}{2}$ is characterized by the conditions
det$(\f)\neq 0$ and $\f_{12}$, $\f_{22}$ are linearly
independent in $V^*$. The same discussion as in the case
$n=2$ of (7.6) shows that the geometric quotient
$W^{ss}(G,\L)/G$ is isomorphic to the universal
quartic in $\P^2 \times \P(S^4 V^*)$. From (5.2) we get:

\noindent \\
{\bf (7.9) Corollary:} \emph{The subset of
\emph{M}$_{\P^2}(4,1)$ given by the conditions
$h^0 (\F(-1))=0$ and $h^1(\F)=1$ is closed and,
equipped with its canonical reduced structure, it is
isomorphic to the universal quartic in
$\P^2 \times \P(S^4 V^*)$.} \\

Let us now consider the simplest case $n=4$ from (6.7).
It concerns morphisms
\bdm
\f : 2\O(-2) \oplus 3\O(-1) \lra \O(-1) \oplus 3\O
\edm
satisfying the conditions: $\f_{12}=0$, $\f_{11}$
has linearly independent entries in $V^*$, $\f_{21}' \neq 0$
for any $\f'$ in the same orbit as $f$, $\f_{22}$
is equivalent to the matrix
\bdm
\left[
\begin{array}{rrr}
-Y & X & 0 \\
-Z & 0 & X \\
0 & -Z & Y
\end{array}
\right] .
\edm
Let $f=0$ be the equation of the support of $\F$.
To be precise,
\bdm
f = \left[
\ba{rrr}
Z & -Y & X
\ea
\right] \f_{21} \left[
\ba{r}
-X_2 \\ X_1
\ea
\right], \quad \text{where} \quad \f_{11} = \left[
\ba{cc}
X_1 & X_2
\ea
\right]
\edm
and $\f_{22}$ is assumed to be the above $3 \times 3$-matrix.
We consider the $G$-invariant map
\bdm
W_o \lra \text{Grass}(2,V^*) \times \P(S^4 V^*) \isom
\P^2 \times \P(S^4 V^*)
\edm
given by
\bdm
\f \lra (\text{span}\{ X_1,X_2 \}, <f>).
\edm
Its image is the universal quartic. To prove that the
map $W_o \lra Q$ is a geometric quotient, we will
construct local sections. We fix a point
$(\text{span}\{ X_1,X_2 \}, <f>)$ in $Q$. We complete
$\{ X_1, X_2 \}$ to a basis $\{ X_1, X_2, X_3 \}$ of
$V^*$. Relative to this basis $f$ can be uniquely written
as
\bdm
f(X_1,X_2,X_3) = -X_2 f_1(X_1,X_2,X_3)+ X_1 f_2 (X_1,X_3).
\edm
Now $f_1$ and $f_2$ can each be uniquely written as
\begin{eqnarray*}
f_1 & = & Z q_{11}(X,Y,Z)-Y q_{21}(X,Y)+ X q_{31}(X), \\
f_2 & = & Z q_{12}(X,Y,Z)-Y q_{22}(X,Y)+ X q_{32}(X).
\end{eqnarray*}
We put
\bdm
\s (\text{span}\{ X_1,X_2 \},<f>)= \left[
\ba{llrrr}
X_1 & X_2 & 0 & 0 & 0 \\
q_{11} & q_{12} & -Y & X & 0 \\
q_{21} & q_{22} & -Z & 0 & X \\
q_{31} & q_{32} & 0 & -Z & Y
\ea
\right] .
\edm
Since all processes involved in defining $\s$ are algebraic,
we see that $\s$ extends to a local section defined on an
open subset of $Q$. From (6.7) we get:

\noindent \\
{\bf (7.10) Corollary:} \emph{The subset of
\emph{M}$_{\P^2}(4,3)$ given by the conditions $h^0(\F(-1))=1$
and $h^1(\F)=0$ is closed and, equipped with its canonical
reduced structure, is isomorphic to the universal quartic in
$\P^2 \times \P(S^4 V^*)$.} \\

We now turn to the moduli spaces M$_{\P^2}(r,\chi)$ for
which $r$ and $\chi$ are not mutually prime. As we shall
see, if we knew the existence of the quotient $W_o//G$,
then we could prove that this quotient is isomorphic to
the corresponding subvariety of the moduli space.
We know the existence of the quotients only in two
cases: for the situation in (4.3)(i) and for $n=3$
in (5.2). In the first case we will use a theorem 
of Dr\'ezet:

Let $m_1, m_2, n_1$ be integers and let us consider
morphisms of sheaves on $\P^n$ of the form
\bdm
m_1 \O(-2) \oplus m_2 \O(-1) \stackrel{\f}{\lra}
n_1 \O.
\edm
We recall from section 3 that a polarization in this
context is a triple $\L = (\l_1,\l_2,\m_1)$ of positive
numbers satisfying the relations $m_1\l_1 + m_2 \l_2
= n_1 \m_1 = 1$. Theorem 6.4 from \cite{drezet-2000}
gives sufficient conditions on $\L$ which assure the
existence of a good quotient.
Below we state part two of the theorem formulated in
the particular case $n=2$ which is of interest to us:

\noindent \\
There exists a good quotient $W^{ss}(G,\L)//G$, which
is a projective variety, if the following four inequalities
are fulfilled:
\begin{eqnarray*}
\l_2 & < & \frac{3}{n_1}, \\
\l_2 & > & \frac{3m_1 + n_1}{3m_1 n_1 + n_1 m_2}, \\
m_2 \l_2 & > & 1 - \frac{3m_1}{n_1 (3m_1-1)} \quad \text{if}
\quad m_1 \le 3,\\
m_2 \l_2 & > & 1 - \frac{3m_1}{8 n_1} \quad \text{if} \quad
m_1 > 3.
\end{eqnarray*}

\noindent \\
Taking $m_1=2$, $m_2=n-2$, $n_1=n$ the above conditions
become
\bdm
\l_1 < \frac{6}{n(n+4)}, \qquad \l_1 < \frac{3}{5n}.
\edm

\noindent \\
{\bf (7.11) Corollary:} \emph{Let $3 \le n \le 7$ be an
integer and let $W^{ss}(G,\L)$ be the space of morphisms
of sheaves on $\P^2$ of the form}
\bdm
2 \O(-2) \oplus (n-2)\O(-1) \stackrel{\f}{\lra} n\O
\edm
\emph{which are semistable with respect to a polarization
$\L$ satisfying $\frac{1}{2n} < \l_1 < \frac{1}{n}$.
Then there exists a good quotient $W^{ss}(G,\L)//G$,
which, moreover, is a projective variety.}

\noindent \\
\emph{Proof:} If $\frac{1}{2n} < \l_1 < \frac{6}{n(n+4)}$,
the statement follows from Dr\'ezet's theorem.
To conclude the proof we only need observe that $W^{ss}(G,\L)$
does not change when $\l_1$ varies in the interval
$\left( \frac{1}{2n}, \frac{1}{n} \right)$.

\noindent \\
{\bf (7.12) Proposition:} \emph{For $3 \le n \le 6$ let
$W_o$ be the subset of $W^{ss}(G,\L)$ from (7.11) given by
the condition det$(\f)\neq 0$. For $n=2$ let $W_o$ be the
space of injective morphisms $2\O(-2) \lra 2\O$. Then
$W_o$ admits a good quotient modulo $G$, which is isomorphic
to the open dense subset of \emph{M}$_{\P^2}(n+2,n)$
given by the conditions}
\bdm
h^0 (\F(-1))=0, \quad h^1(\F)=0,
\quad h^1(\F \tensor \Om^1(1))=0.
\edm
\emph{In particular, the projective varieties
$W^{ss}(G,\L)//G$ and \emph{M}$_{\P^2}(n+2,n)$ are birational.}

\noindent \\
\emph{Proof:} The good quotient $W_o//G$ is an open dense
subset of $W^{ss}(G,\L)//G$. The latter exists by (7.11)
when $n \ge 3$ and by the classical Geometric Invariant
Theory when $n=2$. The map
\bdm
\eta: W_o \lra \text{M}_{\P^2}(n+2,n)
\edm
can be constructed as at (7.6) and has image the open subset
M$_o$ described by the cohomological conditions from the
proposition. By the universal property (7.4)(i) of a good
quotient, $\eta$ factors through a morphism
\bdm
\r : W_o //G \lra \text{M}_o.
\edm
If $n$ is even the injectivity of $\r$ is not as
straightforward as at (7.6) because the fibers of $\eta$
may not be $G$-orbits, as there may occur properly
semistable sheaves. We will prove the injectivity only in
the case $n=2$, the cases $n=4$ and $n=6$ being analogous:

Let $[\f_1]$ and $[\f_2]$ be in $W_o//G$ and assume that
$\F_1 = {\mathcal Coker}(\f_1)$ and $\F_2 = {\mathcal Coker}
(\f_2)$ are properly semistable and stable equivalent.
Thus $\F_1$ and $\F_2$ have the same terms in their
Jordan-H\"older filtrations, say $\A_1$ and $\A_2$.
According to the discussion preceeding (7.8),
$\A_i$ are cokernels of maps $\a_i : \O(-2) \lra \O$.
It is easy to see that, modulo the action of $G$,
$\f_1$ and $\f_2$ are equivalent to matrices
\bdm
\psi_1 = \left[
\ba{cc}
\a_1 & \b_1  \\
0 & \a_2
\ea
\right], \quad \text{respectively} \quad \psi_2 = \left[
\ba{cc}
\a_1 & \b_2 \\
0 & \a_2 
\ea
\right].
\edm
We consider the one-parameter subgroup $\l$ of $G$ given by
\bdm
\l(t) = \left( \left[
\ba{cc}
t & 0 \\
0 & 1
\ea \right], \left[
\ba{cc}
t & 0 \\
0 & 1
\ea \right] \right).
\edm
We have
\bdm
\l(t).\psi_i = \left[
\ba{cc}
\a_1 & t \b_i \\
0 & \a_2 
\ea
\right] \quad \text{forcing} \quad
\lim_{t \to 0} \l(t). \psi_i = \left[
\ba{cc}
\a_1 & 0 \\
0 & \a_2
\ea
\right]
\edm
which we denote by $\psi$. From (7.4)(ii) we get
$[\psi_i]=[\psi]$, so $[\f_i]=[\psi]$, so $\rho$ is
injective.

To finish the proof we only need observe that M$_o$
is smooth. At points represented by stable sheaves this is
already known from (2.12). In general, applying the long
exact sequence of Ext groups to the exact sequence
(4.3)(i), we deduce that for all $\F$ in M$_o$ we have
Ext$^2(\F,\F)=0$. According to Grothendieck's Criterion,
this gives smoothness at the point
in the moduli space represented by $\F$.

Thus far $\r$ is a bijective morphism onto a normal
variety. From Zariski's Main Theorem we conclude that
$\r$ is an isomorphism. \\

We notice that another way of proving that $\r$ is
an isomorphism, which avoids Grothendieck's Criterion
of smoothness and Zariski's Main Theorem,
is exhibited in the proof of (7.14). \\

We do not know if the birational maps constructed above are
isomorphisms. The subsets M$_o$ are open, proper subsets of
M$_{\P^2}(n+2,n)$ because $W_o$ are proper subsets of
the sets of semistable points. For example, in the case
$n=2$, the following matrices are semistable but have
zero determinant:
\bdm
\left[
\ba{cc}
X_1 Y_1 & X_1 Y_2 \\
X_2 Y_1 & X_2 Y_2
\ea
\right]
\edm
where $X_1, X_2$ are linearly independent in $V^*$
and same for $Y_1, Y_2$.

\noindent \\
{\bf (7.13) Proposition:} \emph{For $n=1, 2, 3$ let
$W^{ss}(G,\L)$ be the set of morphisms of sheaves on $\P^2$
of the form}
\bdm
\O(-3) \oplus n\O(-1) \stackrel{\f}{\lra} (n+1)\O,
\edm
\emph{which are semistable with respect to a polarization
$\L$ satisfying $0 < \l_1 < \frac{1}{n+1}$. Then there
exists a geometric quotient $W^{ss}(G,\L)/G$ which
is a fiber bundle with fiber $\P^{4n+9}$ and base a
projective variety of dimension $n^2+n$.}

\noindent \\
The proof is the same as at (7.7). The injectivity of $f$
is clear in the cases $n=1,2$ and follows from remark
(5.4) in the case $n=3$.

\noindent \\
{\bf (7.14) Proposition:} \emph{Let $W_o$ be the open subset
of $W^{ss}(G,\L)$ from (7.13) given by the condition
det$(\f)\neq 0$. Then $W_o$ admits a geometric quotient
modulo $G$ which is isomorphic to the locally closed
subset of \emph{M}$_{\P^2}(n+3,n)$ given by the conditions
$h^0 (\F(-1))=0$ and $h^1(\F)=1$, and equipped with its
canonical reduced structure.}

\noindent \\
\emph{Proof:} The cases of the fine moduli spaces
M$_{\P^2}(4,1)$ and M$_{\P^2}(5,2)$ were discussed
earlier. Assume now that $n=3$. Let $X$ be the subset
of M$_{\P^2}(6,3)$ described by the cohomological
conditions from the proposition.

As at (7.6), there is a morphism $\eta : W_o \lra X$
associated to a flat family on $W_o$
and which factors through a morphism $\r : W_o/G \lra X$.
From (5.3) we
know that all sheaves from $X$ are stable, so we can
repeat the argument from (7.6) proving that the
fibers of $\eta$ are $G$-orbits. Thus $\r$ is bijective.

To prove that $\r$ is an isomorphism we will construct
its inverse.
Let us recall from section 2
that M$_{\P^2}(6,3)$ is the good quotient of
a certain open subset $R$
inside a quotient scheme, modulo the action of
SL$(V)$.
There is a locally closed subvariety
$S$ of $R$, invariant under the action
of the special linear group, such that $X=S//$SL$(V)$.
The existence of $S$ follows from remark 3.4.3 on
p. 54 in \cite{newstead} and from the fact that
in characteristic zero reductive groups are linearly
reductive, cf. p. 50 in loc. cit.
In fact $S$ is the preimage of $X$
under the quotient map $R \lra $ M$_{\P^2}(6,3)$.
Let $\t : S \lra X$ denote the quotient map.

Let $\U$ be the restriction
to $S \times \P^2$ of the universal quotient family
on $R \times \P^2$.
Let $p : S \times \P^2 \lra S$ be the projection
onto the first component.
For an arbitrary point in $s \in S$ we denote by
$\U_s$ the restriction $\U_{|\{ s \} \times \P^2}$.
From (5.2) we know that $\U_s$ has a resolution
\bdm
0 \lra \O(-3) \oplus 3\O(-1) \stackrel{\f}{\lra}
4\O \lra \U_s \lra 0.
\edm
In fact, with the notations from section 5, we have
\bdm
\f_{12}= \r_{32}, \qquad 
\f_{11}= - \r_{31} \psi_{12}^{-1} \left[
\ba{c}
X \\ Y \\ Z
\ea
\right].
\edm
Each $\U_s$ is the middle cohomology of a Beilinson
complex 
\bdm
0 \lra 3\O(-2) \oplus 3\O(-1) \lra 3\O(-1) \oplus 4\O
\lra \O \lra 0
\edm
and $\r_{32}$, $\r_{31}$, $\psi_{12}$ depend algebraically
on the maps in this complex.
We put $\varsigma (s)=\f$ and we claim 
that $\varsigma$ can be extended to a
morphism from a neighbourhood $S_o$ of $s$ in $S$ to $W_o$.

To see this we proceed as in the proof of (7.6).
The higher direct image sheaves
\bdm
p_* (\U(-1)), \quad \quad
R^1 p_* (\U(-1)), \quad \quad 
p_* (\U), \quad \quad
R^1 p_* (\U),
\edm
\bdm
p_* (\U \tensor \Om^1_{S \times \P^2/ S}(1)),
\quad
R^1 p_* (\U \tensor \Om^1_{S \times \P^2/ S}(1))
\edm
are locally free of ranks 0, 3, 4, 1, 3, 3.
They are free on an open neighbourhood $S_o$ of $s$.
Thus $\r_{32}$, $\r_{31}$, $\psi_{12}$ can be made to
depend algebraically on the point in $S_o$.
This allows us to define $\varsigma$ on $S_o$.

We now cover $S$ with such open sets $S_o$ and
we notice that the locally defined maps $\p \circ
\varsigma$ glue together to a globally defined
morphism $\s : S \lra W_o/G$ making the diagram
commute:
\bdm
\xymatrix
{
W_o \ar[d]_{\p} & S_o \ar[l]_{\varsigma} \ar[r]^i
& S \ar[d]^{\t} \ar[dll]_{\s} \\
W_o/G \ar[rr]^{\r} & & X
}.
\edm
Indeed, if $\varsigma_1$ and $\varsigma_2$ are defined on
two distinct neighbourhoods of $s$, then
$\varsigma_1 (s)$ and $\varsigma_2 (s)$ are in the same
$G$-orbit.

Finally, let us observe that $\s$ is constant on the
fibers of $\t$. This is so because if $\t (s_1)
= \t (s_2)$, then the corresponding sheaves $\U_{s_1}$
and $\U_{s_2}$ are isomorphic, so their Beilinson resolutions
are equivalent, i.e. $\varsigma_1 (s)$
and $\varsigma_2 (s)$ are in the same $G$-orbit.
Here $\varsigma_i$ is defined on a neighbourhood
of $s_i$.

By the universal property (7.4)(i) of a good quotient,
the map $\s$ factors through a morphism from $X$
to $W_o/G$. This is the desired inverse of $\r$. Q.e.d. \\

The above proof could be carried out for all locally closed
subsets $X \subset$ M$_{\P^2}(r,\chi)$ occuring at
(4.3), (4.5)(ii), (4.7), (4.8), provided that we knew
the existence of the quotients $W_o//G$. In all cases we
would get the isomorphism $X \isom W_o//G$.
Unfortunately, we do not know how to prove the existence
of $W_o//G$ when $r$ and $\chi$ are not mutually prime
in each of the above cases.
We should mention that an essential ingredient in the
proof is the fact that all sheaves from $S$ have the same
kind of Beilinson complex. This is satisfied because
the cohomological conditions
defining $X$ are closed under stable equivalence.
This fact is easy to check in each case. To give the
flavor of the argument we will just check the case
$n=6$ from (4.3)(i): assume that $\G$ is stable equivalent
to $\F$ and that $\F$ has resolution
\bdm
0 \lra 2\O(-2) \oplus 4\O(-1) \lra 6\O \lra \F \lra 0.
\edm
Assume that $\F$ is properly semistable, so it fits
into an exact sequence
\bdm
0 \lra \F_1 \lra \F \lra \F_2 \lra 0
\edm
with $\F_1$ and $\F_2$ in M$_{\P^2}(4,3)$.
From $h^1 (\F)=0$ and $h^2 (\F_1)=0$ we get
$h^1 (\F_2)=0$. Analogously, from $h^1 (\F \tensor \Om^1(1))
=0$ and from $h^2 (\F_1 \tensor \Om^1(1))=0$ we get
$h^1 (\F_2 \tensor \Om^1(1))=0$. We cannot have
$h^0 (\F_2 (-1)) >0$ because, in view of (6.7), this would
force $h^1 (\F_2 \tensor \Om^1(1))=1$.
Thus $h^0 (\F_2(-1))=0$. From $h^0 (\F(-1))=0$ we
immediately also get $h^0(\F_1(-1))=0$.
In conclusion, both $\F_1$ and $\F_2$ satisfy the
hypotheses of (4.2). We arrive at the resolutions
\bdm
0 \lra \O(-2) \oplus 2 \O(-1) \lra 3\O \lra \F_i \lra 0.
\edm
By hypothesis $\G$ is an extension
\bdm
0 \lra \F_1 \lra \G \lra \F_2 \lra 0
\edm
possibly with $\F_1$ and $\F_2$ interchanged.
By the ``horseshoe lemma'' the resolutions of $\F_1$
and $\F_2$ can be combined to give a resolution
for $\G$ of the same kind as the resolution of $\F$.


\section{Computation of Codimensions}

To find the codimensions of the locally closed subvarieties
of M$_{\P^2}(r,\chi)$ occuring in the previous sections we need
to find the dimensions of the stabilizers of generic points
from $W_o$.
For actions of reductive groups it is known that
a stable point has zero-dimensional isotropy group.
This fact will not remain true in our context.

We begin with a lemma which seems to be known, yet we couldn't
find a reference. Let $V$ be a vector space over $k$ and let
$W$ be the space of $m \times n$-matrices with entries in $V$.
We consider the action by conjugation on $W$ of the reductive
group $G=$GL$(m)\times$GL$(n)/k^*$.

\noindent \\
{\bf (8.1) Lemma:} \emph{The isotropy subgroup of a stable
point from $W$ is trivial.}

\noindent \\
\emph{Proof:} Let $w \in W$ be a stable matrix. Concretely,
what this means, is that no matrix in the same orbit as $w$
can have a zero $p \times q$-submatrix with
$\frac{p}{m} + \frac{q}{n} \ge 1$. We consider an element
in the isotropy group of $w$ represented by $(g,h)$.

As $G$ is reductive, Stab$_G(w)$ is finite, so there are
$t \in k^*$ and an integer $r \ge 1$ such that
$g^r = t I_m$ and $h^r = t I_n$. From this we see that
$g$ and $h$ are diagonalizable matrices. Replacing possibly
$w$ by another point in its orbit, we may assume that
$g$ and $h$ are diagonal matrices. Let us write
\bdm
g=\text{diag}(t_1,\ldots,t_m),
\qquad
h=\text{diag}(s_1,\ldots,s_n).
\edm
From $w=gwh^{-1}$ we see that $w_{ij}=0$ if $t_i \neq s_j$.
Thus, if $t_1, \ldots, t_m, s_1, \ldots s_n$ are not all
equal, then $w$ is a block matrix, say
\bdm
\left[
\ba{cc}
\star & 0 \\
0 & \star
\ea
\right] .
\edm
This contradicts the stability of $w$. In conclusion
$g= t I_m$, $h=t I_n$, i.e. $(g,h)$ represents the identity
of $G$.

\noindent \\
{\bf (8.2) Claim:} \emph{The isotropy group of a generic
semistable morphism}
\bdm
2\O(-2) \oplus (n-1)\O(-1) \stackrel{\f}{\lra} \O(-1) \oplus
n\O, \qquad \f_{12}=0,
\edm
\emph{has dimension $n-1$. The semistability conditions
are understood to be as at (4.3).}

\noindent \\
\emph{Proof:} We choose a morphism $\f$ for which at least
one of the maximal minors of $\f_{22}$ is nonzero.
Let $(g,h)$ be in Stab$_G (\f)$. Keeping the notations
from section 3 we write
\bdm
g^{-1}= \left[
\ba{cc}
g_1 & 0 \\
u & g_2
\ea
\right], \qquad \qquad h=\left[
\ba{cc}
h_1 & 0 \\
v & h_2
\ea
\right].
\edm
We have $\f = h\f g^{-1}$ so $\f_{11}=h_1 \f_{11} g_1$
and $\f_{22}= h_2 \f_{22} g_2$. But $\f_{11}$ and $\f_{22}$
are stable matrices with entries in $V^*$.
From (8.1) we get $h_1=t_1$, $g_1=t_1^{-1}$, $h_2=t_2 I_n$,
$g_2= t_2^{-1}I_{n-1}$. We have
$\f_{21}=v \f_{11}t_1^{-1}+ t_2 \f_{21} t_1^{-1}
+ t_2 \f_{22} u.$
If $t_1 \neq t_2$, then $\f$ is equivalent to a matrix $\f'$
for which $\f_{21}' = 0$. This would contradict the
semistability of $\f$. Thus $t_1=t_2=t$ and
$\f_{22} u=-t^{-2} v \f_{11}$.

Recall that $\f_{11}=[X_1,X_2]$ with linearly independent
$X_1,X_2$ in $V^*$. We put $\psi=[-X_2, X_1]^T$.
From $\f_{22} u \psi =-t^{-2} v \f_{11} \psi =0$
we get $u \psi =0$ because one of the maximal minors of
$\f_{22}$ is nonzero. Thus $u=\a \f_{11}$ with
$\a \in \ $ M$_{n-1,1}(k)$. From
$(t^{-1}v+t\f_{22}\a)\f_{11}=0$ we get 
$v= -t^2 \f_{22} \a$. Thus $(g,h)$ is parametrized by $t$
and by the entries of $\a$, giving the claim. \\

The above proof worked because for $\f_{11}$ there existed
a matrix $\psi$ such that for any $1 \times 2$-matrix $u$
with entries in $V^*$
\begin{align*}
\tag{*}
u \psi = 0 \quad \text{implies that $u$ is a multiple
of $\f_{11}$}.
\end{align*}
For morphisms from (4.7) we can take
\bdm
\f_{11}= \left[
\ba{ccc}
X & Y & Z
\ea 
\right], \qquad \qquad \psi= \left[
\ba{rrr}
-Y & -Z & 0 \\
X & 0 & -Y \\
0 & X & Z
\ea
\right]
\edm
and we see that (*) is true for $1 \times 3$-matrices $u$
with entries in $V^*$. We arrive at:

\noindent \\
{\bf (8.3) Claim:} \emph{The isotropy group of a generic
semistable morphism}
\bdm
3\O(-2) \oplus (n-2)\O(-1) \stackrel{\f}{\lra} \O(-1) \oplus
n\O, \qquad \f_{12}=0,
\edm
\emph{has dimension $n-2$. The semistability conditions
are understood to be as at (4.7).} \\

For morphisms from (4.8) and (4.9) the $2 \times 3$-matrix
$\f_{11}=(f_{ij})_{i=1,2,j=1,2,3}$ with entries in $V^*$
is stable.
Concretely, stability here means that the maximal
minors of $\f_{11}$ are linearly independent in $S^2 V^*$.
We put $f=[f_1,f_2,f_3]^T$, where
\bdm
f_1 = \left|
\ba{cc}
f_{12} & f_{13} \\
f_{22} & f_{23}
\ea
\right| , \qquad 
f_2 = \left|
\ba{cc}
f_{13} & f_{11} \\
f_{23} & f_{21}
\ea
\right| , \qquad 
f_3 = \left|
\ba{cc}
f_{11} & f_{12} \\
f_{21} & f_{22}
\ea
\right| .
\edm
Clearly $\f_{11} f=0$. Our intention is to show that,
for generic $\f_{11}$, and for a $1\times 3$-matrix $u$
with entries in $V^*$,
the equality $uf = 0$ implies that $u$ is a linear combination
of the rows of $\f_{11}$. Indeed, the condition $uf=0$
is the same as saying that the determinant of
\bdm
\psi = \left[
\ba{ccc}
f_{11} & f_{12} & f_{13} \\
f_{21} & f_{22} & f_{23} \\
u_1 & u_2 & u_3
\ea
\right]
\edm
is zero. We need to prove that, modulo operations on
rows and columns, $\psi$ is equivalent to a matrix
having a zero row. For this we will use (5.5),
namely we will exclude the other possibilities listed
there. First we see that, as the columns of $\f_{11}$
are linearly independent, $\psi$ cannot be equivalent
to a matrix having a zero column.
Nor is $\psi$ equivalent to a matrix of the form
\bdm
\left[
\ba{ccc}
0 & 0 & \star \\
0 & 0 & \star \\
\star & \star & \star
\ea
\right] ,
\edm
for if $g \psi h$ has the above form, then all
$2\times 2$-minors positioned on the first two columns
of $\psi h$ are zero. But the matrix obtained by deleting
the third row of $\psi h$ is equivalent to $\f_{11}$,
so it is stable, so its first maximal minor from the left
is nonzero (in fact all its maximal minors are nonzero).

If we choose $\f_{11}$ generic enough, then $\psi$
is equivalent to neither $\psi_1$ nor $\psi_2$ from (5.5).
For instance, if $\f_{11}$, regarded as a map from
$V^* \oplus V^* \oplus V^*$ to $V^* \oplus V^*$, is injective,
then $\psi$ is not equivalent to $\psi_1$.
To rule out $\psi_2$, we need only observe that the
condition det$(\psi_2)\ =0$ defines a thin subset inside
the affine space with coordinates $a_1, \ldots a_5$
(notations as at (5.5)). In conclusion, $\psi$ is equivalent
to a matrix having a zero row.

With the notations from the proof of (8.2), we have
$u=\a \f_{11}$ with $\a \in \ $ M$_{n-1,2}(k)$, and
$v = -t^2 \f_{22} \a$. We arrive at the following:

\noindent \\
{\bf (8.4) Claim:} \emph{The isotropy group of a generic
semistable morphism}
\bdm
3\O(-2) \oplus (n-1)\O(-1) \stackrel{\f}{\lra} 2\O(-1) \oplus
n\O, \qquad \f_{12}=0,
\edm
\emph{has dimension $2n-2$.} \\

Finally, we turn to morphisms from (6.10) and (6.11).

\noindent \\
{\bf (8.5) Claim:} \emph{The isotropy group of a generic semistable
morphism}
\bdm
(n-2)\O(-2) \oplus 3\O(-1) \stackrel{\f}{\lra} (n-3)\O(-1)
\oplus 3\O, \qquad \f_{12}=0, \quad \f_{22}=\psi_1,
\edm
\emph{has dimension $4n-11$.}

\noindent \\
\emph{Proof:} As at (8.2) we have $t_1 = t_2 = t$ and
$v \f_{11} t^{-1}= - t \f_{22} u$. We put $\psi = [Z,-X,Y]$.
From $\psi v \f_{11} = -t^2 \psi \f_{22} u = 0$ we get
$\psi v=0$, because $\f_{11}$ can be chosen generic enough that
one of its maximal minors be nonzero. From $\psi v=0$ we get
$v = \f_{22} \a$ with $\a \in \ $ M$_{3,n-3}$. From
$\f_{22}(\a \f_{11} t^{-1} + t u)=0$ we get
$\a \f_{11} t^{-1} + t u =  [-Y,X,Z]^T \b$.
with $\b \in$ M$_{1,n-2}(k)$. Thus Stab$_G (\f)$ is parametrized
by t, the entries of $\a$ and the entries of $\b$. Q.e.d. \\

Once we know the dimensions of the isotropy groups of
generic points $\f \in W_o$ we can
apply the obvious formula
\bdm
\text{dim}(X)= \text{dim}(W_o)- \text{dim}(G) +
\text{dim}(\text{Stab}_G (\f)).
\edm
We do not carry out here these computations; we refer, insted,
to the table from the introduction where we have recorded the
results.


\section{Duality Results}

In \cite{fr-diplom} one can find a birational map of
fine moduli spaces
\bdm
\text{M}_{\P^2}(r,\chi) \lra \text{M}_{\P^2}(r,r-\chi)
\edm
given by sending a point represented by $\F$ to the point
represented by the dual sheaf $\F^D$. By modifying
slightly the argument from \cite{fr-diplom} we will construct
such birational maps also for those coarse moduli spaces
occuring in section 4. At (9.5) we will also obtain isomorphisms
between dual locally closed subspaces of M$_{\P^2}(r,\chi)$
and M$_{\P^2}(r,r-\chi)$.

In the sequel $\F$ will be a coherent sheaf on $\P^2$
with pure one-dimensiomal support and without zero-dimensional
torsion. We define its \emph{dual} $\F^D$ by
\bdm
\F^D = {\mathcal Ext}^1_{\O_{\P^2}}(\F,\Om^2_{\P^2})(1).
\edm
Clearly $\F^D$ has one-dimensional support, so it has linear
Hilbert polynomial. This can be computed using the following
isomorphisms provided by Serre Duality:
\bdm
H^0 (\F^D (-1)) \isom H^1 (\F)^*, \qquad
H^1 (\F^D (-1)) \isom H^0 (\F)^*,
\edm
\bdm
H^0 (\F^D) \isom H^1 (\F(-1))^*, \qquad
H^1 (\F^D) \isom H^0 (\F(-1))^*.
\edm
Thus, if $P_{\F}(t)=rt+\chi$, then $P_{\F^D}(t)=rt+r-\chi$.
In particular, the slopes of $\F$ and $\F^D$ are related by
$p(\F^D)= 1 - p(\F)$.

\noindent \\
{\bf (9.1) Lemma:} \emph{If $\F$ is Cohen-Macaulay, in particular
if $\F$ is semistable, then $\F^{DD} \isom \F$ and
${\mathcal Ext}^2(\F,\Om^2)=0$.}

\noindent \\
\emph{Proof:} We will apply proposition 1.1.10 from
\cite{hl}. All we need to show is that $\F$ satisfies the
Serre condition S$_{2,1}$:
\bdm
\text{depth}(\F_x) \ge \text{min}\{ 2,\ \text{dim}\,
O_{\P^2, x} -1 \} \quad \text{for all} \quad x \in
\text{Supp}(\F).
\edm
But if $x$ is a closed point in the support of $\F$,
we have depth$(\F_x)=1$ and dim$\, \O_{\P^2, x}=2$.
If $x$ is a generic point of an irreducible component of
Supp$(\F)$, we have dim$\, \O_{\P^2, x}=1$ and the above
inequality is trivially fulfilled.

Finally, we notice that, by virtue of (2.3), semistable sheaves
are Cohen-Macaulay.

\noindent \\
{\bf (9.2) Lemma:} \emph{$\F$ is (semi)stable if and only if
$\F^D$ is (semi)stable.}

\noindent \\
\emph{Proof:} Assume that $\F$ is semistable. Let $\G= \F^D/\K$
be a quotient sheaf of $\F^D$. As $\K$ is a torsion sheaf,
we have ${\mathcal Hom}(\K,\Om^2)=0$. Applying the long
exact sequence in ${\mathcal Ext}$-sheaves to the short
exact sequence
\bdm
0 \lra \K \lra \F^D \lra \G \lra 0
\edm
we see that $\G^D$ is a subsheaf of $\F^{DD} \isom \F$.
Thus
\bdm
1- p(\G)= p(\G^D) \le p(\F) = 1 - p(\F^D), \quad
\text{so} \quad p(\F^D) \le p(\G).
\edm
This proves the semistability of $\F^D$.

Assume that $\F$ is not semistable. Then $\F$ has a quotient
sheaf $\G$ with $p(\G) < p(\F)$. As before, $\G^D$ is a
destabilizing subsheaf of $\F^D$.

\noindent \\
{\bf (9.3) Lemma:} \emph{If $\F$ and $\G$ are semistable
and stable equivalent, then so are $\F^D$ and $\G^D$.}

\noindent \\
\emph{Proof:} Let us consider a Jordan-H\"older filtration
for $\F$:
\bdm
0 = \F_0 \subset \F_1 \subset \ldots \subset \F_{n-1}
\subset \F_n = \F.
\edm
We apply the long exact sequence in ${\mathcal Ext}$-sheaves
to the exact sequences
\bdm
0 \lra \F_i \lra \F_{i+1} \lra \F_{i+1}/\F_i \lra 0.
\edm
As $\F_i$ is a torsion sheaf we have
${\mathcal Hom}(\F_i,\Om^2)=0$. As $\F_{i+1}/\F_i$ is semistable,
we have, by (9.1), ${\mathcal Ext}^2(\F_{i+1}/\F_i, \Om^2)=0$.
We arrive at the exact sequences
\bdm
0 \lra (\F_{i+1}/\F_i)^D \lra \F_{i+1}^D \lra \F_i^D \lra 0.
\edm
Similarly we obtain exact sequences
\bdm
0 \lra (\F/\F_i)^D \lra \F^D \lra \F_i^D \lra 0.
\edm
From these two sequences we conclude that
\bdm
0 = (\F/\F_n)^D \subset (\F/\F_{n-1})^D \subset \ldots \subset
(\F/\F_1)^D \subset (\F/\F_0)^D = \F^D
\edm
is a Jordan-H\"older filtration of $\F^D$ with terms
$(\F_{i+1}/\F_i)^D$, the latter being stable by virtue of (9.2).
The lemma follows.

\noindent \\
{\bf (9.4) Theorem:} \emph{Assume that $\frac{r}{2} \le \chi
\le r$ and that $r$, $\chi$ are mutually prime. Then the
open dense subset of \emph{M}$_{\P^2}(r,\chi)$ given by the
conditions}
\bdm
h^0(\F(-1))=0, \qquad h^1(\F)=0,
\qquad h^1(\F \tensor \Om^1(1))=0,
\edm
\emph{is isomorphic to the 
open dense subset of \emph{M}$_{\P^2}(r,r-\chi)$ given by the
conditions}
\bdm
h^1(\F)=0, \qquad h^0(\F(-1))=0,
\qquad h^0(\F \tensor \Om^1(1))=0.
\edm
\emph{The isomorphism is given by $[\F] \lra [\F^D]$.}

\noindent \\
\emph{Proof:} From (9.1), (9.2) and (9.3) we see that the map
$\d$ given by $[\F] \lra [\F^D]$ is well defined and a 
bijection between the two open sets from the theorem,
which we call M$_o(r,\chi)$ and M$_o(r,r-\chi)$.

Every sheaf $\F$ from M$_o(r,\chi)$ has Beilinson
resolution
\begin{align*}
\tag{*}
0 \lra (r-\chi)\O(-2) \oplus (2\chi -r)\O(-1)
\stackrel{\f}{\lra} \chi \O \lra \F \lra 0.
\end{align*}
The long exact sequence in ${\mathcal Ext}$-sheaves gives the
resolution
\bdm
0 \lra \chi \O(-2) \stackrel{\f^D}{\lra} (r-\chi)\O \oplus
(2\chi -r)\O(-1) \lra \F^D \lra 0,
\edm
where $\f^D$, viewed as a matrix, is simply the transpose
of $\f$. The set of morphisms $\f$ occuring above forms
an open subset $W_o$ inside the vector space of morphisms
\bdm
(r-\chi)\O(-2) \oplus (2\chi -r)\O(-1) \lra \chi \O.
\edm
On $W_o \times \P^2$ we consider the coherent sheaf
$\tilda{\F}$ given by the exact sequence
\bdm
\O_{W_o} \boxtimes (r-\chi)\O_{\P^2}(-2) \oplus
\O_{W_o} \boxtimes (2\chi-r)\O_{\P^2}(-1)
\stackrel{\Phi}{\lra} \chi \O_{W_o \times \P^2}
\lra \tilda{\F} \lra 0.
\edm
On each fiber $\{ \f \}\times \P^2$ the restriction of
$\Phi$ is $\f$. Similarly we construct the dual family as the
cokernel
\bdm
\O_{W_o} \boxtimes \chi \O_{\P^2}(-2)
\stackrel{\Phi^D}{\lra}
\O_{W_o} \boxtimes (r-\chi)\O_{\P^2} \oplus
\O_{W_o} \boxtimes (2\chi -r) \O_{\P^2}(-1)
\lra \tilda{\F}^D \lra 0
\edm
of a morphism $\Phi^D$ which restricts to $\f^D$ on each fiber
$\{ \f \} \times \P^2$. Clearly $\tilda{\F}$ and
$\tilda{\F}^D$ are $W_o$-flat, so they induce morphisms
\bdm
\r : W_o \lra \text{M}_o (r,\chi), \qquad
\r^D : W_o \lra \text{M}_o (r,r-\chi).
\edm
We have $\d \circ \r = \r^D$.

Next we recall from section 2 that M$_{\P^2}(r,\chi)$ is the
good quotient of an open subset $R$ inside a certain
quotient scheme. Let $S$ be the preimage of M$_o(r,\chi)$
under the quotient map $R \lra $ M$_{\P^2}(r,\chi)$.
The map $\p : S \lra $ M$_o(r,\chi)$ ia a good quotient map.
Let $\U$ be the restriction to $S \times \P^2$ of the
universal quotient family on $R \times \P^2$.
From the fact that all restrictions of $\U$ to the fibers
$\{ s \} \times \P^2$, $s \in S$, have Beilinson resolution
(*) we deduce, as in the proof of (7.14), the existence of
locally defined morphisms $\varsigma : S_o \lra W_o$
satisfying $\r \circ \varsigma = \p$. The morphisms
$\r^D \circ \varsigma$ glue to a globally defined morphism
$\p^D$ making the diagram commute:
\bdm
\xymatrix
{
 & S \ar[dddl]_{\p} \ar[dddr]^{\p^D} \\
 & S_o \ar[u]_i \ar[d]^{\varsigma} \\
 & W_o \ar[dl]_{\r} \ar[dr]^{\r^D} \\
\text{M}_o (r,\chi) \ar[rr]^{\d} & & \text{M}_o (r,r-\chi)
}.
\edm
Thus $\d$ is the map induced by $\p^D$ via the universal
property of the quotient map $\p$. As such, $\d$ must be
a morphism. By symmetry, its inverse must be a morphism, too. \\

The above theorem first appeared in \cite{fr-diplom}.
Its proof given there is simpler and makes use of the universal
families on the fine moduli spaces. Our argument, though more
cumbersome, has the following advantage: it works also in the
case when $r$, $\chi$ are not mutually prime, as long as
we know that all sheaves giving a point in M$_o$ have the
same kind of Beilinson complex.

\noindent \\
{\bf (9.5) Theorem:} \emph{Let $X$ be the locally closed
subvariety of \emph{M}$_{\P^2}(r,\chi)$ given by the
conditions}
\bdm
h^0(\F(-1))=a, \qquad h^0(\F)=b, \qquad
h^0(\F \tensor \Om^1(1))=c.
\edm
\emph{Assume that every sheaf giving a point in $X$ satisfies
the above conditions. This is the case, for instance, when
$r$, $\chi$ are mutually prime. Then $X$ is isomorphic to the
locally closed
subvariety $X^D$ of \emph{M}$_{\P^2}(r,r-\chi)$ given by the
conditions}
\bdm
h^1(\F)=a, \qquad h^1(\F(-1))=b, \qquad
h^1(\F \tensor \Om^1(1))=c.
\edm
\emph{The isomorphism is given by $[\F] \lra [\F^D]$.
Here $X$ and $X^D$ are equipped with their canonical reduced
structures.}

\noindent \\
\emph{Proof:} Assume that $X$ is nonempty.
We repeat the arguments from (9.4).

We consider vector bundles $\E^i$ on $\P^2$, $i=-2, -1, 0, 1$,
which are decomposable as direct sums of line bundles.
We assume that for each $\F$ giving a point in $X$
there is a complex
\begin{align*}
\tag{*}
0 \lra \E^{-2} \stackrel{\f}{\lra} \E^{-1} \stackrel{\f'}{\lra}
\E^0 \stackrel{\f''}{\lra} \E^1 \lra 0
\end{align*}
which is exact, except at $\E^0$, where the cohomology is $\F$.
For instance, we could choose $\E^i$ to be the bundles $\C^i$
occuring in the Beilinson complex (4.1).
Let $W_o$ be the set of the above complexes.

$W_o$ will play the same role as in the proof of (9.4).
The existence of $\r : W_o \lra X$ is clear by construction.
To finish the proof, we only need to construct $\r^D$
satisfying $\d \circ \r = \r^D$. For this purpose we will
show that ${\mathcal Hom}(\_,\Om^2)(1)$ applied to (*)
gives a complex
\begin{align*}
\tag{**}
0 \lra \E^{1}_D \lra \E^{0}_D \lra
\E^{-1}_D \lra \E^{-2}_D \lra 0
\end{align*}
which is exact, except at $\E^{-1}_D$,
where the cohomology is $\F^D$.

We consider the long exact sequences in
${\mathcal Ext}(\_,\Om^2)$-sheaves induced by the short exact
sequences
\bdm
0 \lra \E^{-2} \lra \E^{-1} \lra \A \lra 0,
\edm
\bdm
0 \lra \B \lra \E^0 \lra \E^1 \lra 0,
\edm
\bdm
0 \lra \A \lra \B \lra \F \lra 0.
\edm
As $\Ext^j (\O(d),\_)=0$ for $j \ge 1$, we have
$\Ext^j (\E^i,\Om^2)=0$ for $j \ge 1$. The second sequence
gives $\Ext^1 (\B,\Om^2)=0$. In view of (9.1), the semistability
of $\F$ leads to $\Ext^2 (\F,\Om^2)=0$. The third sequence
gives $\Ext^1(\A,\Om^2)=0$. Thus we arrive at the exact
sequences
\bdm
0 \lra \Hom (\A,\Om^2) \lra \Hom (\E^{-1},\O^2) \lra
\Hom (\E^{-2},\Om^2) \lra 0,
\edm
\bdm
0 \lra \Hom (\E^1, \Om^2) \lra \Hom (\E^0,\Om^2) \lra
\Hom (\B,\Om^2) \lra 0,
\edm
\bdm
0 \lra \Hom (\B,\Om^2) \lra \Hom (\A, \Om^2) \lra
\Ext^1 (\F,\Om^2) \lra 0
\edm
which immediately yield (**). \\

We mentioned at the end of section 7 that all locally closed
subvarieties $X$ occuring in section 4, satisfy the hypotheses
of the above theorem. Indeed, it can be verified in each case
that the cohomological properties defining $X$ are closed
under stable equivalence. As a consequence, all locally
closed subvarieties $X \subset $ M$_{\P^2}(r,\chi)$ occuring
in sections 4, 5, 6, with the possible exception of the
subvarieties in M$_{\P^2}(3r,3)$, $r=3,4,5$, occuring in
section 6, are isomorphic to their duals $X^D$.
In particular, (9.4) remains true for the following choices
of multiplicity and Euler characteristic: (6,4), (8,6), (9,6).
We obtain the following:

\noindent \\
{\bf (9.6) Corollary:} \emph{For the following choices of
$(r,\chi)$ the spaces \emph{M}$_{\P^2}(r,\chi)$
and \emph{M}$_{\P^2}(r,r-\chi)$ are birational: (6,4), (8,6),
(9,6).}

\noindent \\
Here is another application of (9.5): the closed subset
of M$_{\P^2}(4,1)$ given by the conditions $h^0(\F(-1))=0$,
$h^1(\F)=1$ (the condition $h^0(\F \tensor \Om^1(1))=1$
is automatically fulfilled), is isomorphic to the closed
subset of M$_{\P^2}(4,3)$ given by the conditions
$h^1(\F)=0$, $h^0(\F(-1))=1$. This we proved earlier at
(7.9) and (7.10) by means of their description as
geometric quotients. \\

Let $W_o^D$ denote the set of complexes (**), i.e. the
set of complexes obtained by applying $\Hom (\_,\Om^2)(1)$
to the complexes from $W_o$. If we identify $W_o$ with
a certain subset of triples of matrices $(\f,\f',\f'')$
inside the vector space
\bdm
W= \text{Hom}(\E^{-2},\E^{-1}) \times \text{Hom} (\E^{-1},\E^0)
\times \text{Hom} (\E^0,\E^1),
\edm
then $W_o^D$ is just the subset of triples of transposed matrices
$({\f''}^{T},{\f'}^{T},\f^T)$ inside the vector space
\bdm
W^D= \text{Hom}(\E^{1}_D,\E^{0}_D) \times
\text{Hom} (\E^{0}_D,\E^{-1}_D)
\times \text{Hom} (\E^{-1}_D,\E^{-2}_D).
\edm
Thus transposition gives an isomorphism of $W_o$ with $W_o^D$,
both equipped with their canonical reduced structures
induced by the ambient spaces $W$ and $W^D$.

On $W_o$ and on $W_o^D$ we have the canonical action of the
(usually nonreductive) algebraic group
\bdm
G= \text{Aut}(\E^{-2}) \times \text{Aut} (\E^{-1}) \times
\text{Aut} (\E^0) \times \text{Aut} (\E^1).
\edm
From the proofs of (9.4) and (9.5) we extract the following:

\noindent \\
{\bf (9.7) Proposition:} \emph{Let $X$ be as in (9.5). Assume that
a good quotient of $W_o$ by $G$ exists and is isomorphic to
$X$. Then a good quotient of $W_o^D$ by $G$ exists and is isomorphic
to $X^D$.} \\

For every subset $X \subset $ M$_{\P^2}(r,\chi)$ described in
section 7 as a good (geometric) quotient, we have a dual description
of $X^D \subset $ M$_{\P^2}(r,r-\chi)$ as a good (geometric)
quotient. For better understanding let us introduce to a
polarization $\L$ of type (2,1) or (2.2) its dual polarization
$\L^D$ of type (1,2), respectively (2.2):
\begin{eqnarray*}
\text{for} \quad \L = (\l_1,\l_2,\m_1) & \text{we put} &
\L^D = (\l_1^D,\m_1^D,\m_2^D)=(\m_1,\l_2,\l_1); \\
\text{for} \quad \L = (\l_1,\l_2,\m_1,\m_2) & \text{we put} &
\L^D = (\l_1^D,\l_2^D,\m_1^D,\m_2^D)=(\m_1,\m_2,\l_2,\l_1).
\end{eqnarray*}
If $W_o$ is defined by semistability conditions
expressed in terms of $\L$,
then $W_o^D$ is defined by semistability conditions expressed in
terms of $\L^D$. We list below the consequences of (9.7)
for the cases of generic sheaves:

\noindent \\
{\bf (9.8) Corollary:} \emph{The open dense
subset of \emph{M}$_{\P^2}(n+1,1)$,
$n \ge 2$, given by the condition $h^1(\F)=0$, is isomorphic to
$W_o/G$, where $W_o$ is the set of injective morphisms}
\bdm
n\O(-2) \stackrel{\f}{\lra} (n-1)\O(-1) \oplus \O, \quad
\f \in W^{ss}(G,\L),\ \L = (\l_1,\m_1,\m_2),\ 0 < \m_2 < \frac{1}{n}.
\edm
\emph{The open dense
subset of \emph{M}$_{\P^2}(n+2,2)$, $n=3,4,5,6$, given
by the conditions $h^0 (\F(-1))=0$, $h^1 (\F)=0$, $h^0(\F \tensor
\Om^1 (1))=0$, is isomorphic to $W_o//G$, where
$W_o$ is the set of injective morphisms}
\bdm
n\O(-2) \stackrel{\f}{\lra} (n-2)\O(-1) \oplus 2\O, \quad
\f \in W^{ss}(G,\L),\ \L = (\l_1,\m_1,\m_2),\
\frac{1}{2n} < \m_2 < \frac{1}{n}.
\edm
\emph{The open dense
subset of \emph{M}$_{\P^2}(n+3,3)$, $n=4,5$, given
by the conditions $h^0 (\F(-1))=0$, $h^1 (\F)=0$, $h^0(\F \tensor
\Om^1 (1))=0$, is isomorphic to $W_o/G$, where
$W_o$ is the set of injective morphisms}
\bdm
n\O(-2) \stackrel{\f}{\lra} (n-3)\O(-1) \oplus 3\O, \quad
\f \in W^{ss}(G,\L),\ \L = (\l_1,\m_1,\m_2),\
\frac{2}{3n} < \m_2 < \frac{1}{n}.
\edm

\noindent \\
One final example of a quotient we were not able to obtain
in section 6: the subset of M$_{\P^2}(6,3)$ given by the
conditions $h^0 (\F(-1))=1$, $h^1 (\F)=0$ (the condition
$h^1 (\F \tensor \Om^1(1))=3$ is automatically fulfilled),
is isomorphic to $W_o/G$, where $W_o$ is the set of injective
morphisms
\bdm
4\O(-2) \stackrel{\f}{\lra} 3\O(-1) \oplus \O(1), \quad
\f \in W^{ss}(G,\L),\ \L= (\l_1,\m_1,\m_2),\ 0 < \m_2 < \frac{1}{4}.
\edm
Applying (9.7) to the quotients from section 6 we get descriptions
for the subsets in M$_{\P^2}(n+3,n)$, $n=4,5,7,8,10,11$,
given by the conditions $h^0 (\F(-1))=0$, $h^1 (\F)=1$.
We omit the details.


\end{document}